\documentclass[10pt]{article}
\usepackage[top=1in, bottom=1in, left=1in, right=1in]{geometry}
\usepackage{srcltx,graphicx}
\usepackage{amsmath, amssymb, amsthm}
\usepackage{color, xcolor, colortbl}
\usepackage{array, multirow}
\usepackage{subfig, overpic}
\usepackage{hyperref}
\usepackage{epstopdf}
\usepackage{algorithm}
\usepackage{algpseudocode}
\usepackage[capitalize,nameinlink]{cleveref} 
\usepackage{ulem}

\allowdisplaybreaks[4]

\graphicspath{{images/}}

\newtheorem{property}{Property}

\newcommand\bbR{\mathbb{R}}

\newcommand\bbN{\mathbb{N}}
\newcommand\bbS{\mathbb{S}}

\newcommand\cH{\mathcal{H}}
\newcommand\cI{\mathcal{I}}
\newcommand\cJ{\mathcal{J}}
\newcommand\cM{\mathcal{M}}
\newcommand\cP{\mathcal{P}}
\newcommand\cC{\mathcal{C}}
\newcommand\cT{\mathcal{T}}
\newcommand\cU{\mathcal{U}}
\newcommand\cK{\mathcal{K}}

\newcommand\ad{\mathrm{ad}}

\newcommand\Ntrainsample{{$N^{\mathrm{train}}_{\mathrm{samples}}$}}
\newcommand\Ntestsample{{$N^{\mathrm{test}}_{\mathrm{samples}}$}}
\newcommand\Nparams{{$N_{\mathrm{params}}$}}
\newcommand\trainerror{{\epsilon_{\mathrm{train}}}}
\newcommand\testerror{{\epsilon_{\mathrm{test}}}}

\newcommand\NL{\mathsf{NL}}
\newcommand\IL{\mathsf{IL}}
\newcommand\LCR{\mathsf{LCR}}
\newcommand\LCK{\mathsf{LCK}}
\newcommand\LCI{\mathsf{LCI}}
\newcommand\CNNR{\mathsf{CR}}
\newcommand\CNNK{\mathsf{CK}}
\newcommand\CNNI{\mathsf{CI}}
\newcommand\Reshape{\mathsf{Reshape}}
\newcommand\ReshapeT{\mathsf{ReshapeT}}
\newcommand\ReshapeM{\mathsf{ReshapeM}}

\newcommand\linear{\mathsf{linear}}
\newcommand\Replicate{\mathsf{Replicate}}

\newcommand\ie{{\it i.e.}~}

\newcommand\dd{\,\mathrm{d}} 
\newcommand\sps[1]{^{(#1)}}

\numberwithin{equation}{section}

\newcommand\revised[2]{#2} 

\newcommand\add[1]{\revised{}{#1}}

\title{
  A multiscale neural network based on hierarchical nested bases
}
\date{}
\author{
Yuwei Fan\thanks{Department of Mathematics, Stanford University,
    Stanford, CA 94305, email: {\tt ywfan@stanford.edu}},~~
Jordi Feliu-Fab{\`a}\thanks{Institute for
    Computational and Mathematical Engineering, Stanford University,
    Stanford, CA 94305, email: {\tt jfeliu@stanford.edu}},~~
Lin Lin\thanks{Department of Mathematics, University of California, Berkeley, and Computational
    Research Division, Lawrence Berkeley National Laboratory, Berkeley, CA 94720,
    email: {\tt linlin@math.berkeley.edu}},~~
Lexing Ying\thanks{Department of Mathematics and Institute for
    Computational and Mathematical Engineering, Stanford University,
    Stanford, CA 94305, email: {\tt lexing@stanford.edu}},~~
Leonardo Zepeda-N\'u\~nez\thanks{Computational
    Research Division, Lawrence Berkeley National Laboratory, Berkeley, CA 94720 , email: {\tt lzepeda@lbl.gov }}
}

\begin{document}
\maketitle
\begin{abstract}

In recent years, deep learning has led to impressive results in many fields. In this paper, 
we introduce a multiscale artificial neural network for high-dimensional nonlinear maps based on
the idea of hierarchical nested bases in the fast multipole method and the $\cH^2$-matrices. This
approach allows us to efficiently approximate discretized nonlinear maps arising from partial differential equations or integral equations. It
also naturally extends our recent work based on the generalization of hierarchical matrices [Fan et
al.~arXiv:1807.01883] but with a reduced number of parameters. In particular, the number of
parameters of the neural network grows linearly with the dimension of the parameter space of the
discretized PDE. We demonstrate the properties of the architecture by approximating the solution maps
of nonlinear Schr{\"o}dinger equation, the radiative transfer equation, and the Kohn-Sham map.
\vspace*{4mm}

\end{abstract}



\section{Introduction}\label{sec:intro}

In recent years, deep learning and more specifically deep artificial neural networks have received
ever-increasing attention from the scientific community. Coupled with a significant increase
in the computer power and the availability of massive datasets, artificial neural networks have fueled
several breakthroughs across many fields, ranging  
from classical machine learning applications such as object recognition
~\cite{Krizhevsky2012,Zeiler2014,Zisserman2014,Szegedy2014}, speech recognition \cite{Hinton2012},
natural language processing \cite{Sarikaya14,Socher12} or text classification \cite{Wang2012} to
more modern domains such as 
language translation \cite{SutskeverNIPS2014}, 
drug discovery \cite{MaSheridan2015}, genomics \cite{Leung2014,Xiong2015}, game playing
\cite{Silver2016}, 
among many others.
For a more extensive review of deep learning, we point the reader to 
\cite{leCunn2015,SCHMIDHUBER2015,Goodfellow-et-al-2016}.

Recently, neural networks have also been employed to solve challenging problems in numerical
analysis and scientific
computing~\cite{Araya-Polo2018,Beck2017,berg2017unified,CHAN2018,Chaudhari2017,E2017,khoo2017solving,PASCHALIS2004211,Raissi2018,Rudd2014,Konstantinos2018}.  While a fully
connected neural network can be theoretically used to approximate very general mappings
\cite{CohenSharir2018,Hornik91,Khrulkov2018,Mhaskar2018}, it may also lead to a prohibitively large
number of parameters, resulting in extremely long training stages and overwhelming memory
footprints. Therefore, it is often necessary to incorporate existing knowledge of the underlying
structure of the problem into the design of the network architecture. One promising and general
strategy is to build neural networks based on a multiscale
decomposition~\cite{fan2018mnn,Yingzhou2018,WangChung2018}. The general idea, often used in image
processing
\cite{badrinarayanan2015segnet,Bruna2012,Chen2018DeepLab,LITJENS201760,Ronneberger2015,Ulyanov2018},
is to learn increasingly coarse-grained features of a complex problem across different layers of the
network structure, so that the number of parameters in each layer can be effectively controlled.

In this paper, we aim at employing neural networks to effectively approximate nonlinear maps of the
form
\begin{equation}\label{eq:non_linear_map_intro}
  u = \cM(v),\quad u,v \in \Omega \subset \mathbb{R}^d,
\end{equation}
\add{which can be viewed as a nonlinear generalization of pseudo-differential operators.}
This type of maps may arise from parameterized and discretized partial differential equations (PDE) or
integral equations (IE), with $u$ being the quantity of interest and $v$ the parameter that serves
to identify a particular configuration of the system.

We propose a neural network architecture based on the idea of {\em hierarchical nested bases} used
in the fast multipole method (FMM) \cite{greengard1987fast} and the $\cH^2$-matrix
\cite{Hackbusch2000H2} to represent nonlinear maps arising in computational physics, motivated by
the favorable complexity of the FMM / $\cH^2$-matrices in the linear setting.  The proposed neural
network, which we call MNN-$\cH^2$, is able to efficiently represent the nonlinear maps benchmarked
in the sequel, in such cases the number of parameters required to approximate the maps can grow
linearly with respect to $N$, the dimension of the parameter space of the discretized PDE. Our
presentation will mostly follow the notation of the $\cH^2$-matrix framework due to its algebraic
nature.

The proposed architecture, MNN-$\cH^2$, is a direct extension of the framework
used to build a multiscale neural networks based on $\cH$-matrices (MNN-$\cH$)~\cite{fan2018mnn} to
$\cH^2$-matrices.
We demonstrate the capabilities of MNN-$\cH^2$ with three classical yet challenging examples in
computational physics: the nonlinear Schr\"odinger equation
\cite{anglin2002bose,pitaevskii1961vortex}, the radiative transfer equation
\cite{klose2002optical,koch2004evaluation,marshak20053d,pomraning1973equations}, and the Kohn-Sham
map \cite{HohenbergKohn1964,KohnSham1965}.  
We find that MNN-$\cH^2$ can yield comparable results to those obtained from
MNN-$\cH$, but with a reduced number of parameters, thanks to the use of hierarchical nested bases.

The outline of the paper is as follows.
\Cref{sec:NNhmatrix} reviews the $\cH^2$-matrices and interprets them within the
framework of neural networks.
\Cref{sec:mnn} extends the neural network representation of $\cH^2$-matrices to the nonlinear case.
\Cref{sec:application} discusses the implementation details and demonstrates the accuracy of the
architecture in representing nonlinear maps, followed by the conclusion and future directions in
\cref{sec:conclusion}.

\section{Neural network architecture for $\cH^2$-matrices}\label{sec:NNhmatrix}
In this section, we reinterpret the matrix-vector multiplication
of $\cH^2$-matrices within the framework of neural networks.  In \cref{sec:H2matreview}, we
briefly review $\cH^2$-matrices for the 1D case, and propose the neural network architecture for the
matrix-vector multiplication of $\cH^2$-matrices in \cref{sec:h2nn}. An extension to the
multi-dimensional setting is presented in \cref{sec:nD}.

\subsection{$\cH^2$-matrices} \label{sec:H2matreview}
The concept of hierarchical matrices ($\cH$-matrices) was first introduced by 
Tyrtyshnikov \cite{tyrtyshnikov1996}, and Hackbusch et
al. \cite{hackbusch1999sparse, hackbusch2000sparse} as an
algebraic formulation of algorithms for hierarchical off-diagonal low-rank matrices. This framework
provides efficient numerical methods for solving linear systems arising from integral equations and
partial differential equations \cite{borm2003introduction} and it enjoys an $O(N\log(N))$ arithmetic
complexity for the matrix-vector multiplication.  By incorporating the idea of hierarchical nested
bases from the fast multipole method \cite{greengard1987fast}, the $\cH^2$-matrices were introduced
in \cite{Hackbusch2000H2} to further reduce the logarithmic factor in the complexity, provided that
a so-called ``consistency condition'' is fulfilled.  In the sequel, we follow the notation
introduced in \cite{fan2018mnn} to provide a brief introduction to the framework of $\cH^2$-matrices
in a simple uniform Cartesian setting.  We refer readers to
\cite{borm2003introduction,Hackbusch2000H2, lin2011fast} for further details.


Consider the integral equation
\begin{equation}\label{eq:integral}
    u(x) = \int_{\Omega}g(x,y)v(y) \dd y, \quad \Omega=[0,1),
\end{equation}
where $u$ and $v$ are periodic in $\Omega$ and $g(x,y)$ is smooth and numerically low-rank away from
the diagonal. A discretization with an uniform grid with $N=2^Lm$ discretization points yields the
linear system given by
\begin{equation}\label{eq:discrete}
    u = A v,
\end{equation}
where $ A\in\bbR^{N\times N}$, and $u, v \in\bbR^N$ are the discrete analogs of $u(x)$ and $v(x)$
respectively.

\begin{figure}[ht!]
    \centering
\begin{minipage}{0.9\textwidth}
    \centering
    \subfloat[Illustration of computational domain for an interior segment (up) and a
    boundary segment (down).]{
    \label{fig:partitionmesh}
    \includegraphics[width=0.6\textwidth]{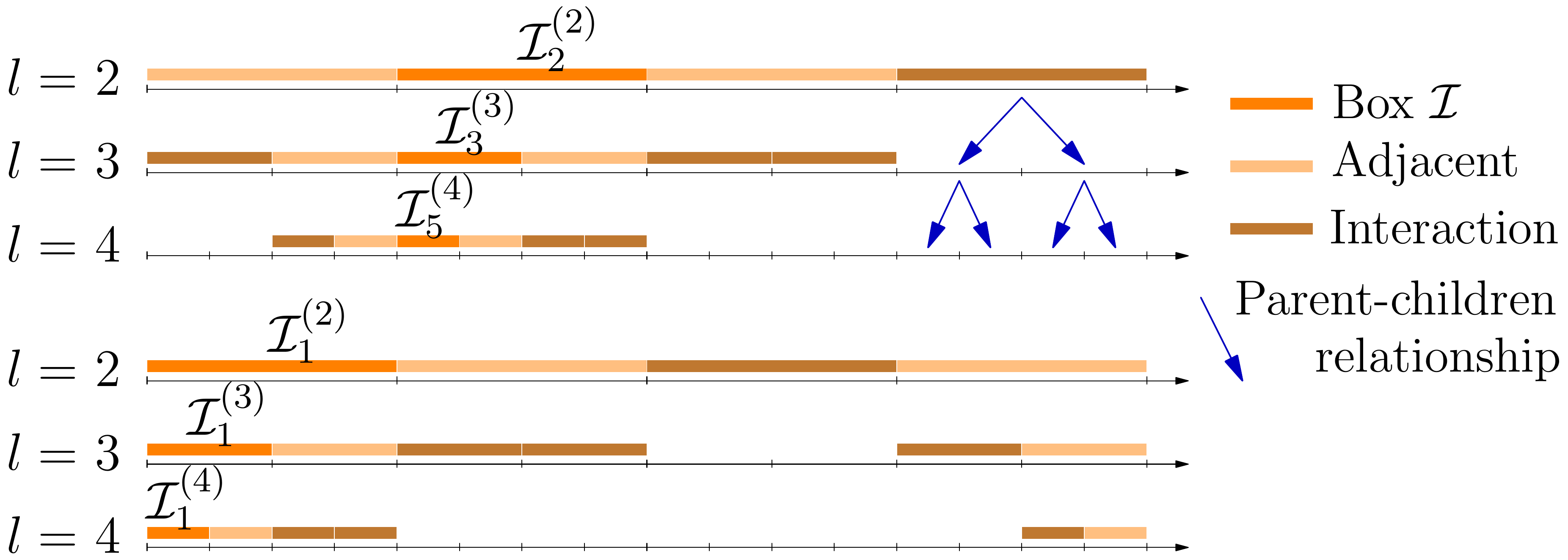}
    }\hspace{0.1\textwidth}
    \subfloat[Hierarchical partition of matrix $A$]{
    \label{fig:partitionmatrix}
    \includegraphics[width=0.22\textwidth]{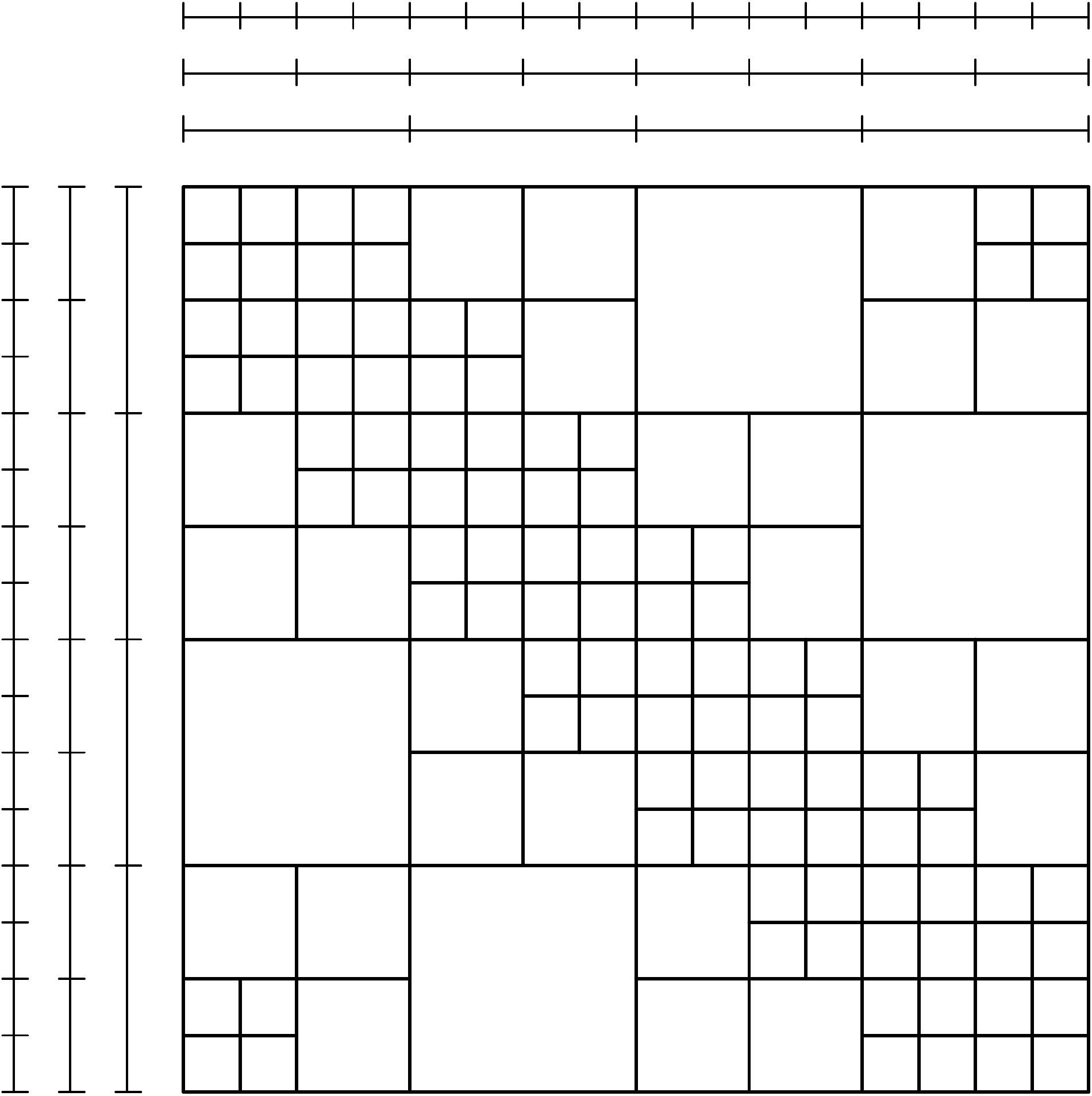}
    }\\
    \subfloat[Decomposition of matrix $A$]{
    \label{fig:decomposition}
    \includegraphics[width=0.95\textwidth]{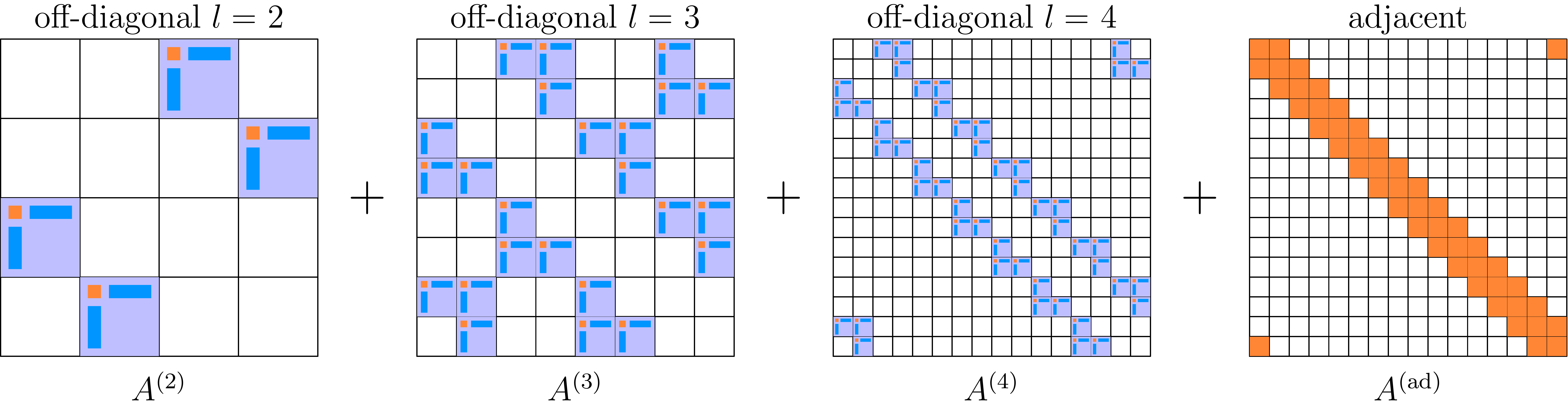}
    }
\end{minipage}
    \caption{\label{fig:partition}Hierarchical partition of computational domain, its corresponding
    partition of matrix $A$ and the decomposition of matrix $A$.}
\end{figure}
A hierarchical dyadic decomposition of the grid in $L+1$ levels can be introduced as follows.
Let $\cI\sps{0}$, the $0$-th level of the decomposition, be the set of all grid points defined as
\begin{equation}\label{eq:cI}
    \cI\sps{0} = \{k/N: k=0,\dots,N-1\}.
\end{equation}
At each level $\ell$ ($1\leq \ell\leq L$), the grid is decomposed in $2^{\ell}$ disjoint
\emph{segments}.
Each segment is defined by $\cI\sps{\ell}_i = \cI\sps{0} \cap[(i-1)/2^{\ell}, i/2^{\ell})$ for
$i=1,\dots,2^{\ell}$. Throughout this manuscript, $\cI\sps{\ell}$(or $\cJ\sps{\ell}$) will denote a
generic segment of a given level $\ell$, and the superscript $\ell$ will be omitted when the level
is clear from the context.

Given a vector $v\in\bbR^N$, we denote $v_{\cI}$ the elements of $v$ indexed by $\cI$; and given a
matrix $A\in\bbR^{N\times N}$, we denote $A_{\cI,\cJ}$ the submatrix of $A$ indexed by $\cI\times
\cJ$.
Following the usual nomenclature in $\cH$-matrices, we define the following relationships between
segments:
\begin{itemize}
    \item[$\cC(\cI)$] \emph{children list} of $\cI$ for $\ell<L$: list of the segments on level
        $\ell+1$ that are subset of $\cI$;
    \item[$\cP(\cI)$] \emph{parent} of $\cI$ for $\ell>0$: set of segments $\cJ$ such that
        $\cI\in\cC(\cJ)$;
    \item[$\NL(\cI)$] \emph{neighbor list} of  $\cI$: list of the segments on level $\ell$ that are
        adjacent to $\cI$ including $\cI$ itself;
    \item[$\IL(\cI)$] \emph{interaction list} of $\cI$ for $\ell\geq2$: set that contains all the 
        segments on level $\ell$ that are children of segments in $\NL(\cP(\cI))$ minus $\NL(\cI)$, \ie
        $\IL(\cI) = \cC(\NL(\cP(\cI))) - \NL(\cI)$.
\end{itemize}


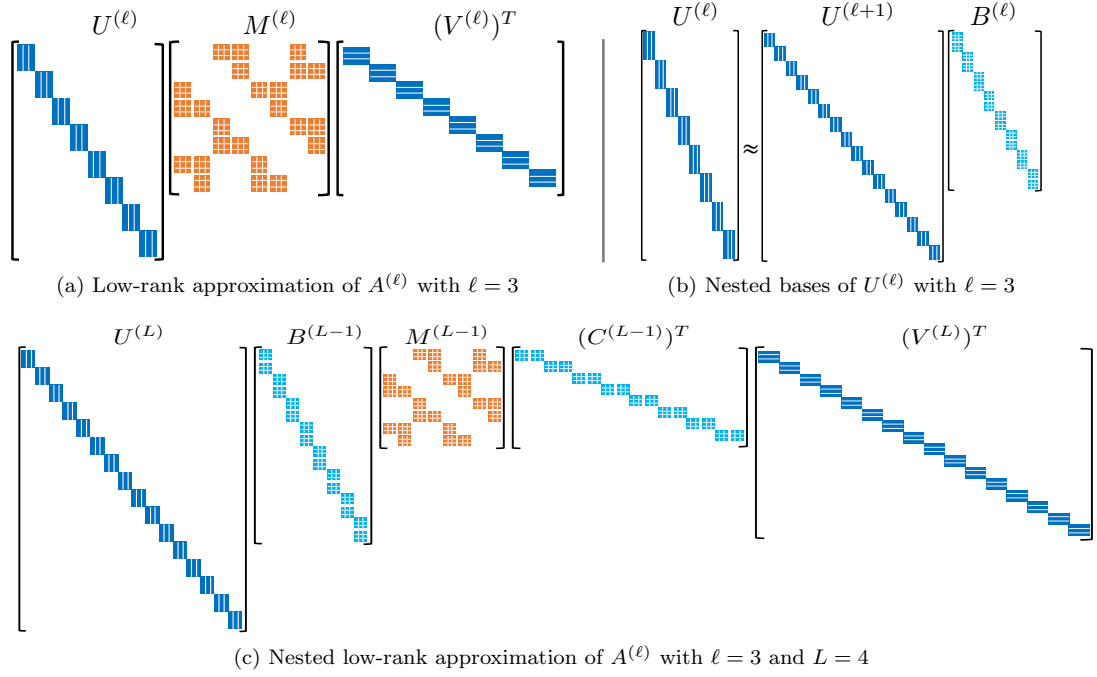
\begin{figure}[ht]
    \centering
\begin{minipage}{0.9\textwidth}
    \subfloat[Low-rank approximation of $A\sps{\ell}$ with $\ell=3$]{
    \label{fig:lowrank}
    \begin{overpic}[width=0.5\textwidth]{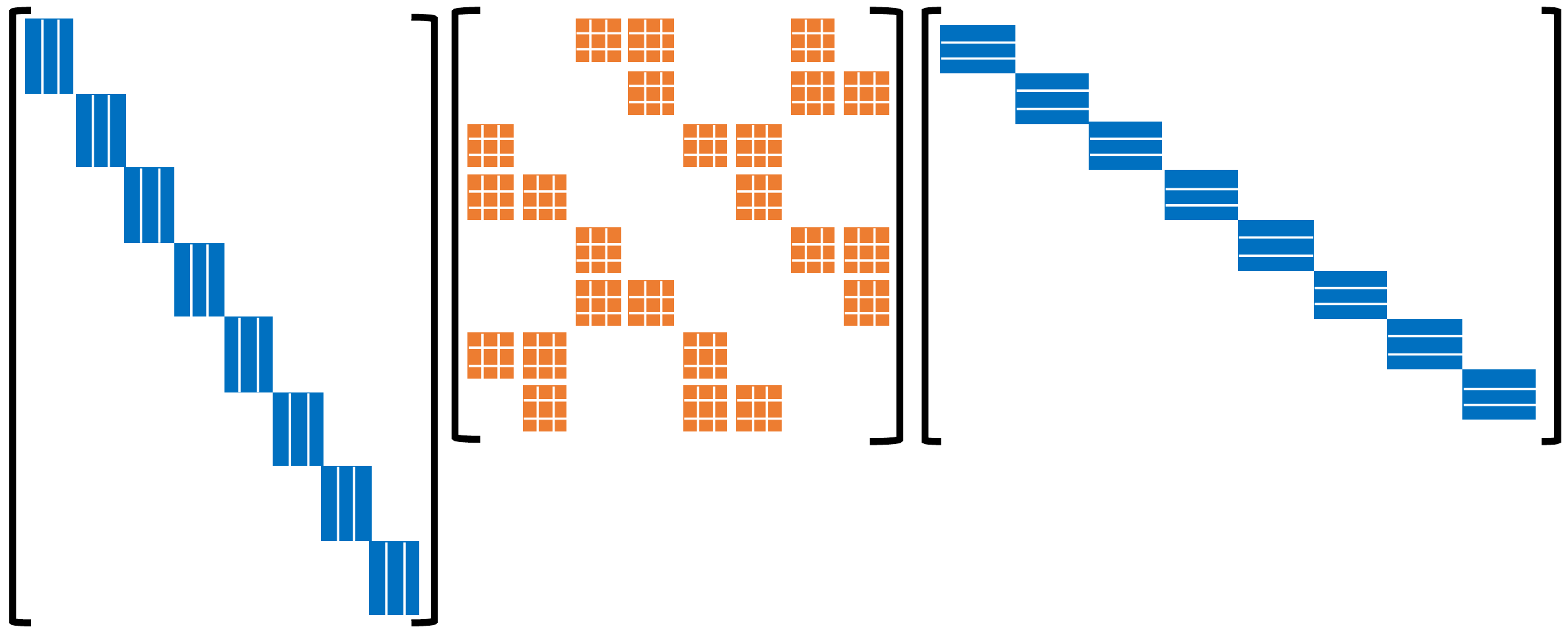}
        \put(15, 41){$U\sps{\ell}$}
        \put(42, 41){$M\sps{\ell}$}
        \put(76, 41){$(V\sps{\ell})^T$}
    \end{overpic}
    }
    \quad {\color{gray}\rule{1pt}{0.20\textwidth}}\quad
    \subfloat[Nested bases of $U\sps{\ell}$ with $\ell=3$]{
    \label{fig:nest}
    \begin{overpic}[width=0.36\textwidth]{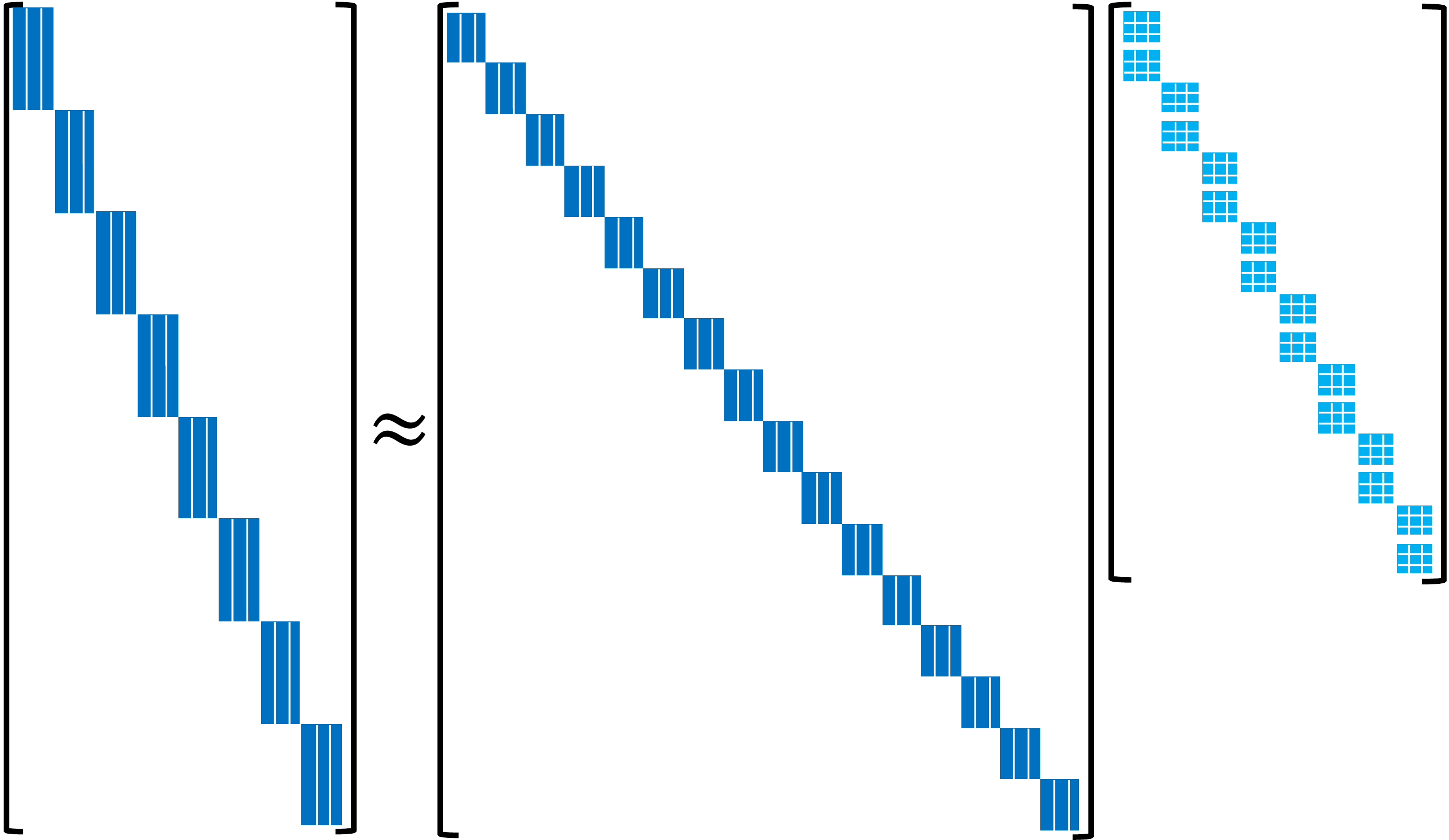}
        \put( 8, 59){$U\sps{\ell}$}
        \put(45, 59){$U\sps{\ell+1}$}
        \put(82, 59){$B\sps{\ell}$}
    \end{overpic}
    }\\[2mm]
    \subfloat[Nested low-rank approximation of $A\sps{\ell}$ with $\ell=3$ and $L=4$]{
    \label{fig:A3}
    \resizebox{0.98\textwidth}{0.26\textwidth}{
    \begin{overpic}[width=\textwidth]{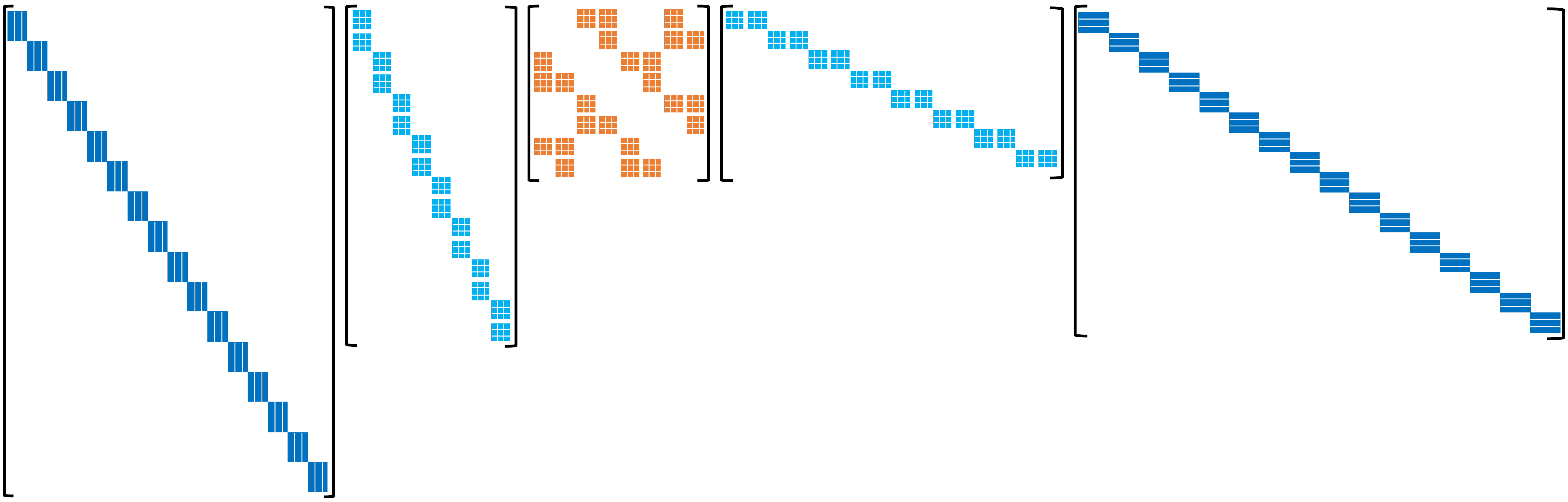}
        \put( 9, 32){$U\sps{L}$}
        \put(25, 32){$B\sps{L-1}$}
        \put(36, 32){$M\sps{L-1}$}
        \put(52, 32){$(C\sps{L-1})^T$}
        \put(82, 32){$(V\sps{L})^T$}
    \end{overpic}
    }}
\end{minipage}
    \caption{\label{fig:factorization}Low rank factorization and nested low rank factorization of
    $A\sps{l}$.}
\end{figure}

\cref{fig:partitionmesh} illustrates this dyadic hierarchical partition of the computational domain,
the parent-children relationship, the \emph{neighbor} list, and \emph{interaction} list on levels
$\ell=2,3,4$. The matrix $A$ can be hierarchically partitioned as illustrated in
\cref{fig:partitionmatrix}. The partition leads to a multilevel decomposition of $A$ shown in
\cref{fig:decomposition}, which can be written as
\begin{equation}\label{eq:decompose}
    A = \sum_{\ell=2}^L A\sps{\ell} + A\sps{\ad},
    \qquad
    \begin{aligned}
    A\sps{\ell}_{\cI,\cJ} &=\begin{cases}
        A_{\cI,\cJ}, & \cI \in \IL(\cJ);\\
        0, & \text{otherwise},
    \end{cases}
    \quad \cI,\cJ \text{ at level } l,\quad 2\leq l\leq L,
    \\
    A_{\cI,\cJ}\sps{\ad} &=\begin{cases}
        A_{\cI,\cJ}, & \cI \in \NL(\cJ); \\
        0, & \text{otherwise},
    \end{cases}
    \quad \cI,\cJ \text{ at level $L$}.
    \end{aligned}
\end{equation}

For simplicity, we suppose that each block has a fixed numerical rank at most $r$, \ie,
\begin{equation} \label{eq:low_rank_fact}
  A\sps{\ell}_{\cI, \cJ}\approx U\sps{\ell}_{\cI}M\sps{\ell}_{\cI,\cJ}(V\sps{\ell}_{\cJ})^{T},
  \quad U\sps{\ell}_{\cI}, V\sps{\ell}_{\cJ}\in\bbR^{N / 2^{\ell}\times r},\quad
  M\sps{\ell}_{\cI,\cJ}\in\bbR^{r\times r},
\end{equation}
where $\cI$ and $ \cJ$ are any interacting segments at level $\ell$. We can approximate
$A\sps{\ell}$ as $A\sps{\ell} \approx U\sps{\ell}M\sps{\ell}(V\sps{\ell})^T$ as depicted in
\cref{fig:lowrank}.  Here $U\sps{\ell}$, $V\sps{\ell}$ are block diagonal matrices with diagonal
blocks $U\sps{\ell}_{\cI}$ and $V\sps{\ell}_{\cI}$ for $\cI$ at level $\ell$, respectively, and
$M\sps{\ell}$ aggregates all the blocks $M\sps{\ell}_{\cI,\cJ}$ for all interacting segments
$\cI,\cJ$ at level $\ell$.

The key feature of $\cH^2$-matrices is that the bases matrices $U\sps{\ell}_{\cI}$ and
$V\sps{\ell}_{\cI}$ between parent and children segments enjoy a {\em nested} low rank
approximation.
More precisely, if $\cI$ is at level $2\leq l<L$ and $\cJ_1, \cJ_2\in\cC(\cI)$ are at level $l+1$,
then $U\sps{l}_{\cI}$ and $V\sps{l}_{\cI}$ satisfy the following approximation:
\begin{equation} \label{eq:coeff_mat}
 U\sps{\ell}_{\cI}\approx
 \begin{pmatrix}
     U\sps{\ell+1}_{\cJ_1} & \\
                           & U\sps{\ell+1}_{\cJ_2}
 \end{pmatrix}
 \begin{pmatrix}
     B\sps{\ell}_{\cJ_1}\\[3mm]
     B\sps{\ell}_{\cJ_2}
 \end{pmatrix},\quad
 V\sps{\ell}_{\cI}\approx
 \begin{pmatrix}
     V\sps{\ell+1}_{\cJ_1} & \\
                           & V\sps{\ell+1}_{\cJ_2}
 \end{pmatrix}
 \begin{pmatrix}
     C\sps{\ell}_{\cJ_1}\\[3mm]
     C\sps{\ell}_{\cJ_2}
 \end{pmatrix},
\end{equation}
where $B\sps{\ell}_{\cJ_i},C\sps{\ell}_{\cJ_i}\in\bbR^{r\times r}$, $i=1,2$. As depicted by
\cref{fig:nest}, if we introduce the matrix $B\sps{l}$ ($C\sps{l}$) that aggregates all the blocks
$B\sps{\ell}_{\cJ_i}$ ($C\sps{\ell}_{\cJ_i}$) for all the parent-children pairs ($\cI$, $\cJ_i$),
\eqref{eq:coeff_mat} can be compactly written as $U\sps{l}\approx U\sps{l+1}B\sps{l}$ and
$V\sps{l}\approx V\sps{l+1}C\sps{l}$. Thus, the decomposition \eqref{eq:decompose} can be further
factorized as
\begin{small}
\begin{equation}\label{eq:H2_A}
    \begin{aligned}
    A &=\sum_{\ell=2}^LA\sps{\ell}+A\sps{\ad}
    \approx \sum_{\ell=2}^L U\sps{L}B\sps{L-1}\cdots B\sps{\ell}M\sps{\ell}
    (C\sps{\ell})^T \cdots (C\sps{L-1})^T(V\sps{L})^T
    + A\sps{\ad}.
    \end{aligned}
\end{equation}
\end{small}%
The matrix-vector multiplication of $A$ with an arbitrary vector $v$ can be approximated by
\begin{equation}\label{eq:H2_Av}
    Av \approx \sum_{\ell=2}^L U\sps{L}B\sps{L-1}\cdots B\sps{\ell}M\sps{\ell}
    (C\sps{\ell})^T \cdots (C\sps{L-1})^T(V\sps{L})^Tv
    + A\sps{\ad}v.
\end{equation}

\begin{algorithm}[htb]
\begin{small}
\begin{center}
\begin{minipage}{0.99\textwidth}
\begin{minipage}{0.6\textwidth}
\begin{algorithmic}[1]
    \State $u\sps{\ad}=A\sps{\ad}v$;
    \State $\xi\sps{L}=(V\sps{L})^Tv$;
    \For {$\ell$ from $L-1$ to 2 by $-1$}
    \State $\xi\sps{\ell}=(C\sps{\ell})^T\xi\sps{\ell+1}$;
    \EndFor
    \For {$\ell$ from $2$ to $L$}
    \State $\zeta\sps{\ell}=M\sps{\ell}\xi\sps{\ell}$;
    \EndFor
    \algstore{break1}
\end{algorithmic}
\end{minipage}
\begin{minipage}{0.35\textwidth}
\begin{algorithmic}[1]
    \algrestore{break1}
    \State $\chi=0$;
    \For {$\ell$ from 2 to $L-1$}
    \State $\chi=\chi+\zeta\sps{\ell}$;
    \State $\chi=B\sps{\ell}\chi$;
    \EndFor
    \State $\chi=\chi+\zeta\sps{L}$;
    \State $\chi=U\sps{L}\chi$;
    \State $u = \chi + u\sps{\ad}$;
\end{algorithmic}
\end{minipage}
\end{minipage}
\end{center}
\end{small}
\caption{Application of $\cH^2$-matrices on a vector $v\in\bbR^N$.}%
\label{alg:H2}%
\end{algorithm}

\cref{alg:H2} provides the implementation of the matrix-vector multiplication of $\cH^2$-matrices.
The key properties of the matrices $U\sps{L}$, $V\sps{L}$, $B\sps{\ell}$, $C\sps{\ell}$,
$M\sps{\ell}$ and $A\sps{\ad}$ are summarized as follows:
\begin{property}\label{pro:matrices}
    The matrices
    \begin{enumerate}
        \item $U\sps{L}$ and $V\sps{L}$ are block diagonal matrices with block size $N/2^L\times r$;
        \label{pro:UV}
        \item $B\sps{\ell}$ and $C\sps{\ell}$, $\ell=2,\cdots,L-1$ are block diagonal matrices with block size $2r\times r$;
        \label{pro:BC}
        \item $M\sps{\ell}$, $\ell=2,\cdots,L$ are block cyclic band matrices with block size
            $r\times r$ and band size $n\sps{\ell}_b$, which is $2$ for $\ell=2$ and $3$ for $\ell>2$;
        \label{pro:M}
        \item $A\sps{ad}$ is a block cyclic band matrix with block size $m\times m$ with band size
            $n_b\sps{\ad}=1$.
        \label{pro:Aad}
    \end{enumerate}
\end{property}

\subsection{Matrix-vector multiplication as a neural network}\label{sec:h2nn}

We represent the matrix-vector multiplication \eqref{eq:H2_Av} using the
framework of neural networks. We first introduce our main tool --- locally connected network --- in
\cref{sec:LC} and then present the neural network representation of \eqref{eq:H2_Av} in
\cref{sec:NNh2}.

\subsubsection{Locally connected network}\label{sec:LC}
\begin{figure}[ht]
    \centering
    \subfloat[$\alpha=\alpha'=1$]{
        \includegraphics[page=1,width=0.25\textwidth]{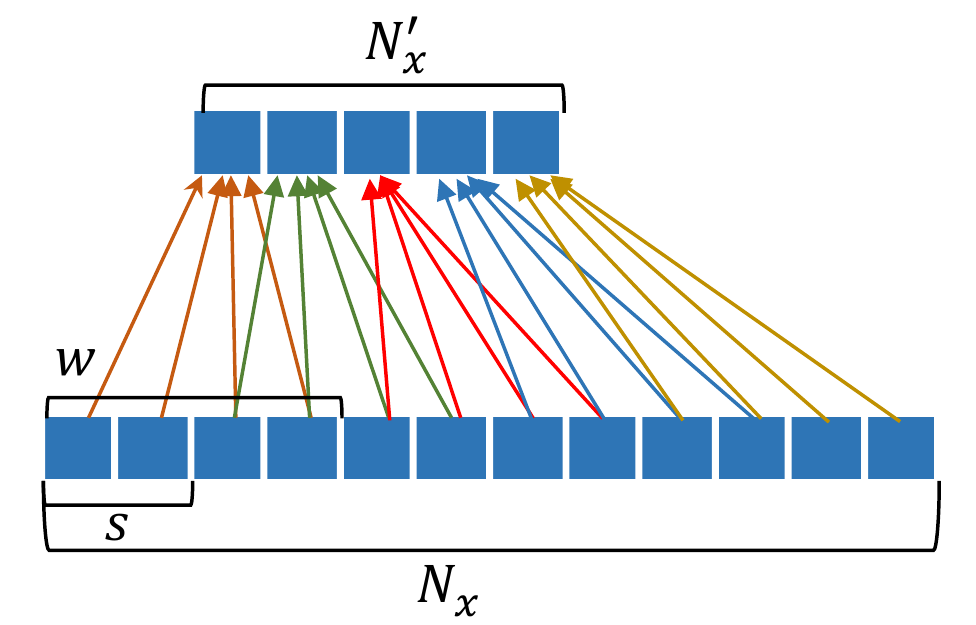}
    }\hspace{0.2\textwidth}
    \subfloat[$\alpha=2$, $\alpha'=3$]{
        \includegraphics[page=2,width=0.25\textwidth]{LC.pdf}
    }
    \caption{\label{fig:lcsample}Sample of LC network with $N_x=12$, $s=2$, $w=4$ and $N_x'=5$.}
\end{figure}

In order to simplify the notation, let us present the 1D case as an example.  In this setup, an NN
layer can be represented by a 2-tensor with size $\alpha\times N_x$, where $\alpha$ is called the
\emph{channel dimension} and $N_x$ is usually called the \emph{spatial dimension}.  A locally
connected network is a type of mapping between two adjacent layers, where the output of each neuron
depends only locally on the input. If a layer $\xi$ with size $\alpha\times N_x$ is connected to a
layer $\zeta$ with size $\alpha'\times N_x'$ by a \emph{locally connected} (LC) network, then
\begin{equation}\label{eq:lc}
    \zeta_{c',i} = \phi\left( \sum_{j=(i-1)s+1}^{(i-1)s+w}\sum_{c=1}^{\alpha}
    W_{c',c;i,j}\xi_{c,j} + b_{c',i}\right),
    \quad i=1,\dots,N_x', ~ c'=1,\dots,\alpha',
\end{equation}
where $\phi$ is a pre-specified function, called \emph{activation}, usually chosen to be e.g. a
linear function, a rectified-linear unit (ReLU) function or a sigmoid function. The parameters $w$
and $s$ are called the \emph{kernel window size} and \emph{stride}, respectively.
\cref{fig:lcsample} presents a sample of the LC network.  Furthermore, we call the layer $\zeta$
\emph{locally connected layer} (LC layer) hereafter.

In \eqref{eq:lc} the LC network is represented using tensor notation; however, we can reshape
$\zeta$ and $\xi$ to a vector by column major indexing and $W$ to a matrix and write \eqref{eq:lc}
into a matrix-vector form as
\begin{equation}\label{eq:lcmat}
    \zeta = \phi(W\xi+b).
\end{equation}
For later usage, we define $\Reshape[n_1,n_2]$ to be the map that reshapes a tensor with size
$n_1'\times n_2'$ to a 2-tensor of size $n_1\times n_2$ such that $n_1n_2=n_1'n_2'$ by column major
indexing.  Here, we implicitly regard a vector with size $n$ as a 2-tensor with size $1\times
n$. 

\begin{figure}[ht]
    \begin{minipage}{\textwidth}
        \begin{tabular}{ >{\centering\arraybackslash}m{0.30\textwidth}|
            >{\centering\arraybackslash}m{0.30\textwidth}|
            >{\centering\arraybackslash}m{0.30\textwidth}}
            \includegraphics[page=1, width=0.30\textwidth]{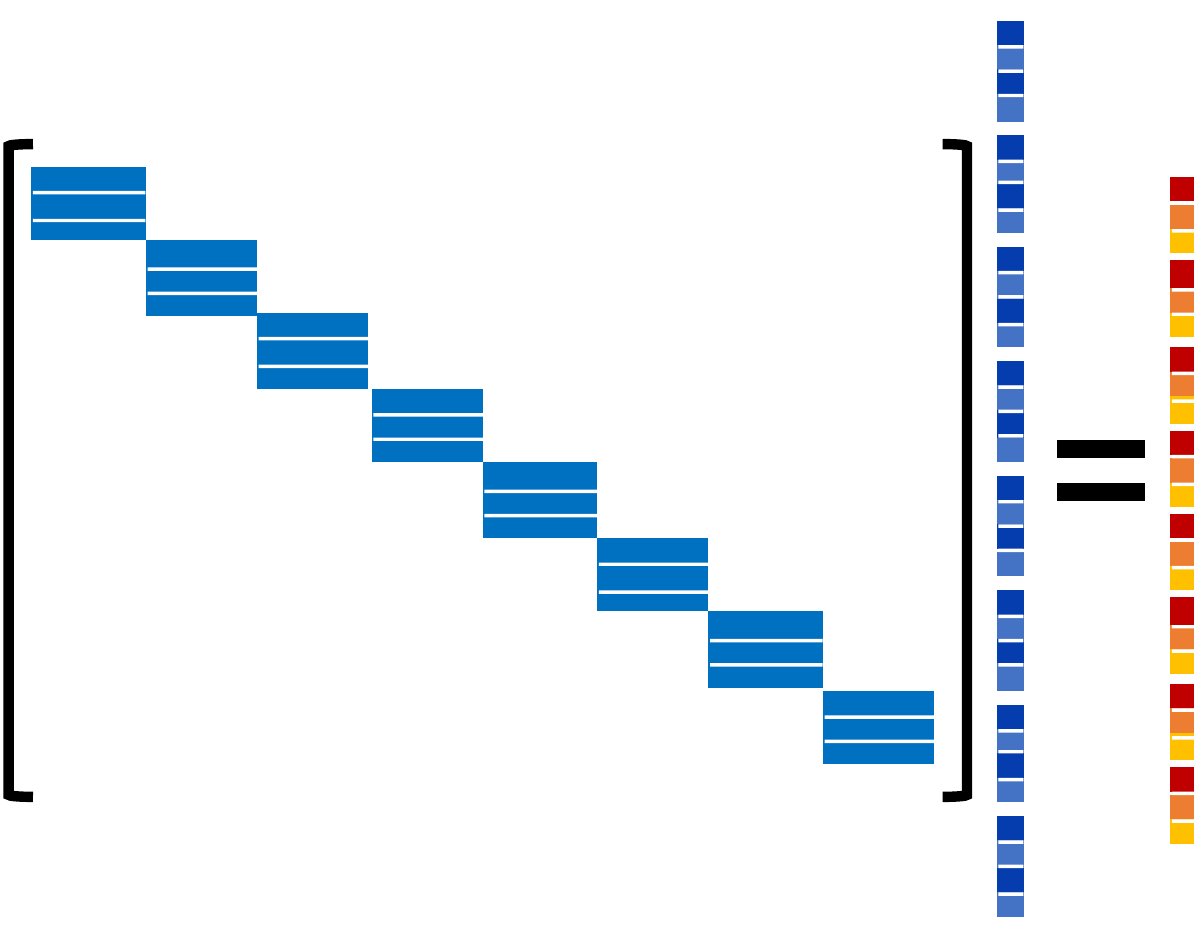} &
            \includegraphics[page=2, width=0.30\textwidth]{LCmats2.pdf} &
            \includegraphics[page=3, width=0.30\textwidth]{LCmats2.pdf} \\
            $s=w=\frac{N_x}{N_x'}$,  &
            $s=1$, $N_x'=N_x$ &
            $s=1$, $w=1$, $N_x'=N_x$ \\[4mm]
            \includegraphics[page=1, width=0.30\textwidth]{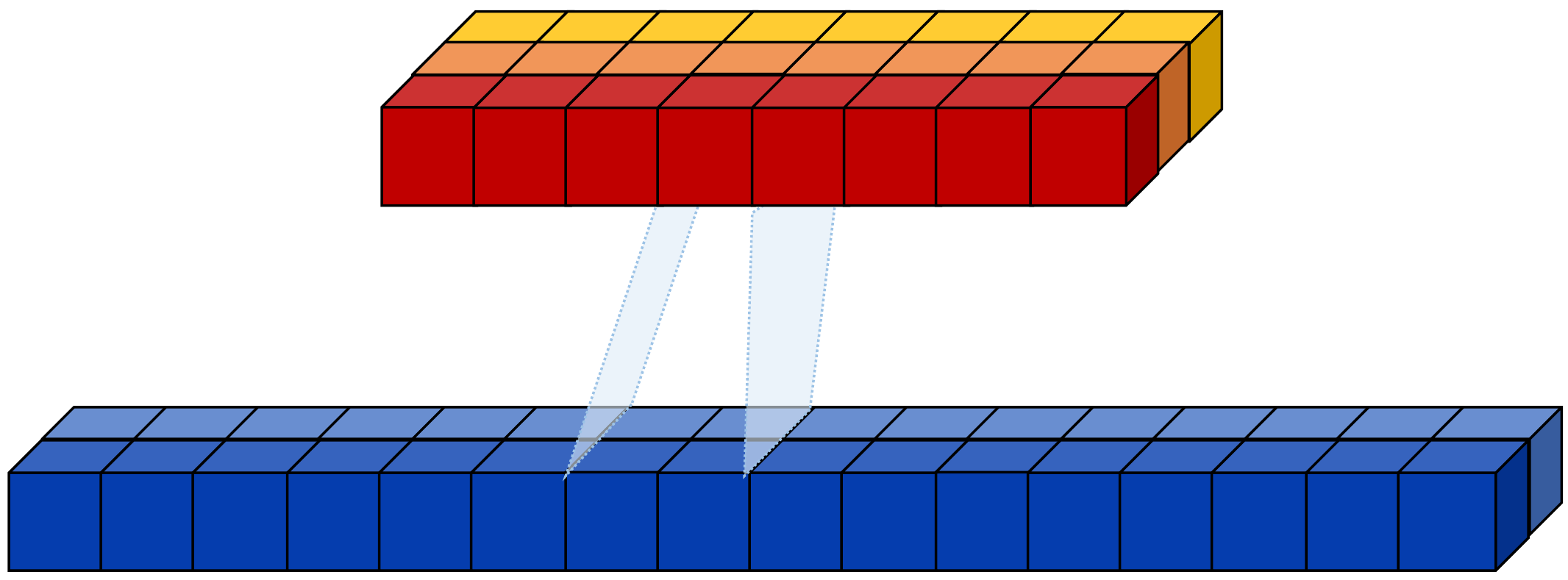} &
            \includegraphics[page=2, width=0.30\textwidth]{LCs.pdf} &
            \includegraphics[page=3, width=0.30\textwidth]{LCs.pdf} \\
            (a) $\LCR[\phi; N_x, \alpha, N_x', \alpha']$ with
            $N_x=16$, $\alpha=2$, $N_x'=8$ and $\alpha'=3$ &
            (b) $\LCK[\phi; N_x, \alpha, \alpha', w]$
            with $N_x=8$, $\alpha=\alpha'=3$ and $w=3$ &
            (c) $\LCI[\phi; N_x, \alpha, \alpha']$
            with $N_x=8$, $\alpha=3$ and $\alpha'=4$ \\
        \end{tabular}
    \end{minipage}
    \caption{\label{fig:LC}Three instances of locally connected networks used to represent the
    matrix-vector multiplication.
    The upper portions of each column depict the patterns of the matrices and the lower portions
    are their respective analogs using locally connect networks.}
\end{figure}

Each LC network has $6$ parameters, $N_x$, $\alpha$, $N_x'$, $\alpha'$, $w$ and $s$. We
define three types of LC networks by specifying some of their parameters.
The upper figures in \cref{fig:LC} depict its corresponding formula in matrix-vector form
\eqref{eq:lcmat}, and the lower figures show a diagram of the map.
\begin{itemize}
    \item[$\LCR$] \emph{Restriction} map: set $s=w=\frac{N_x}{N_x'}$ in LC.
        This map represents the multiplication of a block diagonal matrix with block sizes
        $\alpha'\times s\alpha$ and a vector with size $N_x\alpha$.  We denote this map by
        $\LCR[\phi; N_x, \alpha, N_x', \alpha']$. The application of $\LCR[\linear; 16, 2, 8, 3]$ is
        depicted in \cref{fig:LC}a. 
    \item[$\LCK$] \emph{Kernel} map: set $s=1$ and $N_x'=N_x$.  This map represents the
      multiplication of a periodically banded block matrix (with block size $\alpha'\times \alpha$
      and band size $\frac{w-1}{2}$) with a vector of size $N_x\alpha$.  To account for the
      periodicity, we periodic pad the input layer $\xi_{c,j}$ on the spatial dimension to the size
      $(N_x+w-1)\times \alpha$. We denote this map by $\LCK[\phi; N_x, \alpha, \alpha', w]$, which
      contains two steps: the periodic padding of $\xi_{c,j}$ on the spatial dimension, and the
      application of \eqref{eq:lc}.  The application of $\LCK[\linear; 8,3, 3, 3]$ is depicted in
      \cref{fig:LC}b.
    \item[$\LCI$] \emph{Interpolation} map: set $s=w=1$ and $N_x'=N_x$ in LC. This map represents
        the multiplication of a block diagonal matrix with block size $\alpha'\times \alpha$, times
        a vector of size $N_x\alpha$.  We denote the map by $\LCI[\phi; N_x, \alpha, \alpha']$.
        The application of $\LCI[\linear; 8, 3, 4]$ is depicted in \cref{fig:LC}c.
\end{itemize}

\subsubsection{Neural network representation}\label{sec:NNh2}
We need to find a neural network representation of the following 6 operations in order to perform the matrix-vector multiplication \eqref{eq:H2_Av} for $\cH^2$-matrices: 
\begin{subequations}\label{eq:hm}
        \begin{align}
            \label{eq:Vv} \xi\sps{L} &= (V\sps{L})^T v,\\
            \label{eq:Cv} \xi\sps{\ell} &= (C\sps{\ell})^T\xi\sps{\ell+1}, \quad 2\leq \ell < L,\\
            \label{eq:Mv} \zeta\sps{\ell} &= M\sps{\ell}\xi\sps{\ell},\qquad 2\leq \ell \leq L,\\
            \label{eq:Bv} \chi\sps{\ell} &= B\sps{\ell}\zeta\sps{\ell},\qquad 2\leq \ell < L,\\
            \label{eq:Uv} \chi\sps{L} &= U\sps{L}\zeta\sps{L},\\
            \label{eq:Aadv} u\sps{\ad} &= A\sps{\ad}v.
        \end{align}
\end{subequations}
\begin{figure}[ht]
    \centering
    \begin{overpic}[width=0.95\textwidth]{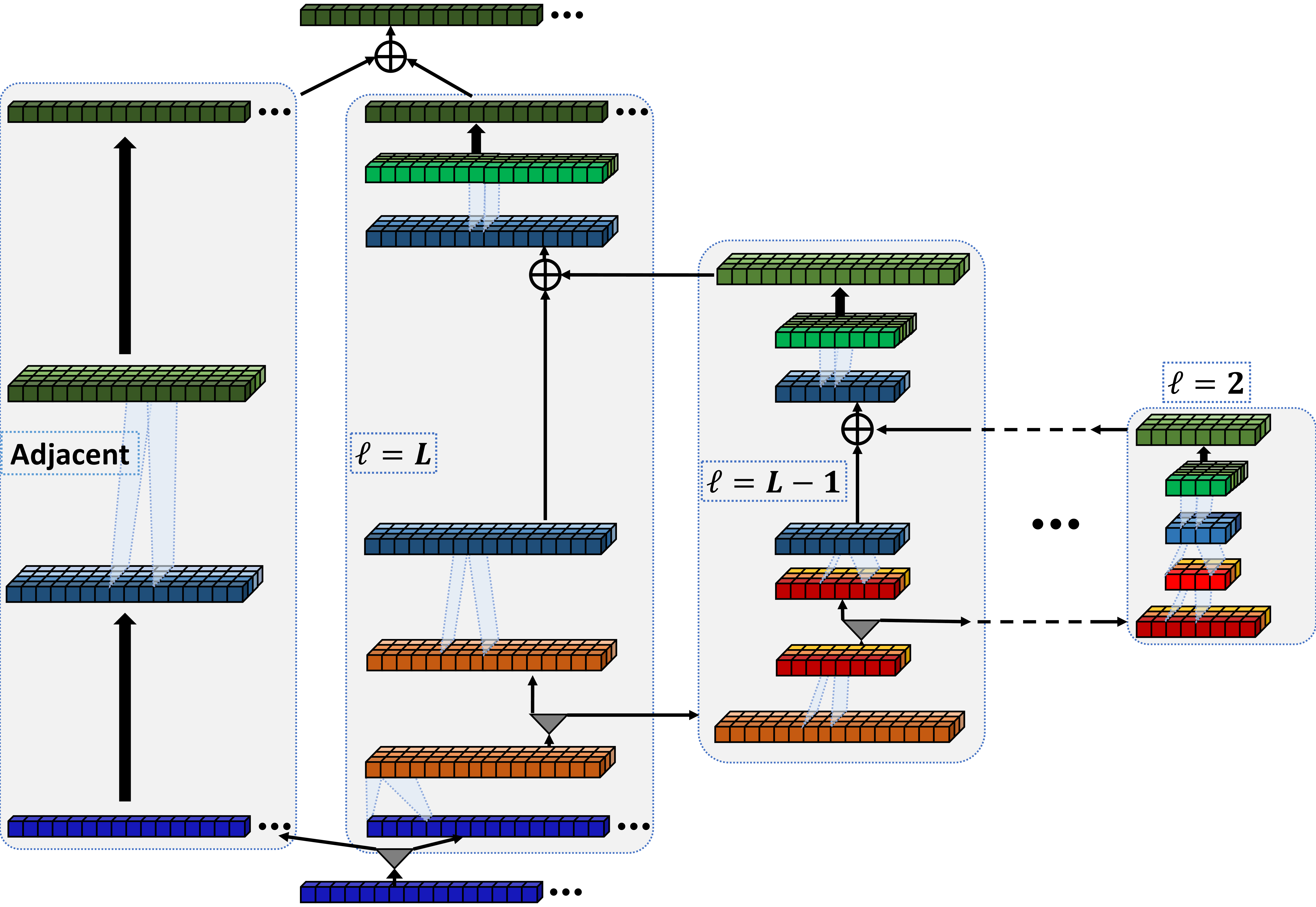}
        \put(12.0, 15.0){\small$\Reshape$}
        \put(12.0, 50.0){\small$\Reshape$}
        \put(11.5, 32.0){\footnotesize$\LCK$-$\linear$}
        \put(34.0,  2.5){\footnotesize$\Replicate$}
        \put(32.0, 64.7){\footnotesize$\mathsf{sum}$}
        \put(38.5, 57.4){\footnotesize$\Reshape$}
        \put(41.0, 15.0){\footnotesize$\Replicate$}
        \put(42.5, 48.7){\footnotesize$\mathsf{sum}$}
        \put(37.0,  7.5){\scriptsize$\LCR$-$\linear$}
        \put(39.0, 22.8){\scriptsize$\LCK$-$\linear$}
        \put(40.0, 52.8){\scriptsize$\LCI$-$\linear$}
        \put(66.0, 22.2){\tiny$\Replicate$}
        \put(67.0, 36.9){\footnotesize$\mathsf{sum}$}
        \put(65.0, 45.4){\scriptsize$\Reshape$}
        \put(67.0, 16.0){\tiny$\LCR$-$\linear$}
        \put(68.0, 26.2){\tiny$\LCK$-$\linear$}
        \put(68.0, 41.2){\tiny$\LCI$-$\linear$}
        \put(93.0, 23.1){\tiny$\LCR$-$\linear$}
        \put(93.0, 26.6){\tiny$\LCK$-$\linear$}
        \put(93.0, 30.1){\tiny$\LCI$-$\linear$}
        \put(94.0, 33.8){\tiny$\Reshape$}
    \end{overpic}
    \caption{\label{fig:h2nn}Neural network architecture for the matrix-vector multiplication of
    $\cH^2$-matrices.}
\end{figure}

%

Following \cref{pro:matrices}.\ref{pro:UV} and the definition of $\LCR$, we can directly represent
\eqref{eq:Vv} as
\begin{equation}\label{eq:LCV}
    \text{\eqref{eq:Vv}}\Rightarrow \xi\sps{L} = \LCR[\linear; N, 1, 2^L, r](v).
\end{equation}
Here we note that the output of $\LCR$ is a 2-tensor, so we should reshape it to a vector. In the
next step, when applying other operations, it is reshaped back to a 2-tensor with same size. These
operations usually do not produce any effect on the whole pipeline, so they are omitted in the
following discussion. Similarly, since all of $V\sps{L}$, $B\sps{\ell}$ and $C\sps{\ell}$ are block
diagonal matrices (\cref{pro:matrices}.\ref{pro:UV} and \cref{pro:matrices}.\ref{pro:BC}),
\begin{equation}\label{eq:LCCBU}
    \begin{aligned}
    \text{\eqref{eq:Cv}} & \Rightarrow
    \xi\sps{\ell} = \LCR[\linear;2^{\ell+1},r,2^{\ell},r](\xi\sps{\ell+1}),\\
    \text{\eqref{eq:Bv}} & \Rightarrow
    \chi\sps{\ell} = \LCI[\linear;2^{\ell},r,2r](\zeta\sps{\ell}),\\
    \text{\eqref{eq:Uv}} & \Rightarrow
    \chi\sps{L} = \LCI[\linear;2^{L},r,m](\zeta\sps{\ell}).
    \end{aligned}
\end{equation}
Analogously, using \cref{pro:matrices}.\ref{pro:M}, \cref{pro:matrices}.\ref{pro:Aad} and the
definition of $\LCK$, 
\begin{equation}\label{eq:LCMA}
    \begin{aligned}
    \text{\eqref{eq:Mv}} & \Rightarrow
    \zeta\sps{\ell} = \LCK[\linear;2^{\ell},r,r,2n_b\sps{\ell}+1](\xi\sps{\ell}),\\
    \text{\eqref{eq:Aadv}} & \Rightarrow
    u\sps{\ad} = \LCK[\linear;2^{L},m,m,2n_b\sps{\ad}+1](v).
    \end{aligned}
\end{equation}

\begin{algorithm}[htb]
\begin{small}
\begin{center}
\begin{minipage}{0.99\textwidth}
\begin{minipage}{0.6\textwidth}
\begin{algorithmic}[1]
    \State $\tilde{v}=\Reshape[m,2^L](v)$;
    \State $\tilde{u}\sps{\ad}=\LCK[\linear; 2^L,m,m,2n_{b}\sps{\ad}+1](\tilde{v})$;
    \State $u\sps{\ad}=\Reshape[1,N](\tilde{u}\sps{\ad})$;
    \State $\xi\sps{L}=\LCR[\linear; N, 1, 2^L, r](v)$;
    \For {$\ell$ from $L-1$ to 2 by $-1$}
    \State $\xi\sps{\ell}=\LCR[\linear; 2^{\ell+1}, r, 2^\ell, r](\xi\sps{\ell+1})$;
    \EndFor
    \For {$\ell$ from $2$ to $L$}
    \State $\zeta\sps{\ell}=\LCK[\linear; 2^\ell, r, r, 2n_b\sps{\ell}+1](\xi\sps{\ell})$;
    \EndFor
    \algstore{break2}
\end{algorithmic}
\end{minipage}
\begin{minipage}{0.38\textwidth}
\begin{algorithmic}[1]
    \algrestore{break2}
    \State $\chi=0$;
    \For {$\ell$ from 2 to $L-1$}
    \State $\chi=\chi+\zeta\sps{\ell}$;
    \State $\chi=\LCI[\linear; 2^\ell, r, 2r](\chi)$;
    \State $\chi=\Reshape[r, 2^{\ell+1}](\chi)$;
    \EndFor
    \State $\chi=\chi+\zeta\sps{L}$;
    \State $\chi=\LCI[\linear; 2^L, r, m](\chi)$;
    \State $\chi=\Reshape[1,N](\chi)$;
    \State $u = \chi + u\sps{\ad}$;
\end{algorithmic}
\end{minipage}
\end{minipage}
\end{center}
\end{small}
\caption{Application of NN architecture for $\cH^2$-matrices on a vector $v\in\bbR^N$.}%
\label{alg:h2NN}%
\end{algorithm}

Combining \eqref{eq:LCV}, \eqref{eq:LCCBU} and \eqref{eq:LCMA} and adding necessary $\Reshape$, we
can now translate \cref{alg:H2} to a neural network representation of the matrix-vector
multiplication of $\cH^2$-matrices in \cref{alg:h2NN}, which is illustrated in \cref{fig:h2nn}.

Let us now calculate the number of parameters used in the network in \cref{alg:h2NN}.  For
simplicity, we ignore the number of parameters in the bias terms $b$ and only consider the ones in
the weight matrices $W$. Given that the number of parameters in an LC layer is $N_x'\alpha\alpha'w$,
the number of parameters for each type of network is:
\begin{equation}
    N_p^{\LCR} = N_x\alpha\alpha',\quad
    N_p^{\LCK} = N_x\alpha\alpha'w,\quad
    N_p^{\LCI} = N_x\alpha\alpha',
\end{equation}
Then the total number of parameters in \cref{alg:h2NN} is
\begin{equation}
    \begin{aligned}
    N_p^{\cH^2} &= 2^Lm^2(2n_{b}\sps{\ad}+1)
    + Nr+2\sum_{\ell=2}^{L-1}2^{\ell+1}r^2
    + \sum_{\ell=2}^{L}2^{\ell}r^2 (2n_b\sps{\ell}+1)
    + 2^Lrm\\
    &\leq Nm(2n_b+1) + 2Nr+2Nr(2n_b+3)\\
    &\leq 3Nm(2n_b+3) = O(N),
    \end{aligned}
\end{equation}
where $n_b=\max(n_{b}\sps{\ad},n\sps{\ell}_b)$, $r\leq m$ and $2^Lm=N$ are used. The calculation
shows that the number of parameters in the neural network scales linearly in $N$ and is therefore of
the same order as the memory storage in $\cH^2$-matrices. This is lower than the quasilinear order
$O(N\log(N))$ of $\cH$-matrices and its neural network generalization.

layers, which reduces the computational and storage cost for large systems. Numerical results


\subsection{Multi-dimensional case}\label{sec:nD}
Following the discussion in the previous section, \cref{alg:h2NN} can be easily extended to the
$d$-dimensional case by performing a tensor-product of the one-dimensional case. In this subsection,
we consider $d=2$ for instance, and the generalization to the $d$-dimensional case becomes
straightforward.  For the integral equation
\begin{equation}\label{eq:ie2d}
    u(x) = \int_{\Omega}g(x,y)v(y)\dd y,\quad \Omega=[0,1)\times [0,1),
\end{equation}
we discretize it with an uniform grid with $N\times N$, $N=2^Lm$, grid points 
and denote the resulting matrix obtained from the discretization of \eqref{eq:ie2d} by $A$.
Conceptually \cref{alg:h2NN} required the following $3$ components:
\begin{enumerate}
    \item multiscale decomposition of the matrix $A$, given by \eqref{eq:decompose};
    \item nested low-rank approximation of the far-field blocks of $A$, given by \eqref{eq:coeff_mat} and
        \cref{pro:matrices} for the resulting matrices; 
    \item definition of LC layers and theirs relationship \eqref{eq:LCV},\eqref{eq:LCCBU} and \eqref{eq:LCMA} with
        the matrices in \cref{pro:matrices}.  
\end{enumerate}
We briefly explain how each step can be seamlessly extended to the higher dimension in what follows.
 
\paragraph{Multiscale decomposition.}
The grid is hierarchically partitioned into $L+1$ levels, in which each \emph{box} is defined by
$\cI\sps{d,\ell}_{i}=\cI\sps{\ell}_{i_1}\otimes\cI\sps{\ell}_{i_2}$, where $i=(i_1,i_2)$ is a
multi-dimensional index, $\cI\sps{\ell}_{i_1}$ identifies the segments for 1D case and $\otimes$ is
the tensor product. The definitions of the children list, parent, neighbor list and
interaction list can be easily extended. Each box $\cI$ with $\ell<L$ has $4$ children.  Similarly,
the decomposition \eqref{eq:decompose} on $A$ can also be extended.

\paragraph{Nested low-rank approximation.}
Following the structure of $\cH^2$-matrices, the nonzero blocks of $A\sps{\ell}$ can be approximated
by
\begin{equation}
  A\sps{\ell}_{\cI, \cJ}\approx U\sps{\ell}_{\cI}M\sps{\ell}_{\cI,\cJ}(V\sps{\ell}_{\cJ})^{T},
    \quad U\sps{\ell}_{\cI}, V\sps{\ell}_{\cJ}\in\bbR^{(N / 2^{\ell})^2\times r},\quad
    M\sps{\ell}_{\cI,\cJ}\in\bbR^{r\times r},
\end{equation}
and the matrices $U\sps{\ell}$ satisfy the consistency condition, \ie
\begin{equation} \label{eq:U2d}
 U\sps{\ell}_{\cI}\approx
 \begin{pmatrix}
     U\sps{\ell+1}_{\cJ_1} & & &\\
                           & U\sps{\ell+1}_{\cJ_2} & & \\
                           & & U\sps{\ell+1}_{\cJ_3} & \\
                           & & & U\sps{\ell+1}_{\cJ_4} \\
 \end{pmatrix}
 \begin{pmatrix}
     B\sps{\ell}_{\cJ_1}\\[2mm]
     B\sps{\ell}_{\cJ_2}\\[2mm]
     B\sps{\ell}_{\cJ_3}\\[2mm]
     B\sps{\ell}_{\cJ_4}
 \end{pmatrix},
\end{equation}
where $\cJ_j$ are children of $\cI$, and $B\sps{\ell}_{\cJ_j}\in\bbR^{r\times r}$, $j=1,\dots,4$.
Similarly, the matrices $V\sps{\ell}$ also have the same nested relationship.

We denote an entry of a tensor $T$ by $T_{i,j}$, where $i$ is $2$-dimensional index $i=(i_1,i_2)$.
Using the tensor notations, $U\sps{L}$ and $V\sps{L}$ in \eqref{eq:H2_Av} can be treated as
4-tensors of dimension $N\times N \times 2^L r\times 2^L$, while $B\sps{\ell}$ and $C\sps{\ell}$ in
\eqref{eq:H2_Av} can be treated as 4-tensors of dimension $2^{\ell+1}r\times 2^{\ell+1}\times
2^{\ell}r\times 2^{\ell}$. We generalize the notion of band matrix $A$ to {\em band tensors} $T$ by satisfying  
\begin{equation}
    T_{i,j}=0,\quad \text{if } |i_1-j_1|>n_{b,1} ~\text{ or }~ |i_2-j_2| > n_{b,2},
\end{equation}
where $n_b=(n_{b,1}, n_{b,2})$ is called the band size for tensor. Thus \cref{pro:matrices} can be
extended to
\begin{property}\label{pro:matrices2d}
    The 4-tensors
    \begin{enumerate}
        \item $U\sps{L}$ and $V\sps{L}$ are block diagonal tensors with block size $N/2^L\times N /
          2^L \times r\times 1$.
        \item $B\sps{\ell}$ and $C\sps{\ell}$, $\ell=2,\cdots,L-1$ are block diagonal tensors with
          block size $2r\times 2\times r\times 1$ 
        \item $M\sps{\ell}$, $\ell=2,\cdots,L$ are block cyclic band tensors with block size
            $r\times 1\times r\times 1$ and band size $n\sps{\ell}_b$, which is $(2,2)$ for $\ell=2$
            and $(3,3)$ for $\ell>2$;
        \item $A\sps{ad}$ is a block cyclic band matrix with block size $m\times m\times m\times m$ and band size
            $n_b\sps{\ad}=(1,1)$.
    \end{enumerate}
\end{property}

\paragraph{LC layers.}
An NN layer for 2D can be represented by a 3-tensor of size $\alpha\times N_{x,1}\times N_{x,2}$,
where $\alpha$ is the channel dimension and $N_{x,1}$, $N_{x,2}$ are the spatial dimensions. If a
layer $\xi$ with size $\alpha\times N_{x,1}\times N_{x,2}$ is connected to a locally connected layer
$\zeta$ with size $\alpha'\times N_{x,1}'\times N_{x,2}'$, then
\begin{footnotesize}
\begin{equation}\label{eq:lc2d}
    \zeta_{c',i} = \phi\left( \sum_{j=(i-1)s+1}^{(i-1)s+w}\sum_{c=1}^{\alpha}
    W_{c',c;i,j}\xi_{c,j} + b_{c',i}\right),
    \quad i_1=1,\dots,N_{x,1}', i_2=1,\dots,N_{x,2}',~ c'=1,\dots,\alpha',
\end{equation}
\end{footnotesize}%
where $(i-1)s= ( (i_1-1)s_1, (i_2-1)s_2)$. As in the 1D case, the channel dimension corresponds to
the rank $r$, and the spatial dimensions correspond to the grid points of the discretized domain.
Analogously to the 1D case, we define the LC networks $\LCR$, $\LCK$ and $\LCI$ and use them to
express the 6 operations in \eqref{eq:hm} that constitute the building blocks of the neural network.
The parameters $N_x$, $s$ and $w$ in the one-dimensional LC networks are replaced by their
2-dimensional counterpart $N_x=(N_{x,1}, N_{x,2})$, $s=(s_1,s_2)$ and $w=(w_{1},w_{2})$,
respectively.  We point out that $s=w=\frac{N_x}{N_x'}$ for the 1D case is replaced by
$s_j=w_{j}=\frac{N_{x,j}}{N_{x,j}'}$, $j=1,2$ for the 2D case in the definition of LC.

\begin{figure}[ht]
    \centering
    \subfloat[{Diagram of $\ReshapeT[2, 1,3,3]$ and $\ReshapeM[2,1,3,3]$}%
    ]{
        \begin{overpic}[width=0.45\textwidth,page=1]{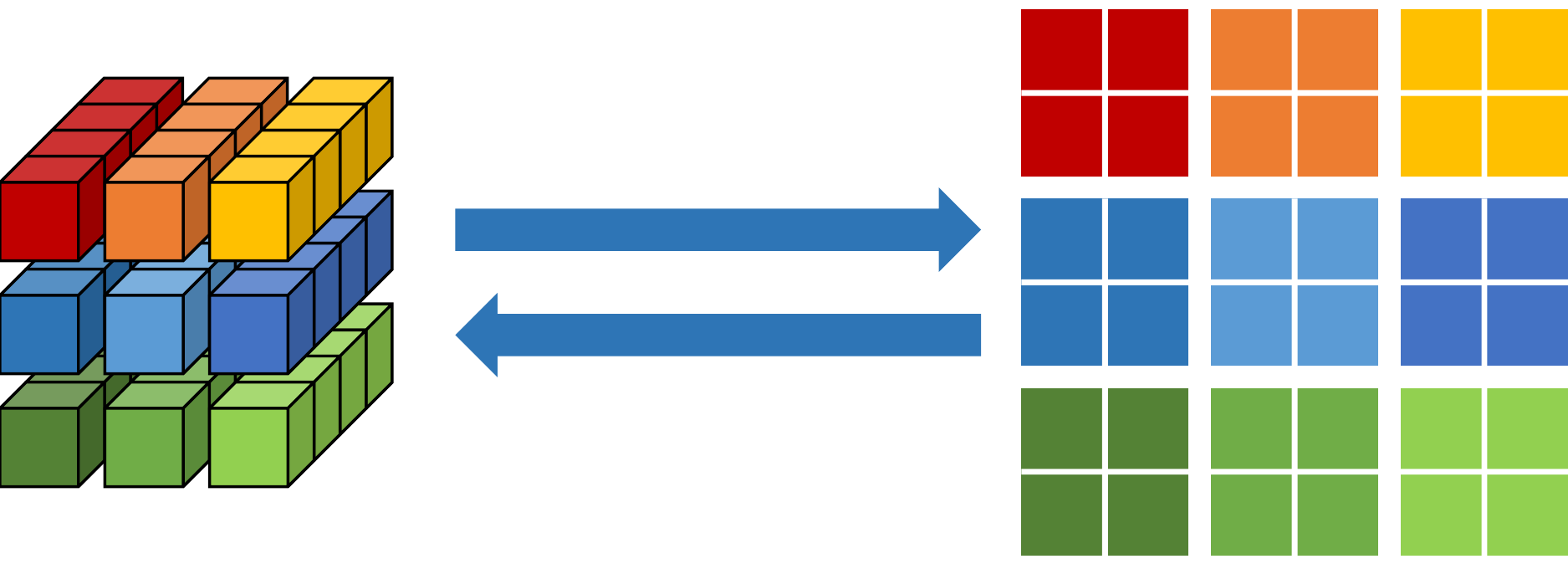}
            \put(36, 8){\small$\ReshapeT$}
            \put(36,25){\small$\ReshapeM$}
        \end{overpic}
    }
    \quad {\color{gray}\rule{1pt}{0.15\textwidth}}\quad
    \subfloat[{Diagram of $\ReshapeM[2,2,3,3]$}]{
        \begin{overpic}[width=0.45\textwidth,page=2]{reshape2d.pdf}
            \put(38,22){\small$\ReshapeM$}
        \end{overpic}
    }
    \caption{\label{fig:reshape2d}Diagram of $\ReshapeT$ and $\ReshapeM$ in \cref{alg:h2NN2d}.}
\end{figure}
Using the notations above we extend \cref{alg:h2NN} to the 2D case in \cref{alg:h2NN2d}. It is crucial to
note that the $\ReshapeT$ and $\ReshapeM$ functions in \cref{alg:h2NN2d} 
are \textit{not} the usual column major based reshaping operations. $\ReshapeM[a, r, n_1, n_2]$
reshapes a 3-tensor $T$ with size $a^2r\times n_1\times n_2$ to a 3-tensor with size $r\times an_1
\times an_2$ by reshaping each row $T_{\cdot, j,k}$ to a 3-tensor with size $r\times a\times a$ and 
joining them to a large 3-tensor. $\ReshapeT[a,r,n_1,n_2]$ is the inverse of $\ReshapeM[a,r,n_1,n_2]$. 
\cref{fig:reshape2d} 
diagrams these two reshape functions.

\begin{algorithm}[htb]
\begin{footnotesize}
\begin{center}
\begin{minipage}{\textwidth}
\begin{minipage}{0.59\textwidth}
\begin{algorithmic}[1]
    \State $\tilde{v}=\ReshapeT[m,1,2^L,2^L](v)$;
    \State $\tilde{u}\sps{\ad}=\LCK[\linear; (2^L,2^L),m^2,m^2,2n_{b}\sps{\ad}+1](\tilde{v})$;
    \State $u\sps{\ad}=\ReshapeM[m,1,2^L,2^L](\tilde{u}\sps{\ad})$;
    \State $\xi\sps{L}=\LCR[\linear; (N,N), 1, (2^L,2^L), r](v)$;
    \For {$\ell$ from $L-1$ to 2 by $-1$}
    \State $\xi\sps{\ell}=\LCR[\linear; (2^{\ell+1},2^{\ell+1}), r, (2^\ell,2^\ell), r](\xi\sps{\ell+1})$;
    \EndFor
    \For {$\ell$ from $2$ to $L$}
    \State $\zeta\sps{\ell}=\LCK[\linear; (2^\ell,2^\ell), r, r, 2n_b\sps{\ell}+1](\xi\sps{\ell})$;
    \EndFor
    \algstore{break2}
\end{algorithmic}
\end{minipage}
\begin{minipage}{0.40\textwidth}
\begin{algorithmic}[1]
    \algrestore{break2}
    \State $\chi=0$;
    \For {$\ell$ from 2 to $L-1$}
    \State $\chi=\chi+\zeta\sps{\ell}$;
    \State $\chi=\LCI[\linear; (2^\ell,2^\ell), r, 2r](\chi)$;
    \State $\chi=\ReshapeM[2,r, 2^{\ell}, 2^{\ell}](\chi)$;
    \EndFor
    \State $\chi=\chi+\zeta\sps{L}$;
    \State $\chi=\LCI[\linear; (2^L,2^L), r, m^2](\chi)$;
    \State $\chi=\ReshapeM[m,1,2^L,2^L](\chi)$;
    \State $u = \chi + u\sps{\ad}$;
\end{algorithmic}
\end{minipage}
\end{minipage}
\end{center}
\end{footnotesize}
\caption{Application of NN architecture for $\cH^2$-matrices on a vector $v\in\bbR^{N^2}$.}%
\label{alg:h2NN2d}%
\end{algorithm}

\section{Multiscale neural network}\label{sec:mnn}
The nonlinear map in the form
$u=\cM(v)$ with $u,v\in\bbR^{N^d}$\add{, which can be viewed as a nonlinear generalization
of pseudo-differential operators,} is ubiquitous from integral equations and partial differential
equations in practical applications. In general, to evaluate such nonlinear maps, one needs to 
use iterative methods that may require a large number of iterations, and at each iteration one may 
need to solve the underlying equation several times, resulting in computational expensive algorithms.
Instead, we propose to bypass this endeavor by leveraging the
ability of NNs to represent high-dimensional nonlinear maps. In this
section, we construct a hierarchical approximation of such a nonlinear map by extending the architectures 
provided in \cref{alg:h2NN} and \cref{alg:h2NN2d} to the nonlinear case. We refer to the resulting 
NN architecture as {\em multiscale neural network-$\cH^2$} (MNN-$\cH^2$) due to its multiscale structure 
inspired by $\cH^2$-matrices.

To simplify the notation, we focus on the 1D case in this section. The following presentation can be
readily extended to the multi-dimensional case by following the discussion in \cref{sec:nD}.


\subsection{Algorithm and architecture}\label{sec:mnnAlgorithm}
\begin{algorithm}[htb]
\begin{small}
\begin{center}
\begin{minipage}{0.99\textwidth}
\begin{minipage}{0.6\textwidth}
\begin{algorithmic}[1]
    \State $\xi_0=\Reshape[m,2^L](v)$;
    \For {$k$ from $1$ to $K$ do}
    \State $\xi_k=\LCK[\phi; 2^L,m,m,2n_{b}\sps{\ad}+1](\xi_{k-1})$; \label{alg:line_ad}
    \EndFor
    \State $u\sps{\ad}=\Reshape[1,N](\xi_K)$;
    \State $\zeta_0\sps{L}=\LCR[\linear; N, 1, 2^L, r](v)$;
    \For {$\ell$ from $L-1$ to 2 by $-1$}
    \State $\zeta_0\sps{\ell}=\LCR[\phi; 2^{\ell+1}, r, 2^\ell, r](\zeta_0\sps{\ell+1})$;
    \EndFor
    \For {$\ell$ from $2$ to $L$}
    \For {$k$ from $1$ to $K$}
    \State $\zeta_k\sps{\ell}=\LCK[\phi; 2^\ell, r, r, 2n_b\sps{\ell}+1](\zeta_{k-1}\sps{\ell})$;
    \EndFor
    \EndFor
    \algstore{break3}
\end{algorithmic}
\end{minipage}
\begin{minipage}{0.35\textwidth}
\begin{algorithmic}[1]
    \algrestore{break3}
    \State $\chi=0$;
    \For {$\ell$ from 2 to $L-1$}
    \State $\chi=\chi+\zeta_K\sps{\ell}$;
    \State $\chi=\LCI[\phi; 2^\ell, r, 2r](\chi)$;
    \State $\chi=\Reshape[r, 2^{\ell+1}](\chi)$;
    \EndFor
    \State $\chi=\chi+\zeta_K\sps{L}$;
    \State $\chi=\LCI[\linear; 2^L, r, m](\chi)$;
    \State $\chi=\Reshape[1,N](\chi)$;
    \State $u = \chi + u\sps{\ad}$;
\end{algorithmic}
\vspace{15mm}
\end{minipage}
\end{minipage}
\end{center}
\end{small}
\caption{Application of MNN-$\cH^2$ to a vector $v\in\bbR^N$.}
\label{alg:mnnh2}%
\end{algorithm}

\begin{figure}[ht]
    \centering
    \begin{overpic}[width=0.95\textwidth]{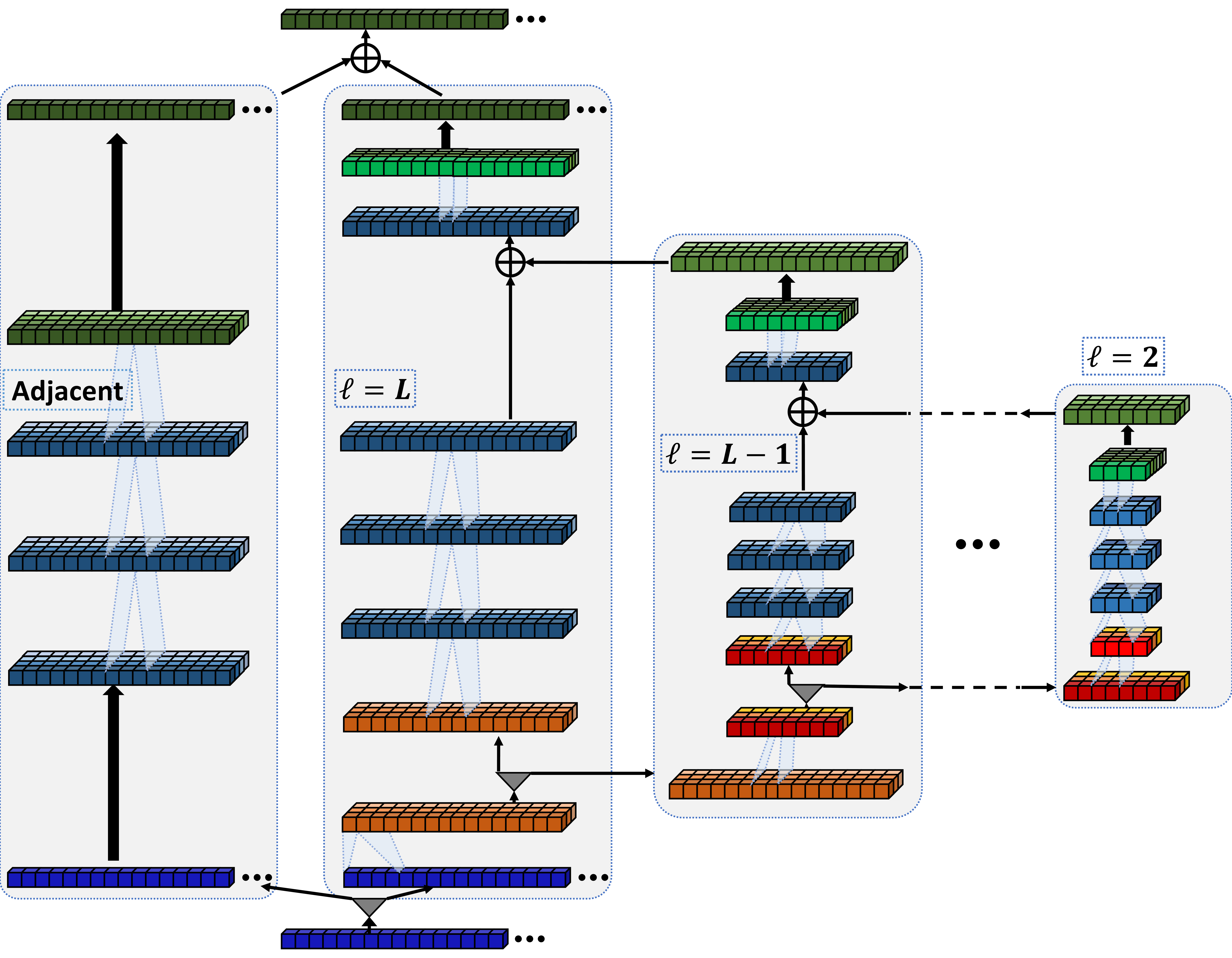}
        \put(11.0, 14.0){\small$\Reshape$}
        \put(11.0, 60.0){\small$\Reshape$}
        \put(13.5, 27.0){\footnotesize$\LCK$-$\phi$}
        \put(13.5, 37.0){\footnotesize$\LCK$-$\phi$}
        \put(11.5, 46.0){\footnotesize$\LCK$-$\linear$}
        \put(33.0,  3.0){\footnotesize$\Replicate$}
        \put(32.0, 73.5){\footnotesize$\mathsf{sum}$}
        \put(38.5, 66.2){\footnotesize$\Reshape$}
        \put(41.0, 15.5){\footnotesize$\Replicate$}
        \put(43.8, 57.1){\footnotesize$\mathsf{sum}$}
        \put(35.0,  8.0){\scriptsize$\LCR$-$\linear$}
        \put(40.5, 23.0){\scriptsize$\LCK$-$\phi$}
        \put(40.5, 31.0){\scriptsize$\LCK$-$\phi$}
        \put(40.5, 38.5){\scriptsize$\LCK$-$\phi$}
        \put(39.0, 61.5){\scriptsize$\LCI$-$\linear$}
        \put(66.0, 22.5){\tiny$\Replicate$}
        \put(66.5, 45.5){\footnotesize$\mathsf{sum}$}
        \put(65.0, 53.9){\scriptsize$\Reshape$}
        \put(66.0, 16.0){\tiny$\LCR$-$\phi$}
        \put(70.0, 26.2){\tiny$\LCK$-$\phi$}
        \put(70.0, 30.5){\tiny$\LCK$-$\phi$}
        \put(70.0, 34.8){\tiny$\LCK$-$\phi$}
        \put(70.0, 50.0){\tiny$\LCI$-$\phi$}
        \put(94.0, 24.0){\tiny$\LCR$-$\phi$}
        \put(94.0, 27.5){\tiny$\LCK$-$\phi$}
        \put(94.0, 30.8){\tiny$\LCK$-$\phi$}
        \put(94.0, 34.6){\tiny$\LCK$-$\phi$}
        \put(94.0, 38.3){\tiny$\LCI$-$\phi$}
        \put(93.0, 41.9){\tiny$\Reshape$}
    \end{overpic}
    \caption{\label{fig:mnnh2}Neural network architecture for MNN-$\cH^2$.}
\end{figure}
 


Similar to \cite{fan2018mnn}, we extend \cref{alg:h2NN} to the nonlinear case by replacing the
linear activation function by a nonlinear one, and extend one $\LCK$ layer to $K\in\bbN$ nonlinear
$\LCK$ layers. \cref{alg:h2NN} is then revised to \cref{alg:mnnh2}.
Following \cite{fan2018mnn}, the last layer corresponding to the adjacent part, the layer
corresponding to $(V\sps{L})^Tv$ and $U\sps{L}\zeta$ are set to linear layers.  In addition, the
layer in line 3 of \cref{alg:mnnh2} is a linear layer when $k=K$, and the activation $\phi$ in
\cref{alg:mnnh2} can be any nonlinear or linear activation function depending of the target
application.  \cref{fig:mnnh2} illustrates the architecture of MNN-$\cH^2$.

Similarly to the linear case, we compute the number of parameters of MNN-$\cH^2$ to obtain
\begin{equation}\label{eq:np_mnnh2}
    \begin{aligned}
        N_{p,LC} &= 2^Lm^2K(2n_{b}\sps{\ad}+1)
    + Nr+2\sum_{\ell=2}^{L-1}2^{\ell+1}r^2
    + K\sum_{\ell=2}^{L}2^{\ell}r^2 (2n_b\sps{\ell}+1) 
    + 2^Lrm\\
    &\leq 2Nr+2Nr(2+K(2n_b+1))+NmK(2n_b+1)\\
    &\leq 3NmK(2n_b+3) \sim O(N).
    \end{aligned}
\end{equation}
Here the number of parameters in $b$ from \eqref{eq:lc} is also ignored.  Compared to
$\cH$-matrices, the main saving of the arithmetic complexity of $\cH^2$-matrices is its nested
structure of $U\sps{\ell}$ and $V\sps{\ell}$. Therefore, accordingly compared to MNN-$\cH$ in
\cite{fan2018mnn}, the main saving on the number of parameters of MNN-$\cH^2$ comes from the nested
structure of $\LCR$ and $\LCI$ layers in \cref{alg:mnnh2}. \add{For a system with large $N$, 
MNN-$\cH^2$ has fewer parameters, thus it can reduce the computational cost and storage cost.}



\subsection{Translation invariant case}\label{sec:invariant}
For the linear system \eqref{eq:integral}, if the kernel is of convolution type, \ie $g(x,y) =
g(x-y)$, then the matrix $A$ is a Toeplitz matrix. As a result, the matrices $M\sps{\ell}$,
$A^{(\ad)}$, $U\sps{L}$, $V\sps{L}$, $B\sps{\ell}$ and $C\sps{\ell}$ are all block cyclic matrices.
In the more general nonlinear case, the operator $\cM$ is {\em translation invariant} (or more
accurately translation equivariant) if 
\begin{equation}\label{eq:invariant}
  \cT\cM(v) = \cM(\cT v)
\end{equation}
holds for any translation operator $\cT$. This indicates that the weights $W_{c',c;i,j}$ and bias
$b_{c,i}$ in \eqref{eq:lc} can be independent of index $i$. This is the case of a {\em convolutional
  neural network} (CNN):
\begin{equation}\label{eq:cnn}
    \zeta_{c',i} = \phi\left( \sum_{j=(i-1)s+1}^{(i-1)s+w}\sum_{c=1}^{\alpha} W_{c',c;j}\xi_{c,j} +
    b_{c'}\right), \quad i=1,\dots,N_{x}', ~ c'=1,\dots,\alpha',
\end{equation}
Note that the difference between this and an LC network is that here $W$ and $b$ are independent of
$i$. In this convolutional setting, we shall instead refer to the LC layers $\LCR$, $\LCK$, and
$\LCI$ as $\CNNR$, $\CNNK$, and $\CNNI$, respectively. By replacing the LC layers in
\cref{alg:mnnh2} with the corresponding CNN layers, we obtain the neural network architecture for
the translation invariant kernel. It is easy to calculate that the number of parameters of $\CNNR$,
$\CNNK$ and $\CNNI$ are
\begin{equation}
    N_p^{\CNNR} = \frac{N_x}{N_x'}\alpha',\quad
    N_p^{\CNNK} = \alpha\alpha'w,\quad
    N_p^{\CNNI} = \alpha\alpha'.
\end{equation}
Thus, the number of parameters in \cref{alg:mnnh2} implemented by CNN is $O(\log(N))$ as shown
below:
\begin{equation}
    \begin{aligned}
        N_{p,\mathrm{CNN}} &= \frac{N}{2^L}r
        +2\sum_{\ell=2}^{L-1}2r^2
        + K\sum_{\ell=2}^{L}r^2 (2n_b\sps{\ell}+1) 
        +rm+m^2K(2n_{b}\sps{\ad}+1)\\
        &\leq 2mr+4(L-3)r^2+Kr^2(2n_b+1)(L-2)+Km^2(2n_b+1)\\
        &\leq m^2(4L+K(2n_b+1)(L-1)) = O(\log(N)).
    \end{aligned}
\end{equation}

\paragraph{Mixed model for the non-translation invariant case.}
Note that the number of parameters in the translation invariant case is much lower compared to the
non-invariant case. In addition, the constant $3mK(2n_b+3)$ in \eqref{eq:np_mnnh2} is usually a
large number for practical applications. For example, if $m=5$, $K=5$, $n_b=3$, the constant is
$675$. To reduce the number of parameters in MNN-$\cH^2$, we propose a mixed model to replace some
of the LC layers by CNN layers even in the non-translation invariant setting. For example, in one of
the numerical applications in \cref{sec:application}, we use LC layers for the $\LCR$ and $\LCI$
layers and for the last layer of the adjacent part, while using $\CNNK$ for the remaining layers.
We will verify the effectiveness of this heuristic mixed model in \cref{sec:RTE}.

\section{Applications}\label{sec:application}
In this section we study the performance of the MNN-$\cH^2$ structure using three examples: the
nonlinear Schr{\"o}dinger equation (NLSE) in \cref{sec:NLSE}, the steady-state radiative transfer
equation (RTE) in \cref{sec:RTE}, and the Kohn-Sham map in \cref{sec:KSMap}.

The MNN-$\cH^2$ structure was implemented in Keras \cite{keras}, a high-level neural network
application programming interface (API) running on top of TensorFlow \cite{tensorflow}, which is an
open source software library for high performance numerical computation.  The loss function is
chosen as the mean squared error. The optimization is performed using the Nadam optimizer
\cite{dozat2015incorporating}.  The weights in MNN-$\cH^2$ are initialized randomly from the normal
distribution and the batch size is always set between $1 / 100$th and $1 / 50$th of the number of
training samples.  As discussed in Section \ref{sec:invariant}, if the operator $\cM$ is translation
invariant, all the layers are implemented using CNN layers, otherwise we use LC layers or a mixture
of LC and CNN layers.

In all the tests, the band size is chosen as $n_{b,\ad}=1$ and $n_{b}\sps{\ell}$ is $2$ for $\ell=2$ and
$3$ otherwise. The activation function in $\LCR$ and $\LCI$ is chosen to be linear, while ReLU is
used in $\LCK$. All the tests are run on GPU with data type \verb|float32|.
The selection of parameters $r$ (number of channels), $L$ ($N=2^Lm$) and $K$ (number of layers in
Algorithm \ref{alg:mnnh2}) are problem dependent.

The training and test errors are measured by the relative error with respect to $\ell^2$ norm
\begin{equation}\label{eq:relativeerror}
  \epsilon = \frac{||u-u_{NN}||_{\ell^2}}{||u||_{\ell^2}}.
\end{equation}
where $u$ is the target solution generated by numerical discretization of PDEs and $u_{NN}$ is the
prediction solution by the neural network. We denote by $\trainerror$ and $\testerror$ the average
training error and average test error within a given set of samples, respectively. Similarly, we denote by
$\sigma_{\mathrm{train}}$ and $\sigma_{\mathrm{test}}$ the estimated standard deviation of the
training and test errors within the given set of samples. The numerical results presented in this
section are obtained by repeating the training a few times, using different random seeds. 




\subsection{NLSE with inhomogeneous background potential}\label{sec:NLSE}
The nonlinear Schr{\"o}dinger equation (NLSE) is a widely used model in quantum physics to study
phenomenon such as the Bose-Einstein condensation~\cite{anglin2002bose,pitaevskii1961vortex}. It has
been studied in \cite{fan2018mnn} using the MNN-$\cH$ structure.
In this work, we use the same example to compare the results from MNN-$\cH^2$ with those from MNN-$\cH$.
Here we study the NLSE with inhomogeneous background potential $V(x)$
\begin{equation}\label{eq:nlse}
    \begin{aligned}
        &-\Delta u(x) + V(x)u(x) + \beta u(x)^3=E u(x),\quad x\in [0,1]^d,\\
        &\text{ s.t.} \int_{[0,1]^d}u(x)^2\dd x=1, \text{ and } \int_{[0,1]^d}u(x)\dd x>0,
    \end{aligned}
\end{equation}
with periodic boundary conditions, to find its ground state $u_G(x)$.
We consider a defocusing cubic Schr{\"o}dinger equation with a strong nonlinear term $\beta=10$.
The normalized gradient flow method in \cite{bao2004computing} is employed for the numerical
solution of NLSE.

In this work, we use neural networks to learn the map from the background potential to the ground
state
\begin{equation}
    V(x) \rightarrow u_G(x).
\end{equation}
Clearly, this map is translation invariant, and thus MNN-$\cH^2$ is implemented using CNN rather
than LC network. In the following, we study MNN-$\cH^2$ on 1D and 2D cases, respectively.

In order to compare with MNN-$\cH$ in \cite{fan2018mnn}, we choose the same potential $V$ as in
\cite{fan2018mnn}
\begin{equation}
    V(x) = -\sum_{i=1}^{n_g}\sum_{j_1,\dots,j_d=-\infty}^{\infty}\frac{\rho^{(i)}}{(2\pi T)^{d/2}}\exp\left(
    -\frac{|x-j-c^{(i)}|^2}{2T} \right),
\end{equation}
where the periodic summation imposes periodicity on the potential, and the parameters
$\rho^{(i)}\sim \cU(1,4)$, $c^{(i)}\sim \cU(0,1)^d$, $i=1,\dots,n_g$ and $T\sim \cU(2,4)\times
10^{-3}$.

\subsubsection{One-dimensional case}\label{sec:nlse1d}
\begin{figure}[htb]
  \centering
  \includegraphics[width=0.3\textwidth,clip]{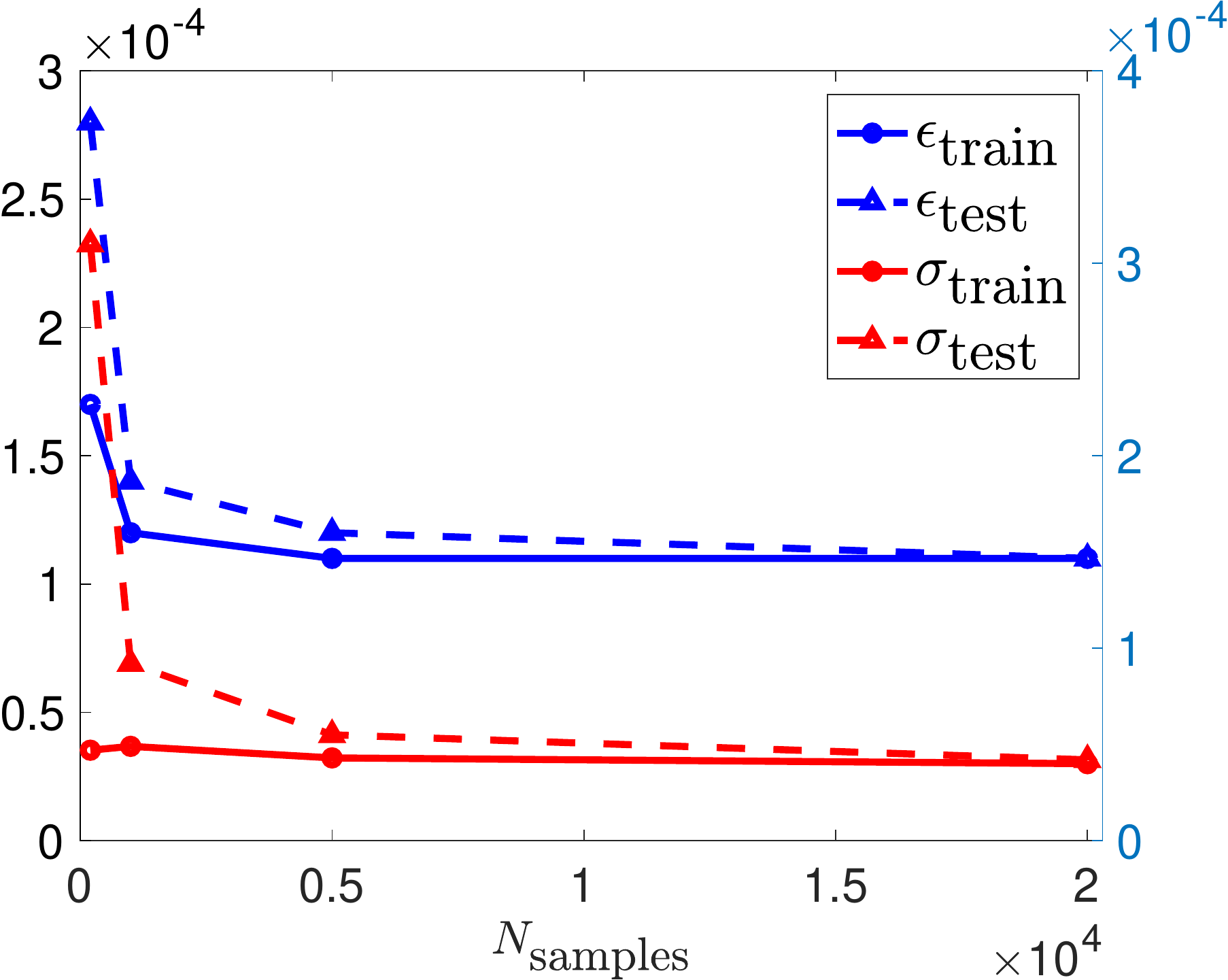}
  \caption{\label{tab:nlse1d_samples}The mean ($\epsilon$ with respect to the left y-axis) and
    standard deviation ($\sigma$ with respect to the right y-axis)  of the relative error in
    approximating the ground state of NLSE for different number of samples {\Ntrainsample} for 1D case
    with $r=6$, $K=5$ and $n_g=2$. In this case, \Nparams=7209.}
\end{figure}

\begin{figure}[htb]
    \centering
    \subfloat[\label{tab:nlse1d_channel} $K=5$, $n_g=2$]{
    \includegraphics[width=0.31\textwidth, clip]{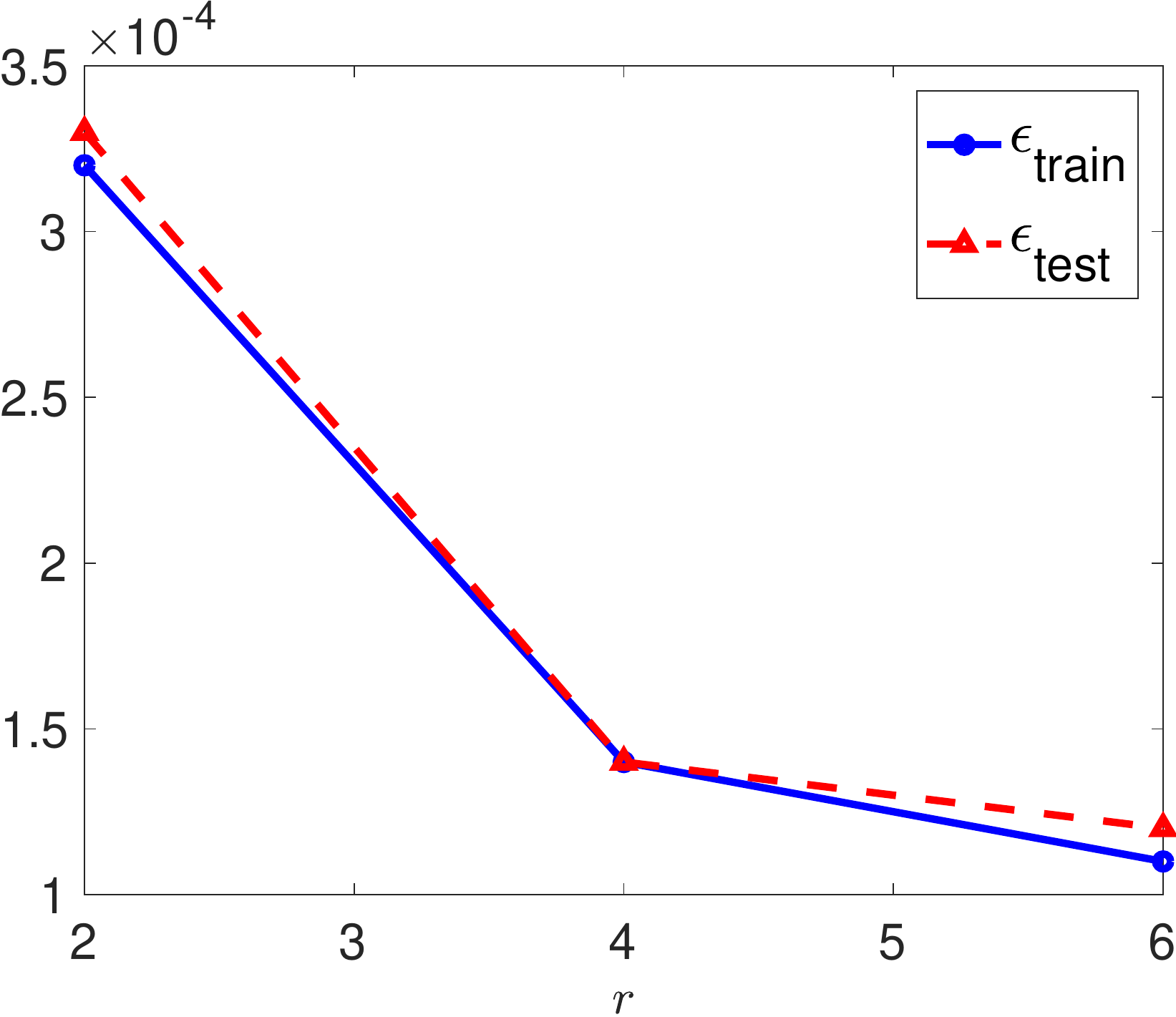}
    }
    \subfloat[\label{tab:nlse1d_layers} $r=6$, $n_g=2$]{
    \includegraphics[width=0.3\textwidth, clip]{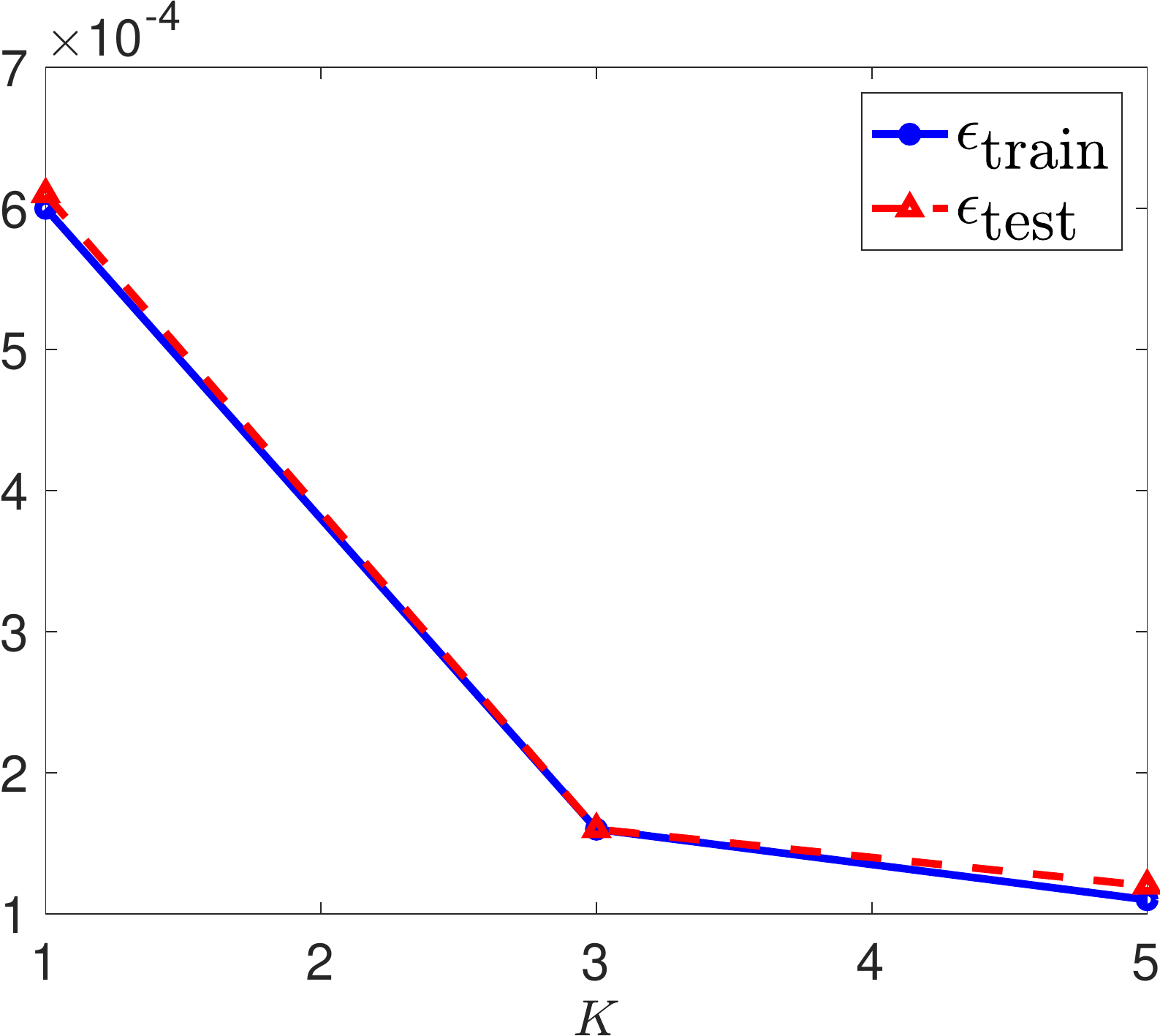}
    }
    \subfloat[\label{tab:nlse1d_ng} $K=5$, $r=6$]{
    \includegraphics[width=0.3\textwidth, clip]{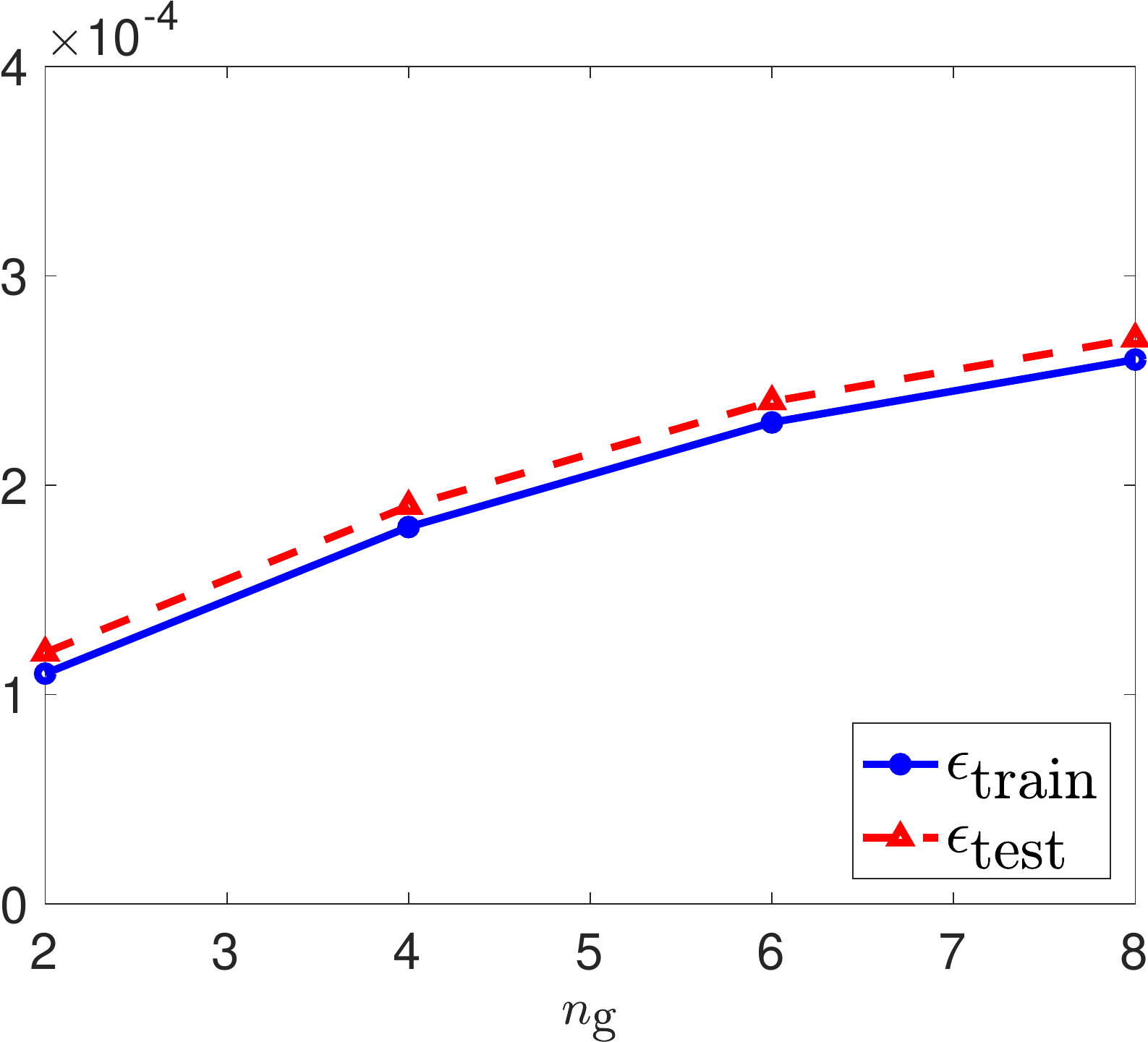}
    }
    \caption{\label{tab:nlse1d}Relative error in approximating the ground state of 1D NLSE for
    different number of channels $r$, different number of $\CNNK$ layers $K$ and different number of
    Gaussians $n_g$ with {\Ntrainsample} = \Ntestsample $=5000$.}
\end{figure}

For the one-dimensional case, we choose the number of discretization points $N=320$, and set $L=7$
and $m=5$. The numerical experiments performed in this section use the same datasets as those in
\cite{fan2018mnn}.  In that context, we study how the performance of MNN-$\cH^2$ depends on the
number of training samples {\Ntrainsample} (\cref{tab:nlse1d_samples}), the number of channels $r$
(\cref{tab:nlse1d_channel}), the number of $\CNNK$ layers $K$ (\cref{tab:nlse1d_layers}), and the
number of Gaussians $n_g$ (\cref{tab:nlse1d_ng}).

\begin{figure}[htb]
    \centering
    \subfloat[test error\label{fig:nlse1d_comE}]{
    \includegraphics[width=0.3\textwidth,clip]{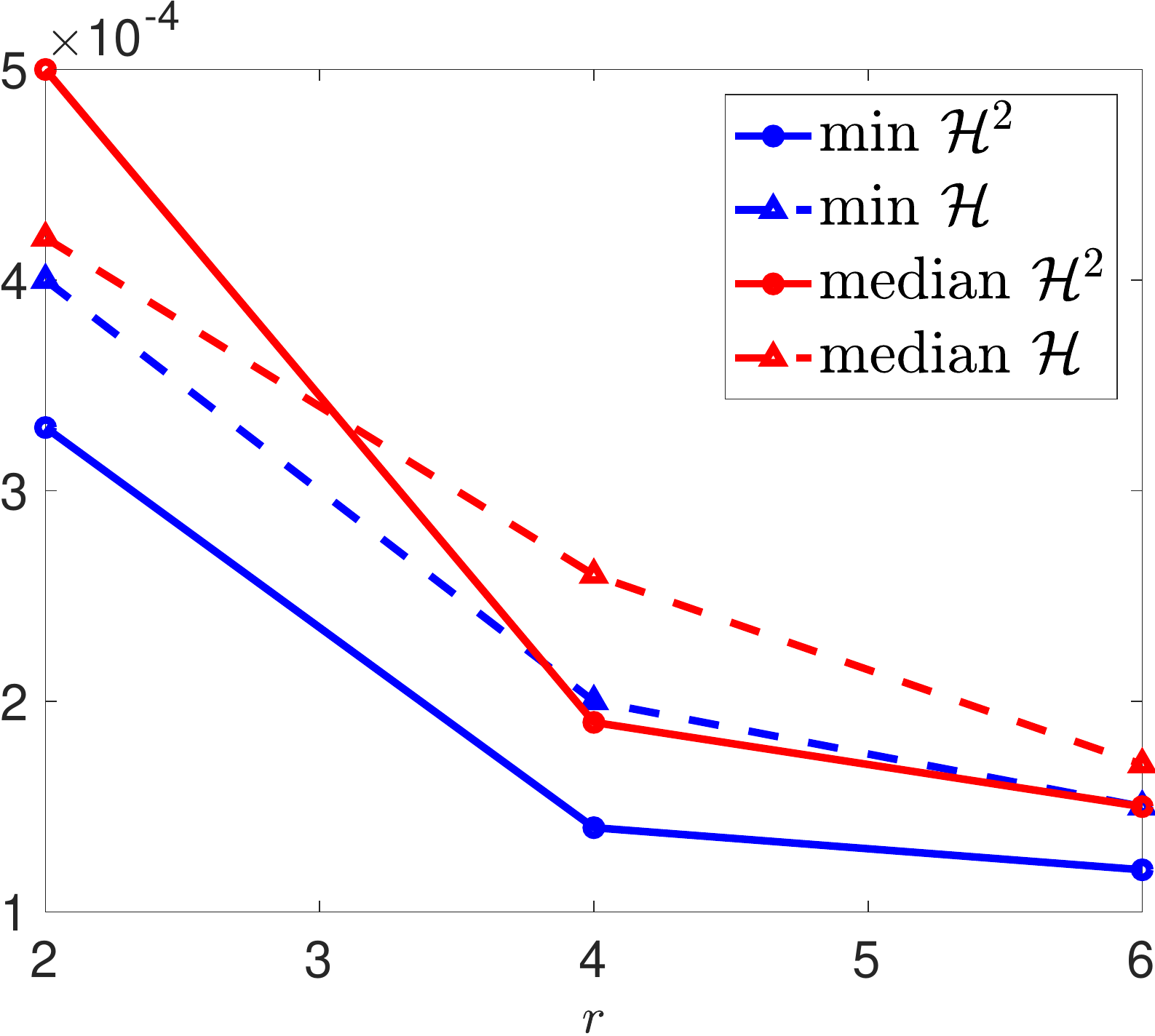}
    }\qquad
    \subfloat[\Nparams\label{fig:nlse1d_comNp}]{
    \includegraphics[width=0.33\textwidth,clip]{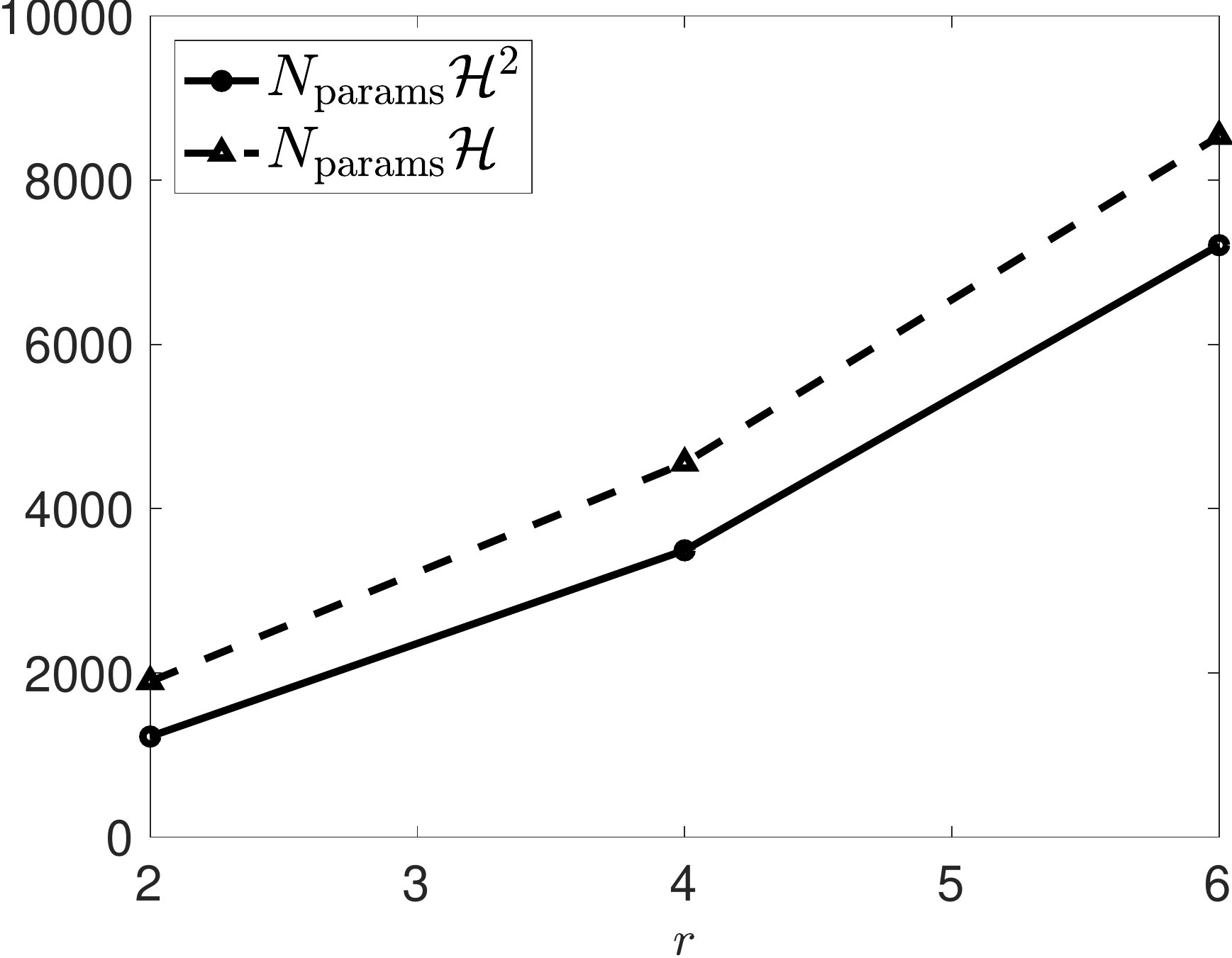}
    }
    \caption{\label{tab:nlse1d_com}
    Numerical results of MNN-$\cH$ / MNN-$\cH^2$ for the minimum and median $\epsilon_{\mathrm{train}}$ for 1D
    NLSE with random initial seed.  The ``min'' and ``median'' stand for the test error
    corresponding to the minimum and median training data cases, respectively, and $H$ and $H^2$
    stand for MNN-$\cH$ and MNN-$\cH^2$, respectively. The setup of MNN-$\cH^2$ is $K=5$, $n_g=2$,
    and {\Ntrainsample} = \Ntestsample $=5000$.}
\end{figure}

\begin{figure}[htb]
    \centering
    \includegraphics[width=1.0\textwidth,clip]{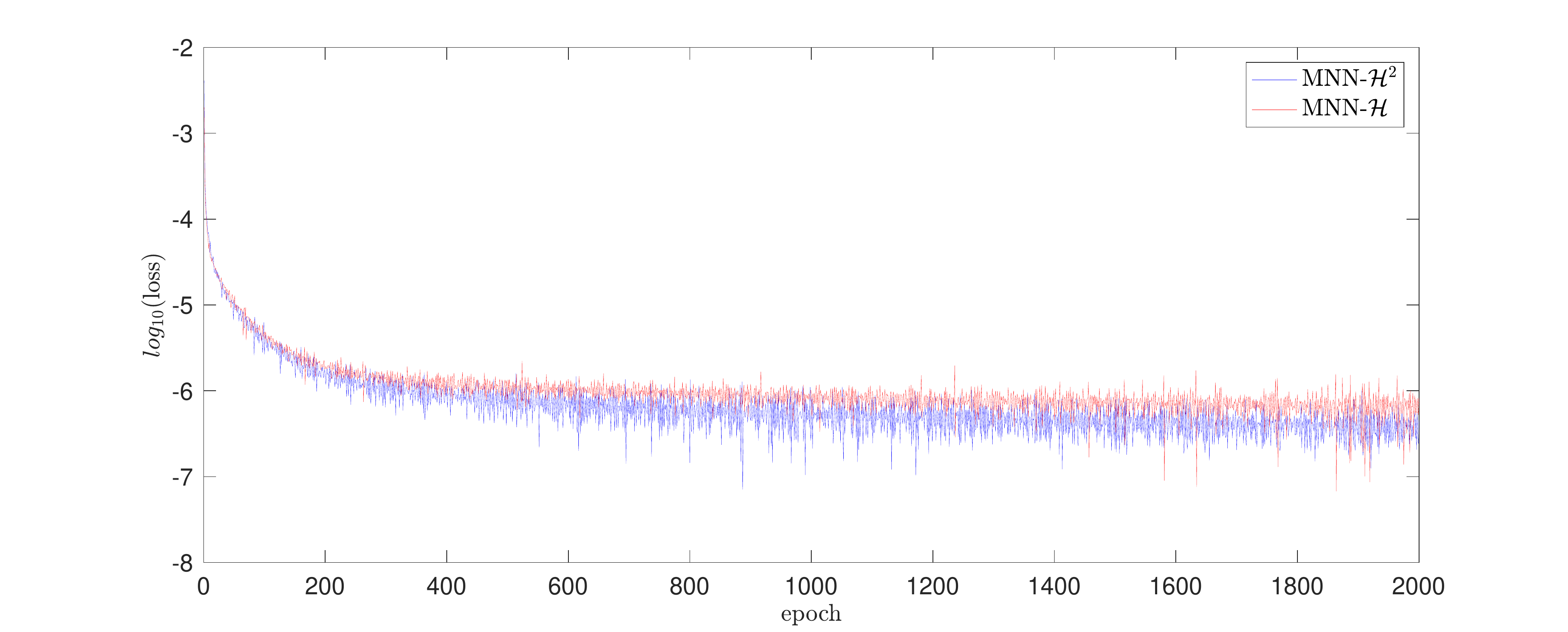}
    \caption{\label{tab:nlse1d_Err}
    Behavior of the loss function of MNN-$\cH^2$ and MNN-$\cH$ with $r=6$ and $K=5$ in the training
    procedure.}
\end{figure}

\cref{tab:nlse1d_samples} shows that MNN-$\cH^2$ can achieve small training error with as few as 200
training samples, which is much smaller than the number of parameters used in the example
(\Nparams=7209). To see why this is possible, let us consider first the linear system $u=Av$ with
$A\in\bbR^{N\times N}$. In order to determine the matrix $A$ using matrix-vector products, we need
at most $O(N)$ independent samples of the form $(u,v)$. Furthermore, if $A$ is an $\cH$-matrix (resp.
$\cH^2$-matrix), the number of parameters in $A$ is reduced to $O(N\log N)$ (resp. $O(N)$).  Hence
only $O(\log(N))$ (resp. $O(1)$) samples of the form $(u,v)$ are sufficient to determine
$A$~\cite{hackbusch1999sparse, hackbusch2000sparse,lin2011fast}.  We expect that similar results can
be generalized to the MNN-$\cH^2$ network, i.e.  the number of samples of the form $(u,v)$ should also be
proportional to the ratio of the number of degrees of freedom in the network and $N$.  For instance,
the neural network used in \cref{tab:nlse1d_samples}, $\frac{N_{\mathrm{params}}}{N} =
\frac{7209}{320}\approx 22.5$, which is much smaller than the number of training samples used in the
simulation.


For the case {\Ntrainsample} $=200$, the test error is slightly larger than the training error, and
the standard deviation within the set of test samples $\sigma_{\mathrm{test}}$ is relatively large. As
{\Ntrainsample} increases to $1000$, the test error is reduced by a factor of $2$, and
$\sigma_{\mathrm{test}}$ is reduced by a factor of $3$. When {\Ntrainsample} increases to $5000$
and $20000$, the test error remains nearly unchanged while $\sigma_{\mathrm{test}}$ continues to
decrease.  For the nonlinear map $u=\cM(v)$, $v\in\Omega\subset\bbR^N$, a large number of samples is
required to obtain an accurate approximation.  Furthermore, we do not observe overfitting in
\cref{tab:nlse1d_samples}.


\cref{tab:nlse1d} presents the numerical results for different choices of channels $r$, $\CNNK$
layers $K$, and Gaussians $n_g$.  As $r$ or $K$ increases,
\cref{tab:nlse1d_channel,tab:nlse1d_layers} show that the error decreases and then stagnates. The
choice of $r=6$ and $K=5$ is used for the 1D NLSE below as a balance of efficiency and accuracy.
In \cref{tab:nlse1d_ng} we find that increasing the number of wells and hence the complexity of the
input field, only leads to marginal increase of the training and test errors.


As demonstrated in the complexity analysis earlier, because of the hierarchical nested bases used in
the restriction and interpolation layers, MNN-$\cH^2$ should use a fewer number of parameters than
MNN-$\cH$ for the same parameter setup, which can be seen in \cref{fig:nlse1d_comNp}.

\cref{fig:nlse1d_comE} compares MNN-$\cH^2$ and MNN-$\cH$ in terms of the minimum error and median
error of the networks by performing the training procedure for a few times with different random
seeds. These results are reported for different number of channels $r$ ranging from $2$ to $6$.
\add{The setup of the training is the same for both MNN-$\cH$ and MNN-$\cH^2$. The learning rate is
  $10^{-3}$, the number of epochs is $6000$, \Ntrainsample=\Ntestsample=$5000$, and the batch size
  is $50$.}  We find that the errors of both networks are comparable for all values of $r$, both in
terms of the minimum and the median. Thus, the reduction of the number of parameters in MNN-$\cH^2$
does not sacrifice accuracy as compared with MNN-$\cH$.  \add{Concerning the training and the test
  procedures, the training and test times of MNN-$\cH^2$ are a bit smaller than those of
  MNN-$\cH$. For example, for the setup with $r=6$ and $K=5$, the training time for MNN-$\cH^2$ on a
  GPU Tesla P100 is $4.5$ hours, while that for MNN-$\cH$ with same setup is $5.5$ hours. The
  prediction time on $10000$ samples is $0.26$ seconds for MNN-$\cH^2$ compared to $0.29$ seconds
  for MNN-$\cH$. The behavior of the loss function for MNN-$\cH^2$ and MNN-$\cH$ during training is
  depicted in \cref{tab:nlse1d_Err}. Here we only present the results for the first 2000 epochs because
  the loss function hardly decreases in the remaining epochs.  One can see that the loss functions
  for both networks exhibit similar behavior and the loss of the MNN-$\cH^2$ is relatively smaller.
  In comparison to MNN-$\cH$, MNN-$\cH^2$ has fewer parameters, trains faster, and yields smaller
  prediction time.} This behavior is also consistently observed in other examples in this section.

\subsubsection{Two-dimensional case}
\begin{figure}[htb]
    \centering
    \subfloat[\label{tab:nlse2d_channel} $K=5$, $n_g=2$]{
    \includegraphics[width=0.3\textwidth, clip]{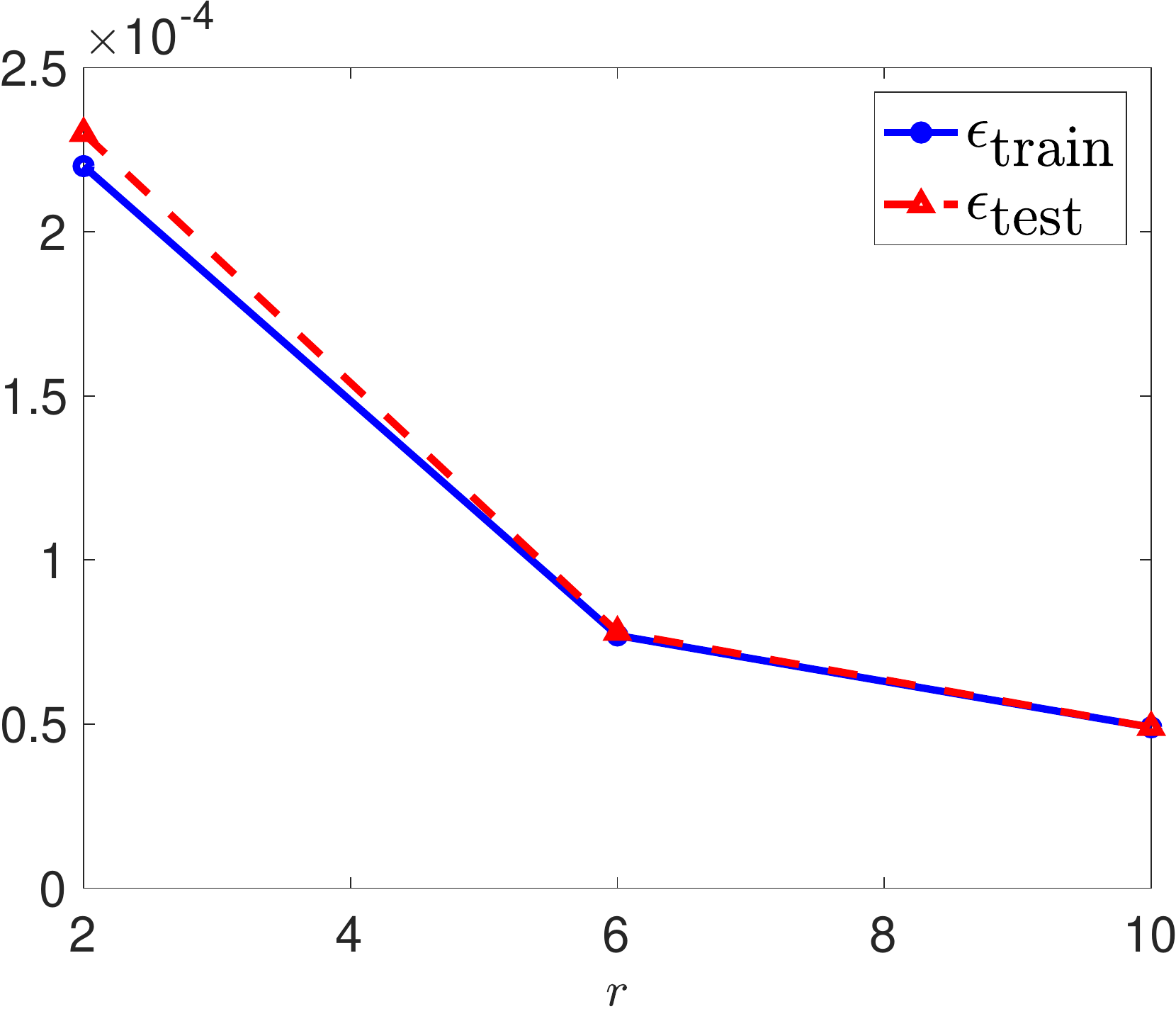}
    }
    \subfloat[\label{tab:nlse2d_layers} $r=6$, $n_g=2$]{
    \includegraphics[width=0.3\textwidth, clip]{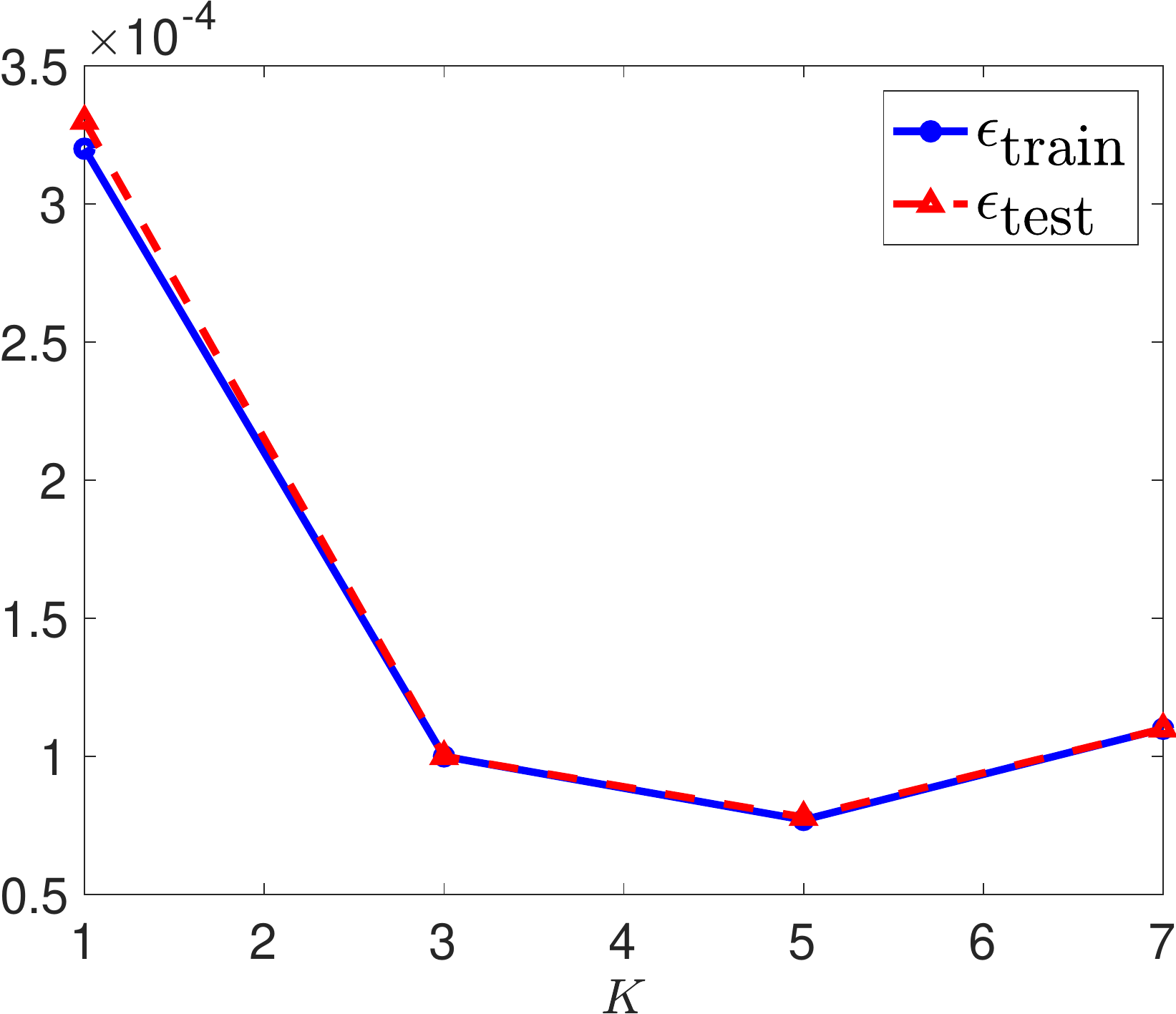}
    }
    \subfloat[\label{tab:nlse2d_ng} $K=5$, $r=6$]{
    \includegraphics[width=0.28\textwidth, clip]{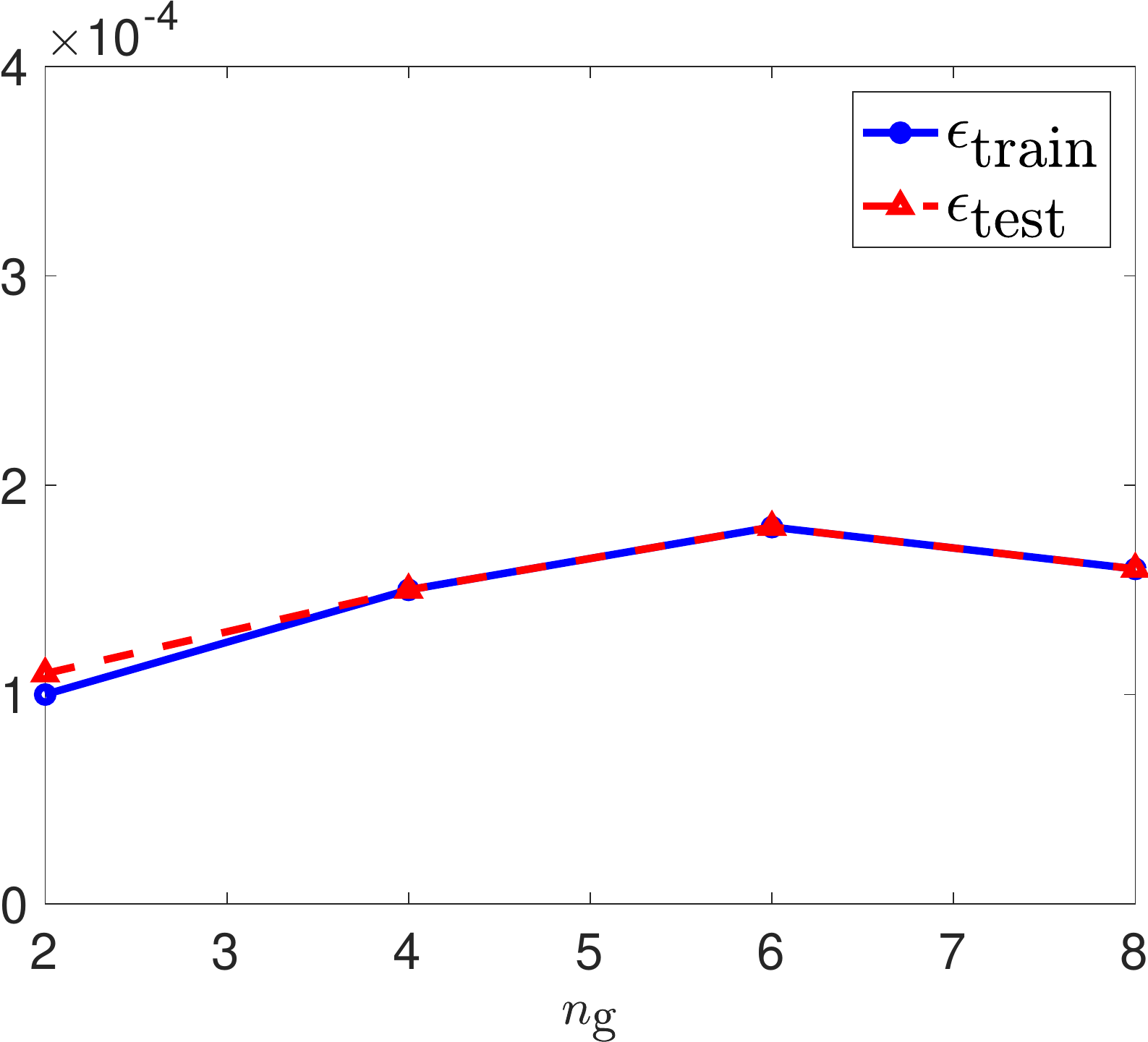}
    }
    \caption{\label{tab:nlse2d}Relative error in approximating the ground state of 2D NLSE for
    different number of channels $r$, different number of $\CNNK$ layers and different number of
    Gaussians $n_g$ for the 2D case with {\Ntrainsample} = \Ntestsample $=20000$.}
\end{figure}

\begin{figure}[htb]
    \centering
    \subfloat[test error]{
    \includegraphics[width=0.3\textwidth,clip]{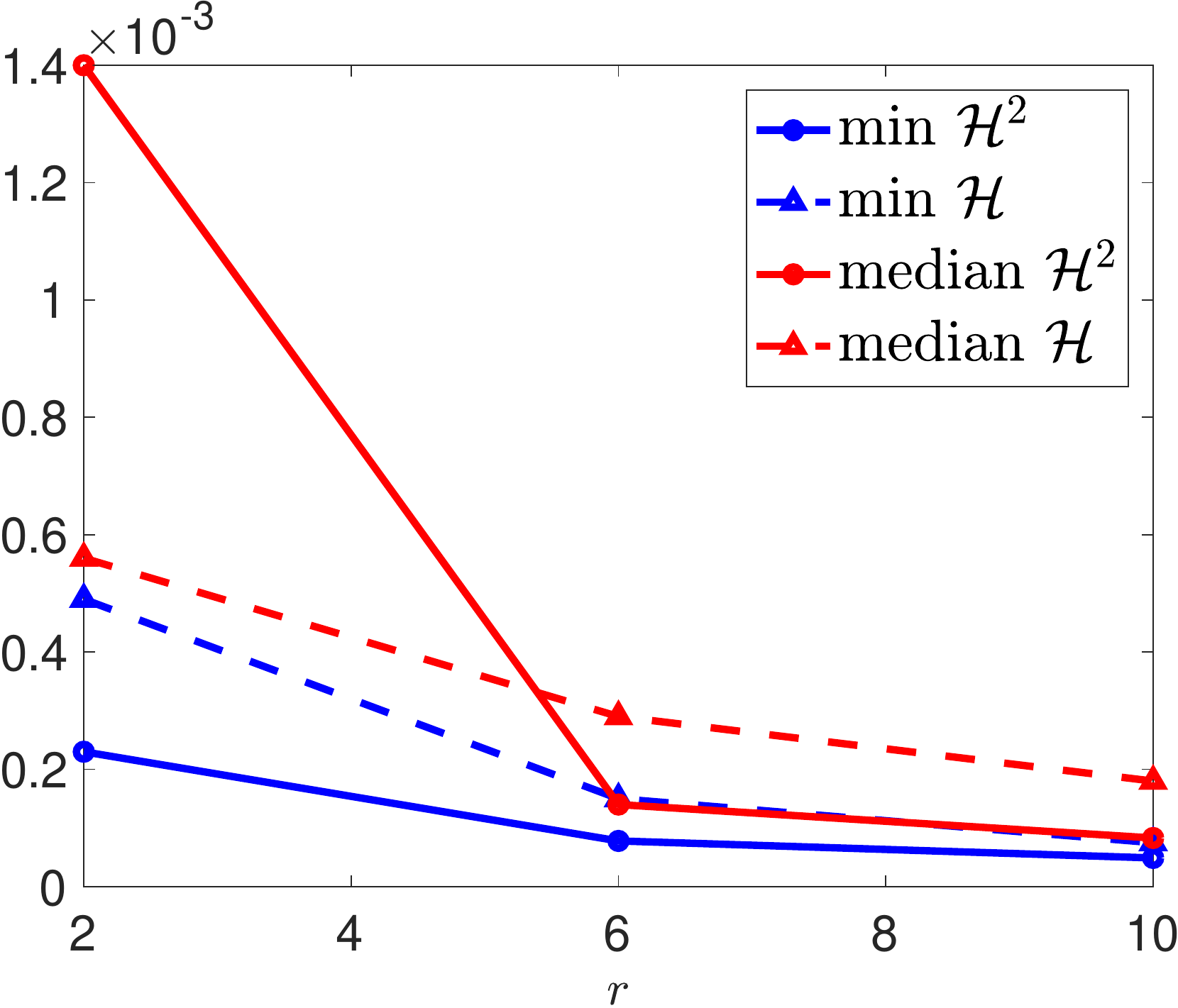}
    }\qquad
    \subfloat[\Nparams]{
    \includegraphics[width=0.3\textwidth,clip]{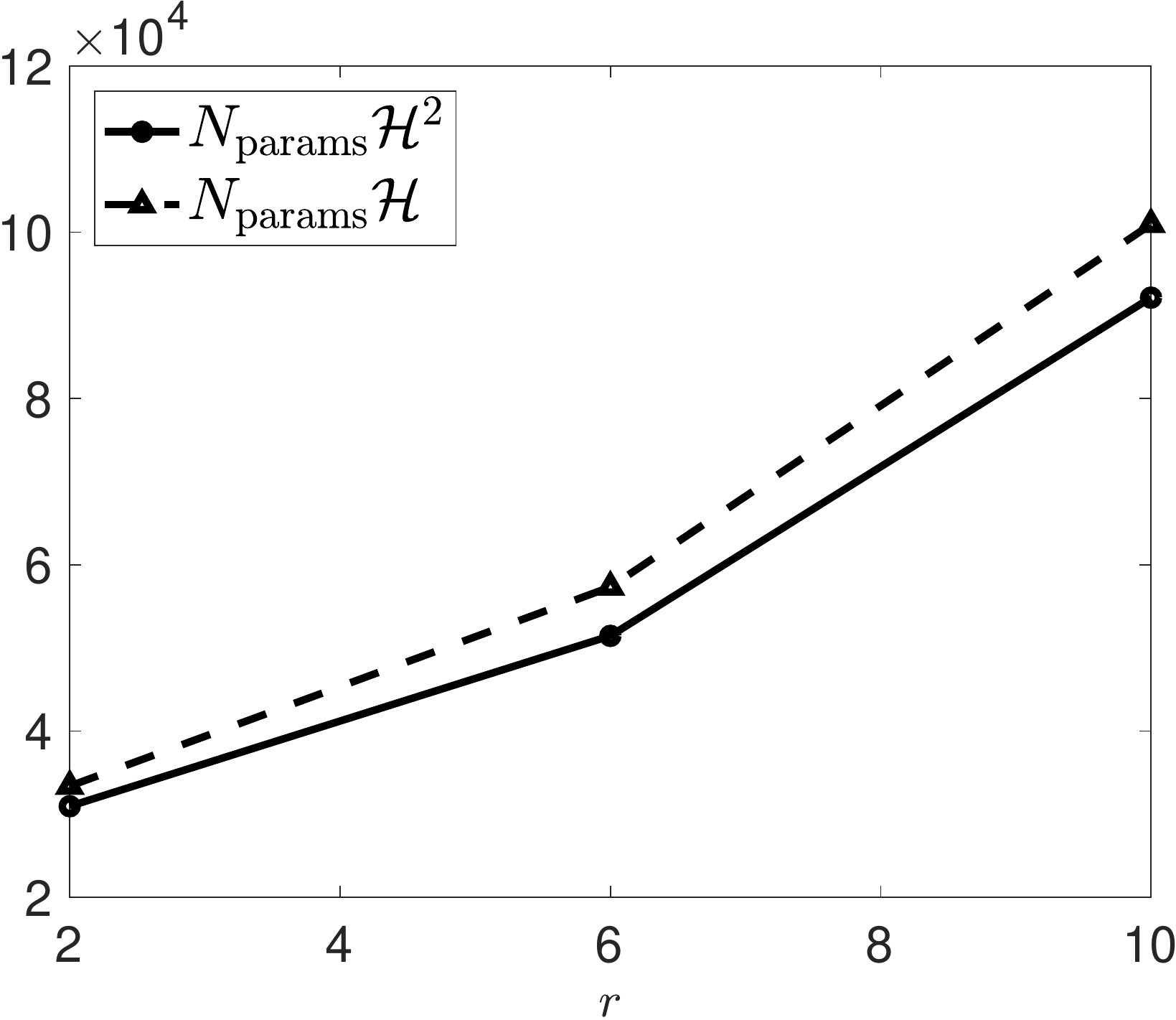}
    }
    \caption{\label{tab:nlse2d_median2}
    Numerical results of MNN-$\cH$ / MNN-$\cH^2$ for the minimum and median $\epsilon_{\mathrm{train}}$ for 1D
    NLSE with random initial seed.  The ``min'' and ``median'' stand for the test error
    corresponding to the minimum and median training data cases, respectively, and $H$ and $H^2$
    stand for MNN-$\cH$ and MNN-$\cH^2$, respectively. The setup of MNN-$\cH^2$ is $K=5$, $n_g=2$
    and {\Ntrainsample} = \Ntestsample $=20000$.}
\end{figure}

For the two-dimensional example, we choose the number of discretization $N$ in each dimension to be
$80$ and set $L=4, m=5$.  The datasets in \cite{fan2018mnn} were used for the $2$D experiments. We
study the behavior of MNN for: different number of channels, $r$ (see \cref{tab:nlse2d_channel} for
the best results and \cref{tab:nlse2d_median2} for the median error); different number of $\CNNK$
layers, $K$ (\cref{tab:nlse2d_layers}); and different number of Gaussians, $n_g$
(\cref{tab:nlse2d_ng}).

Due to the increase of the number of parameters in the 2D networks, we set
{\Ntrainsample}=\Ntestsample = $20000$. From \cref{tab:nlse2d,tab:nlse2d_median2}, we arrive at
similar conclusions as the 1D case: (a) no overfitting is observed for all the tests; (b) the error
first decreases and then stagnates as $r$ or $K$ increases; (c) MNN-$\cH^2$ is not sensitive to the
complexity of the input; and (d) MNN-$\cH^2$ uses fewer number of parameters and obtains a 
comparable error as MNN-$\cH$.

\subsection{Radiative transfer equation}\label{sec:RTE}
Radiative transport equation (RTE) is the widely used tool for describing particle propagation in
many different fields, such as neutron transport in reactor physics \cite{pomraning1973equations},
light transport in atmospheric radiative transfer \cite{marshak20053d}, heat transfer
\cite{koch2004evaluation}, and optical imaging \cite{klose2002optical}. Here we consider the
steady-state RTE in the homogeneous scattering regime
\begin{equation}\label{eq:RTE}
  \begin{aligned}
      v\cdot\nabla_{x}\varphi(x,v)+\mu_t(x)\varphi(x,v)
      &=\mu_s(x) u(x)+f(x),\quad \text{ in } \Omega\times \bbS^{d-1},\quad \Omega\subset\bbR^d,\\
      \varphi(x,v) &= 0,\quad \text{ on } \{(x,v)\in\partial\Omega\times \bbS^{d-1}: n(x)\cdot
      v<0\},\\
      u(x) &= \frac{1}{4\pi}\int_{\bbS^{d-1}}\varphi(x,v)\dd v,
  \end{aligned}
\end{equation}
where $d$ is the dimension, $\varphi(x,v)$ denotes the photon flux that depends on both space $x$
and angle $v$, $f(x)$ is the light source, $\mu_s(x)$ is the scattering coefficient, and $\mu_t(x)$
is the total absorption coefficient. In most applications, one can assume that $\mu_t(x)$ is equal
to $\mu_s(x)$ plus a constant background. The mean density $u(x)$ is uniquely determined by $\mu_s$,
$\mu_t$, and $f$ \cite{fan2018fast}. In this homogeneous regime, by eliminating $\varphi(x,v)$ from
the equation and keeping only $u(x)$ as unknown, one can rewrite RTE as an integral equation
\begin{equation}\label{eq:RTE_IE}
    u = \left( \cI -\cK\mu_s \right)^{-1}\cK f,
\end{equation}
with the operator $\cK$ defined as
\begin{equation}
    \cK f = \int_{y\in\Omega}K(x,y)f(y)\dd y,\quad
    K(x,y) = \frac{\exp\left( -|x-y|\int_0^1\mu_t(x-s(x-y))\dd s \right)}{4\pi |x-y|^{d-1}}.
\end{equation}

In practical applications such as inverse problems, either \eqref{eq:RTE} or \eqref{eq:RTE_IE} is
often solved repetitively, which can be quite expensive even if the fast algorithms for example in
\cite{fan2018fast,ren2016fast} are used. Here, we use MNN-$\cH^2$ to learn the map
\begin{equation}
    \mu_s(x)\rightarrow u(x)
\end{equation}
from the scattering coefficient $\mu_s$ to the mean density $u(x)$.

\subsubsection{One-dimensional slab geometry case}

We first study the one-dimensional slab geometry case for $d=3$, \ie the parameters are homogeneous
on the direction $x_2$ and $x_3$.  With slight abuse of notations, we denote $x_1$ by $x$ in this
subsection. Then, \eqref{eq:RTE_IE} turns to
\begin{equation}\label{eq:RTE_IE1d}
    u(x)=\left( \cI-\cK_1\mu_s \right)^{-1}\cK_1 f(x),
\end{equation}
where the operator $\cK_1$ is defined as
\begin{equation}
    \begin{aligned}
    \cK_1 f(x) &= \int_{y\in \Omega}K_1(x,y)f(x)\dd y,\\
    K_1(x,y) &= \frac{1}{2}\mathrm{Ei}\left(-|x-y|\int_{0}^1\mu_t(x-s(x-y))\dd s\right),
    \end{aligned}
\end{equation}
and $\mathrm{Ei}(\cdot)$ is the exponential integral.

\begin{figure}[htb]
    \centering
    \subfloat[\label{tab:rte1d_channel}$K=5$, $n_g=2$]{
    \includegraphics[width=0.45\textwidth, clip]{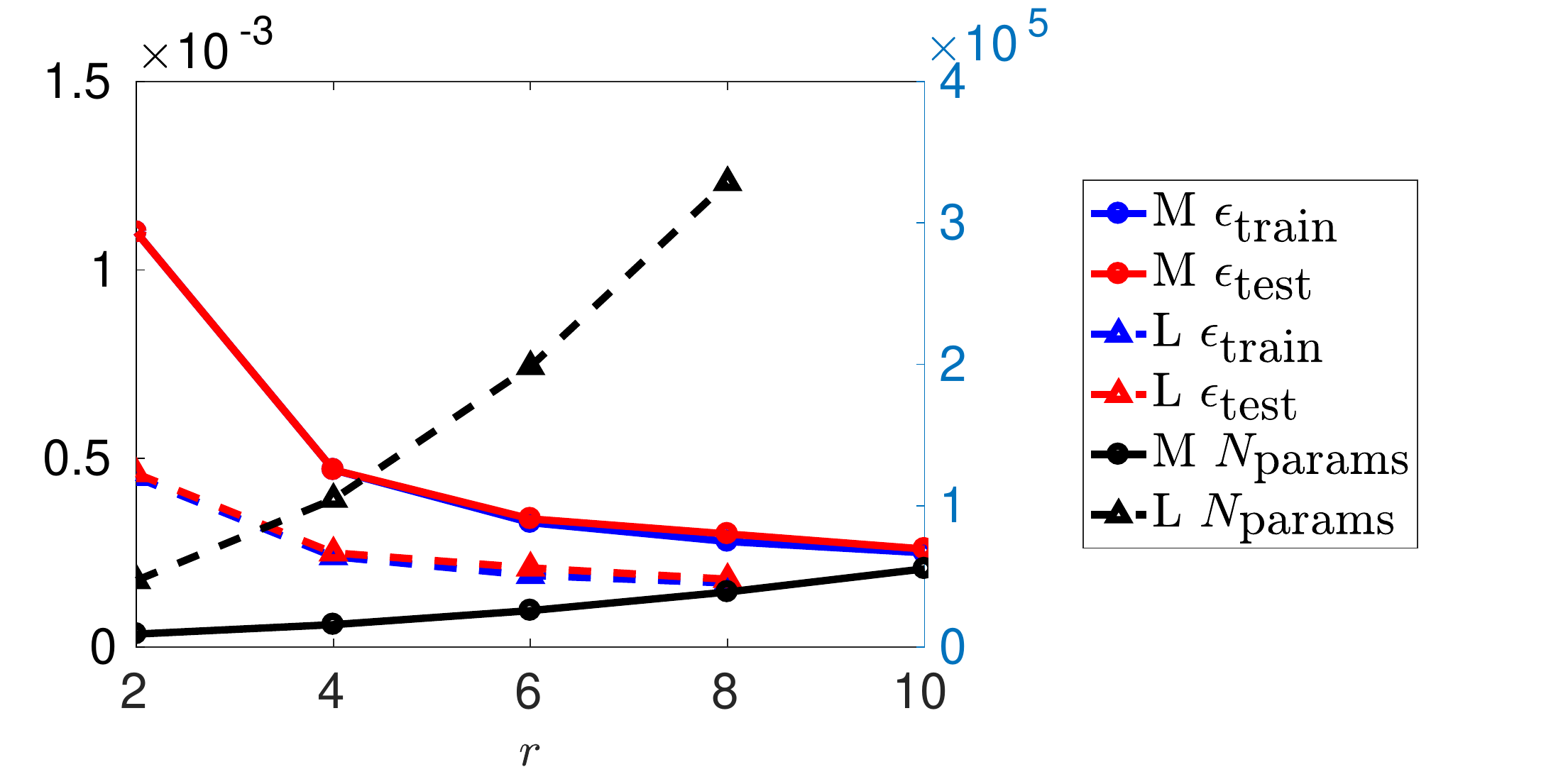}
    }
    \subfloat[\label{tab:rte1d_layers}$n_g=2$]{
    \includegraphics[width=0.3\textwidth, clip]{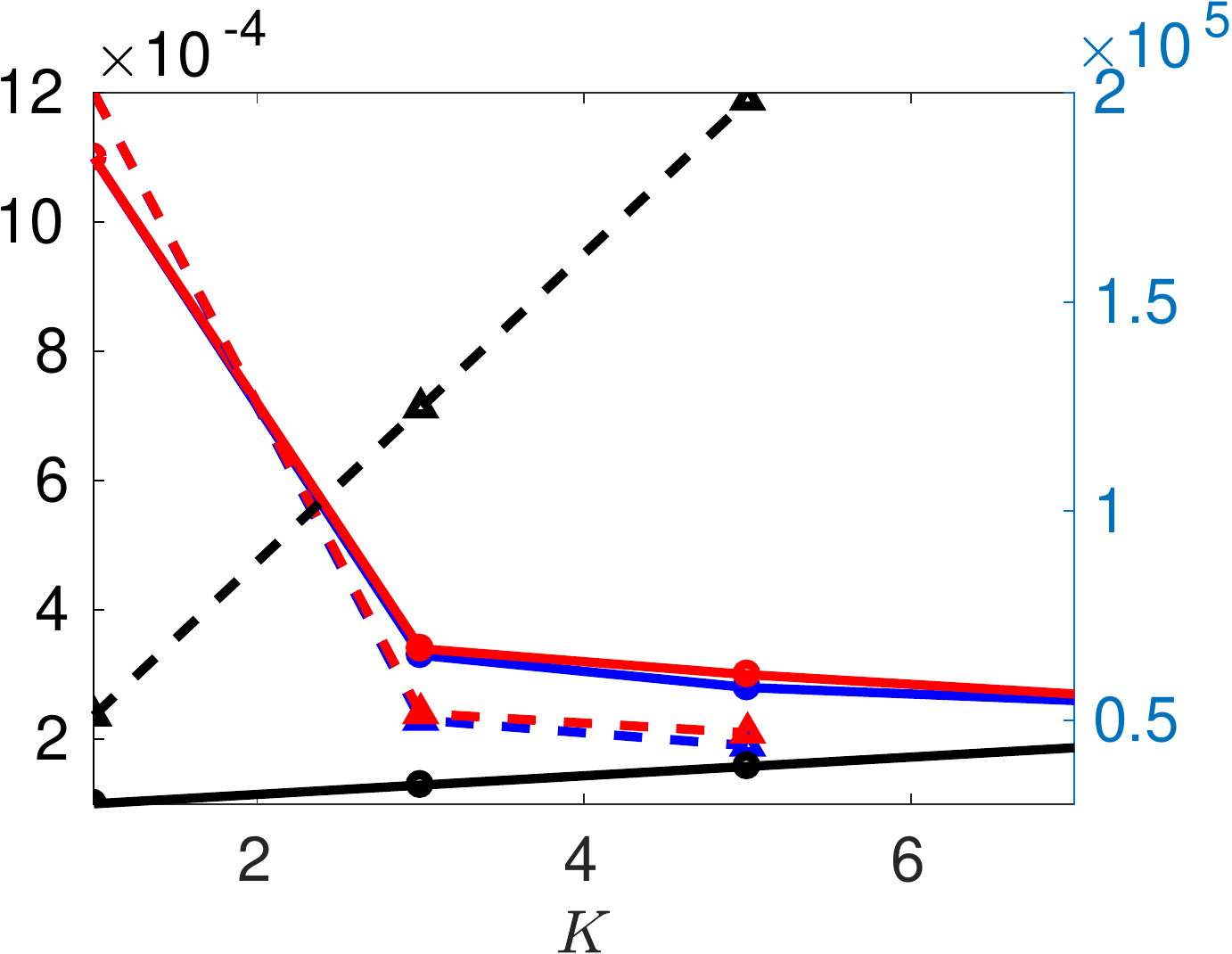}
    }
    \caption{\label{tab:rte1d}Relative error in approximating the density of RTE for 1D case for
    different number of channels $r$ and different number of $\CNNK$/$\LCK$ layers $K$ with
    {\Ntrainsample}=\Ntestsample $=20000$. ``M'' and ``L'' stands for MNN-$\cH^2$-Mix and
    MNN-$\cH^2$-LC, respectively. (b) the number of channel $r$ is 8 for MNN-$\cH^2$-Mix and is 6
    for MNN-$\cH^2$-LC.}
\end{figure}

\begin{figure}[htb]
    \centering
    \includegraphics[width=0.3\textwidth]{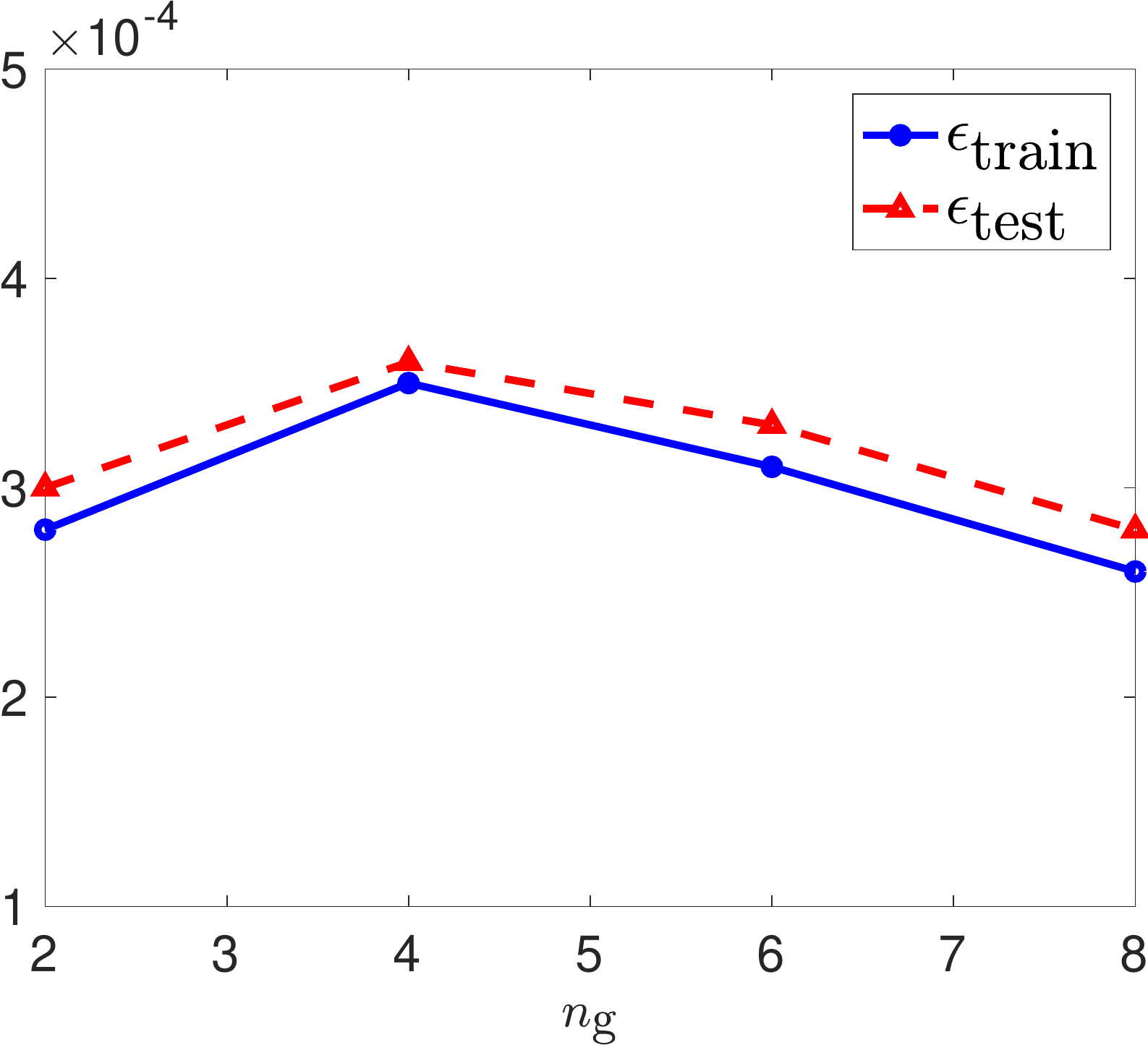}
    \caption{\label{tab:rte1d_ng}Relative error in approximating the density of RTE for 1D case for
    different number of Gaussians $n_g$ for MNN-$\cH^2$-Mix with $K=5$, $r=8$ and {\Ntrainsample} =
    \Ntestsample $=20000$.}
\end{figure}

Here we set $f(x)=1$, and $\mu_a(x)=\mu_t(x)-\mu_s(x)=0.2$, $x\in\Omega$, and the scattering
coefficient has the form
\begin{equation}
    \mu_s(x) = \sum_{i=1}^{n_g}\frac{\rho^{(i)}}{\sqrt{2\pi T}}\exp\left(
    -\frac{|x-c^{(i)}|^2}{2T} \right),
\end{equation}
where the parameters $\rho^{(i)}\sim \cU(0.1, 0.3)$, $c^{(i)}\sim
\cU(0.2,0.8)$, $i=1,\dots,n_g$ and $T\sim \cU(2,4)\times 10^{-3}$.
The numerical samples are generated by solving \eqref{eq:RTE_IE1d}.

Because the map $\mu_s\rightarrow u$ is not translation invariant, MNN-$\cH^2$ cannot be implemented
using CNNs as before. As discussed at the end of \cref{sec:invariant}, we can combine LC layers and
CNN layers together to reduce the number of parameters. The resulting neural network is denoted by
MNN-$\cH^2$-Mix. As a reference, we implement MNN-$\cH^2$ by LC network and it is denoted by
MNN-$\cH^2$-LC. Note that since both $\mu_s$ and $u$ are not periodic the periodic padding in
$\LCK$/$\CNNK$ should be replaced by zero padding.

The number of discretization points is $N=320$, and $L=6$, $m=5$.  We perform numerical experiments
to study the numerical behavior for different number of channels (\cref{tab:rte1d_channel}) and
different number of $\CNNK$/$\LCK$ layers $K$ (\cref{tab:rte1d_layers}).  For both MNN-$\cH^2$-Mix
and MNN-$\cH^2$-LC, as $r$ or $K$ increase, the errors first decrease and then stagnate. We use
$r=8$ and $K=5$ for MNN-$\cH^2$-Mix in the following.  For the same setup, the error of
MNN-$\cH^2$-LC is somewhat smaller and the number of parameters is quite larger than that of
MNN-$\cH^2$-Mix. Thus, MNN-$\cH^2$-Mix serves as a good balance between the number of parameters and
the accuracy.

\cref{tab:rte1d_ng} summarizes the results of MNN-$\cH^2$-Mix for different $n_g$ with $K=5$ and
$r=8$. Numerical results show that MNN-$\cH^2$-Mix is not sensitive to the complexity of the input.

\subsubsection{Two-dimensional case}
Here we set $f(x)=1$ and $\mu_a(x)=\mu_t(x)-\mu_s(x)=0.2$ for $x\in\Omega$. The scattering
coefficient takes the form
\begin{equation}
  \mu_s(x) = \sum_{i=1}^{2}\frac{\rho^{(i)}}{{2\pi T}}\exp\left(-\frac{|x-c^{(i)}|^2}{2T} \right),
\end{equation}
where $x=(x_1,x_2)$ and the parameters $\rho^{(i)}\sim \cU(0.01, 0.03)$,
$c^{(i)}\sim\cU(0.2,0.8)^2$, $i=1,2$ and $T\sim \cU(2,4)\times 10^{-3}$.  The numerical samples are
generated by solving \eqref{eq:RTE_IE}.

\begin{figure}[htb]
    \centering
    \includegraphics[width=0.3\textwidth]{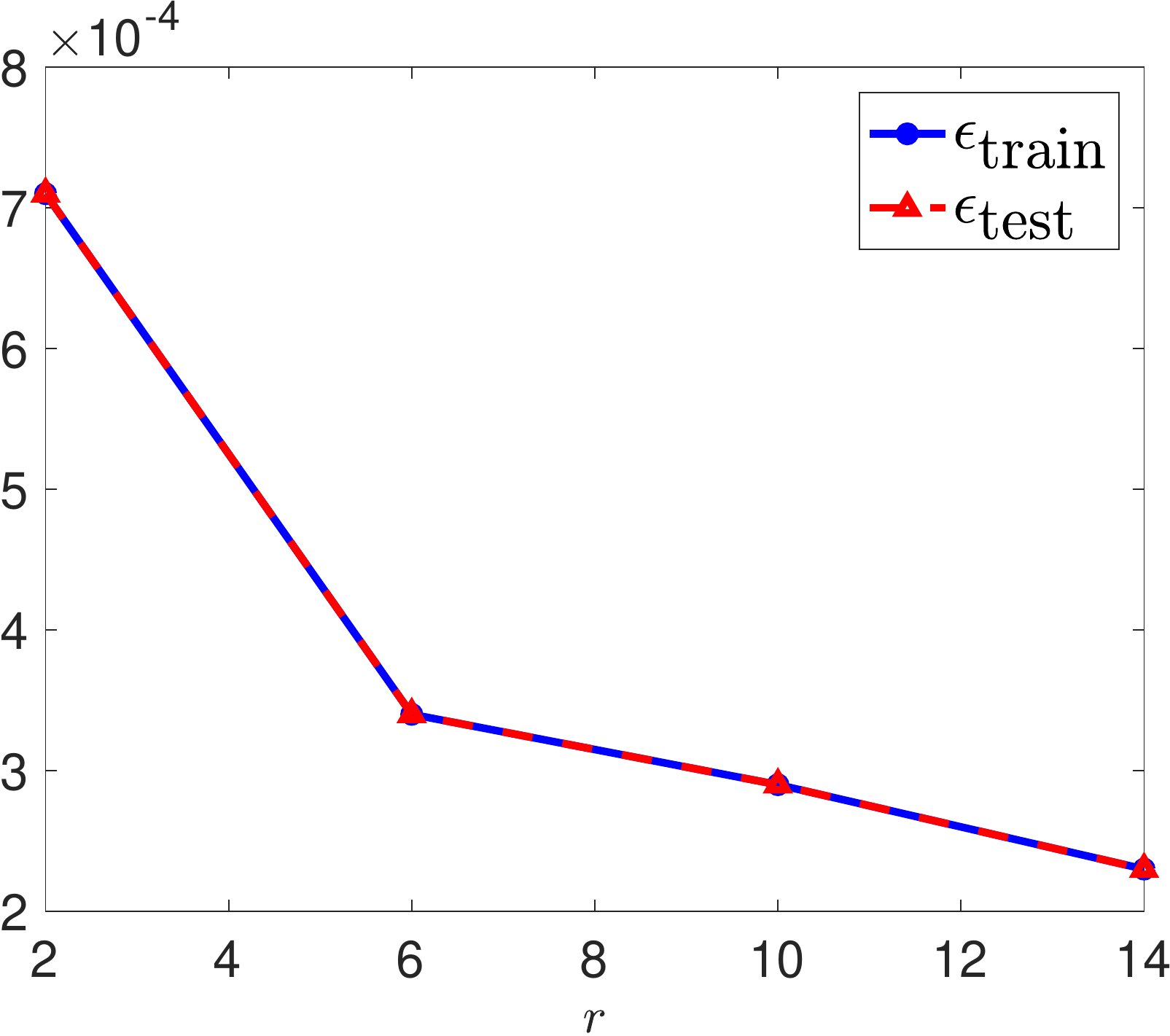}
    \caption{\label{tab:rte2d_al}Relative error in approximating the density of RTE for the 2D case for
    different number of channels for MNN-$\cH^2$-Mix with $K=5$ and {\Ntrainsample} = \Ntestsample
    $=20000$.}
\end{figure}

Because the map $\mu_s\rightarrow u$ is not translation invariant, we implement the MNN-$\cH^2$-Mix
architecture as the 1D case. Considering that the adjacent part takes a large number of parameters
for the 2D case, we implement the adjacent part by the $\CNNK$ layers. \cref{tab:rte2d_al} gathers the
results for different number of channels $r$. Note that, similar to the 1D case, there is no
overfitting for all the tests and the relative error decreases as $r$ increases.

\subsection{Kohn-Sham map}\label{sec:KSMap}

In the Kohn-Sham density functional theory~\cite{HohenbergKohn1964,KohnSham1965}, one needs to solve
the following nonlinear eigenvalue
equations  (spin degeneracy omitted):
\begin{equation}\label{eqn:KSDFT}
  \begin{split}
    & \left(-\frac{1}{2} \Delta + V(x)\right) \psi_{i}(x) =
    \varepsilon_{i} \psi_{i}(x), \,\, x \in \Omega=[-1,1)^d \\
    & \int_{\Omega} \psi_{i}(x) \psi_{j}(x) d x = \delta_{ij},
    \quad \rho(x) = \sum_{i=1}^{n_e} |\psi_i(x)|^2,
   \end{split}
\end{equation}
where $n_e$ is the number of electrons, $d$ is the spatial dimension,
and $\delta_{ij}$ stands for the
Kronecker delta. All eigenvalues $\{\varepsilon_{i}\}$ are real and ordered
non-decreasingly. The electron density $\rho(x)$ satisfies the constraint
\begin{equation}
  \rho(x)\ge 0,\quad \int_{\Omega} \rho(x) \dd x = n_{e}.
  \label{eqn:rho_constraint}
\end{equation}
In this subsection, we employ the multiscale neural networks to approximate the 
Kohn-Sham map
\begin{equation}
    \mathcal{F}_{\text{KS}}: V \to \rho.
\end{equation}
The potential function $V$ is given by
\begin{equation} \label{eqn:gaussian_wells}
    V(x) = -\sum_{i=1}^{n_e} \sum_{j\in \mathbb{Z}^d} \rho^{(i)} \exp\left(-\frac{(x - c^{(i)}-2j)^2}{2 \sigma^2}\right),
    \qquad  x \in [-1,1)^d,
\end{equation}
where $c^{(i)} \in [-1,1)^d$ and $\rho^{(i)}\in\cU(0.8,1.2)$. We set $\sigma = 0.05$ for 1D and $\sigma=0.2$
for the 2D case. 
The centers of the Gaussian wells $c^{(i)}$ are chosen randomly under the
constraint that $|c^{(i)}-c^{(j)}| > 2 \sigma$.
The Kohn-Sham map is discretized using a pseudo-spectral method \cite{Trefethen2000}, and solved by
a standard eigensolver.

\subsubsection{One-dimensional case}
For the one-dimensional case, we choose $N=320$, $L=7$ and $m=5$, and use the same datasets as in
\cite{fan2018mnn} to study the numerical behavior of MNN-$\cH^2$ for different $n_e$, $r$ and $K$.

\begin{figure}[htb]
    \centering
    \subfloat[\label{tab:ks_1d_channel} $K=6$, $n_e=2$]{
    \includegraphics[width=0.3\textwidth]{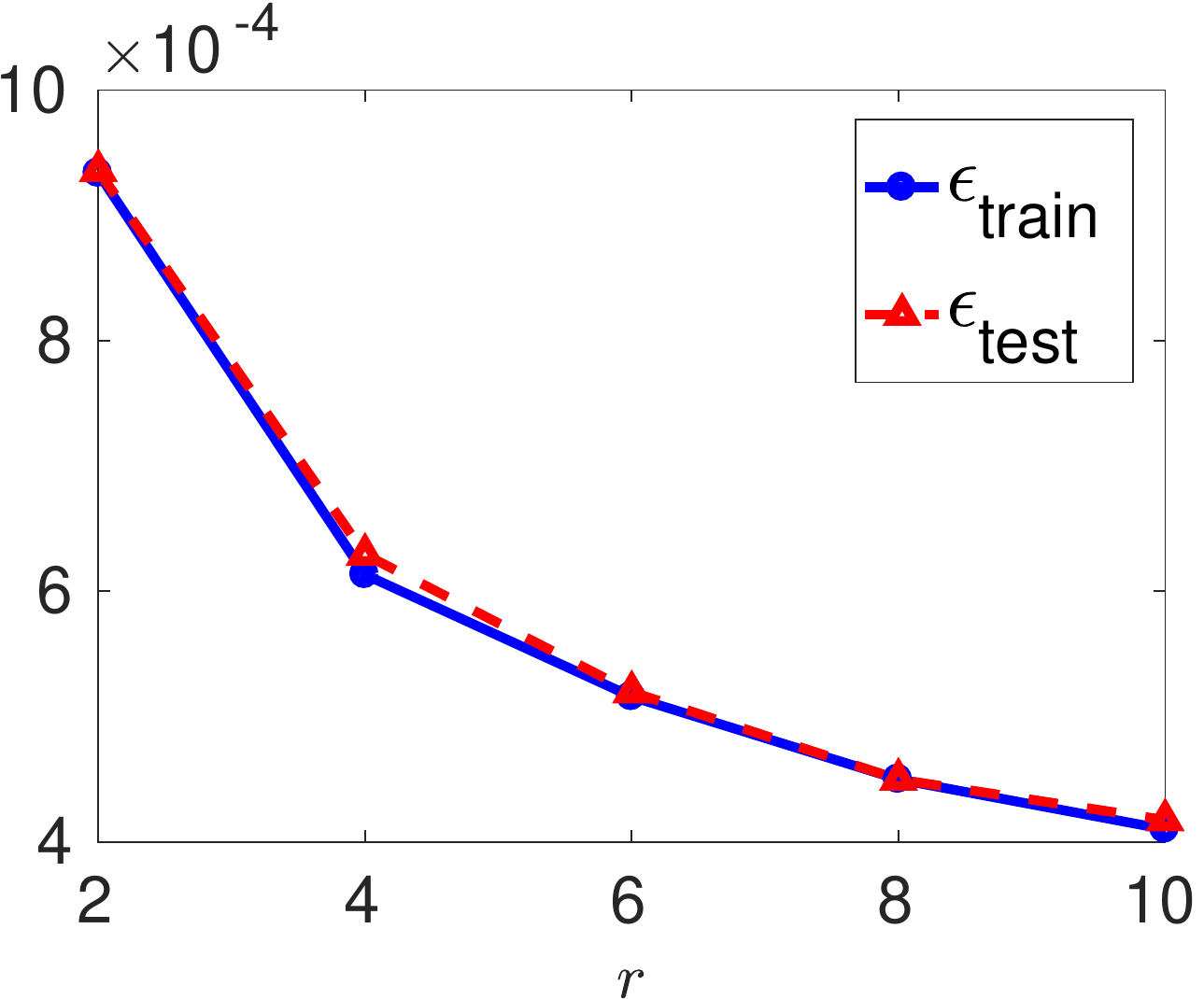}
    }
    \subfloat[\label{tab:ks_1d_K} $r=8$, $n_e=2$]{
    \includegraphics[width=0.3\textwidth]{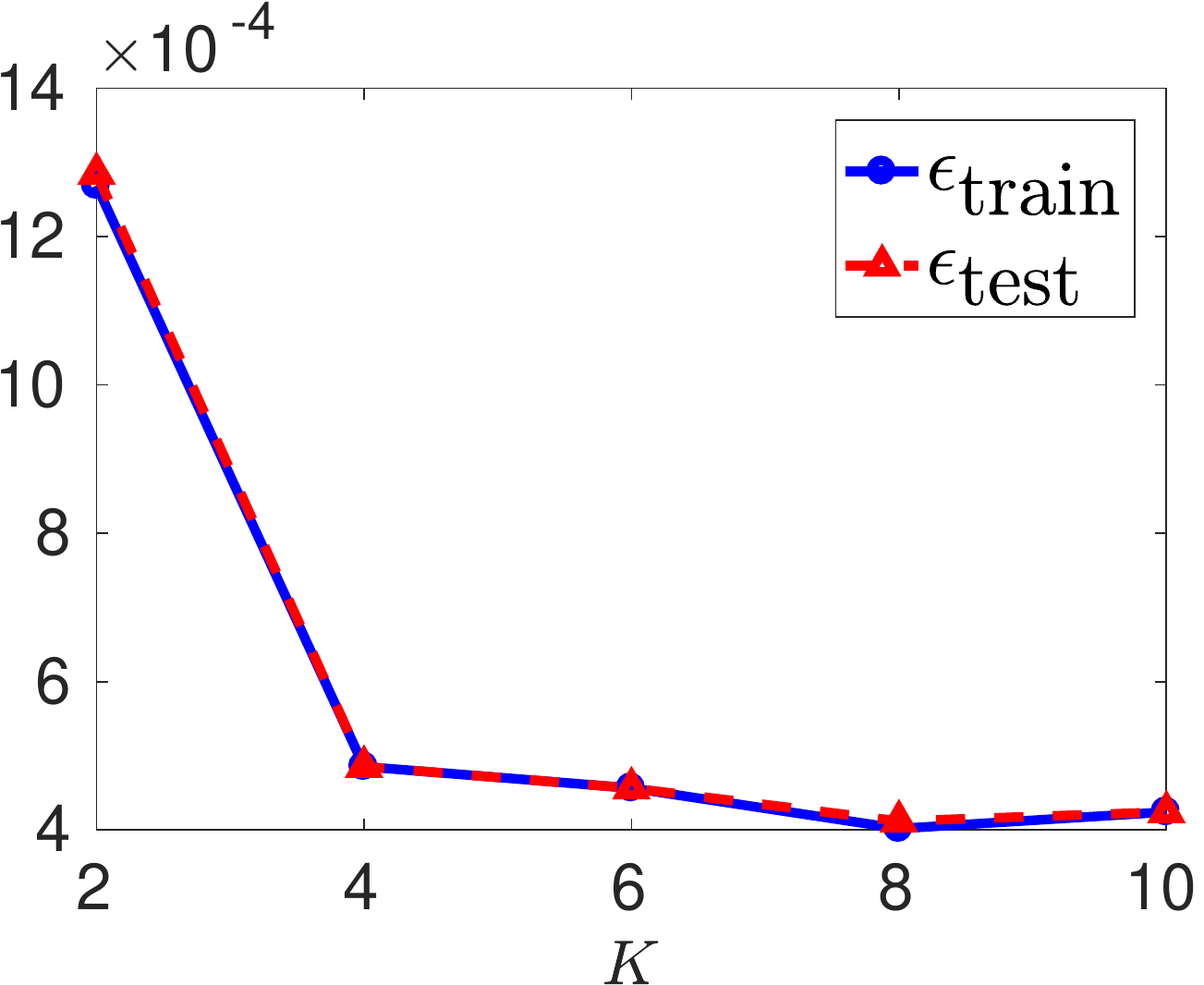}
    }
    \subfloat[\label{tab:ks_1d_ng} $K=6$, $r=5$]{
    \includegraphics[width=0.3\textwidth]{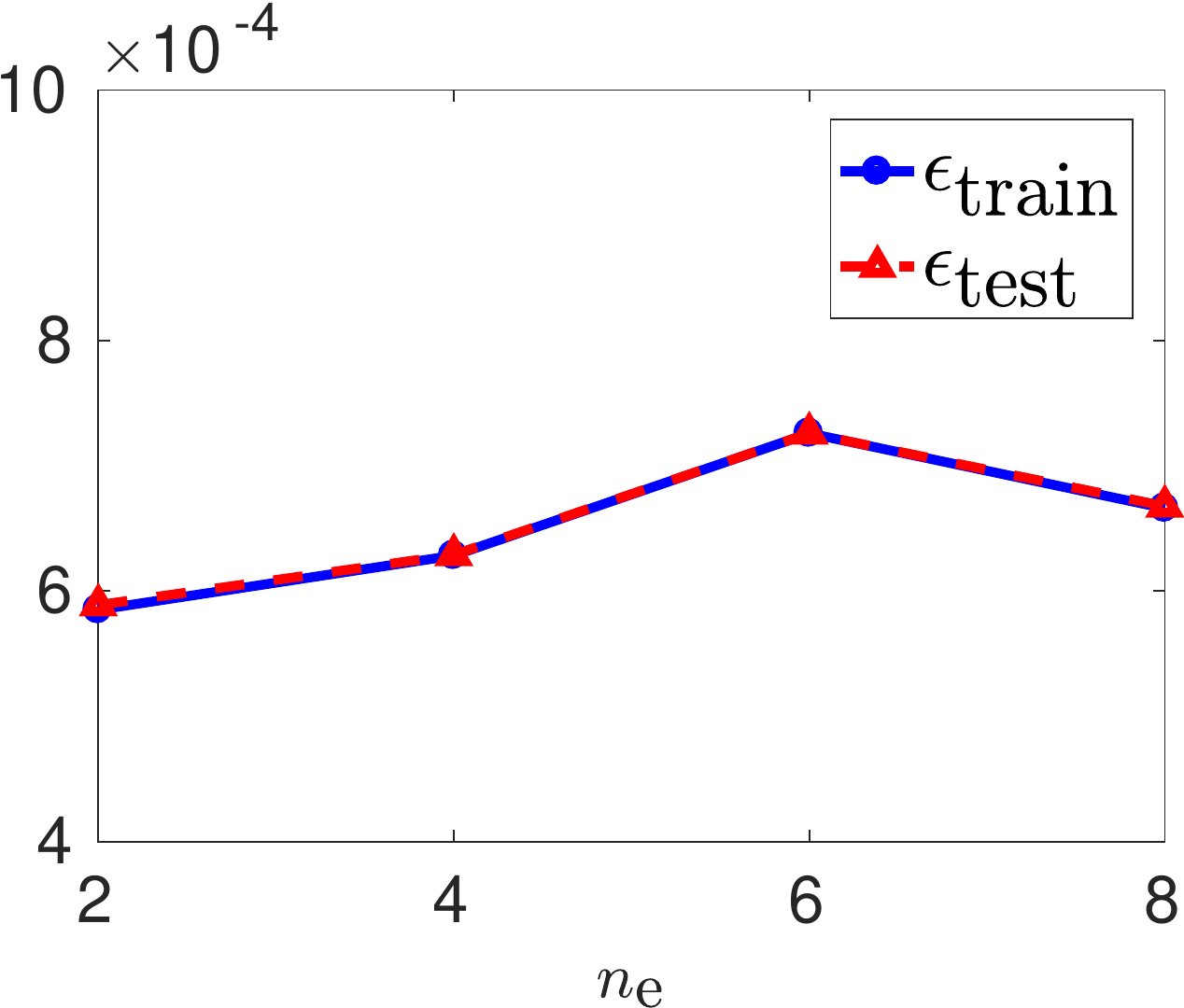}
    }
    \caption{ Relative error on the approximation of the Kohn-Sham map for different $r$, $K$, and $n_g$ {\Ntrainsample} =$16000$, and \Ntestsample =$4000$.}
\end{figure}
\begin{figure}[htb]
    \centering
    \subfloat[test error]{
    \includegraphics[width=0.3\textwidth,clip]{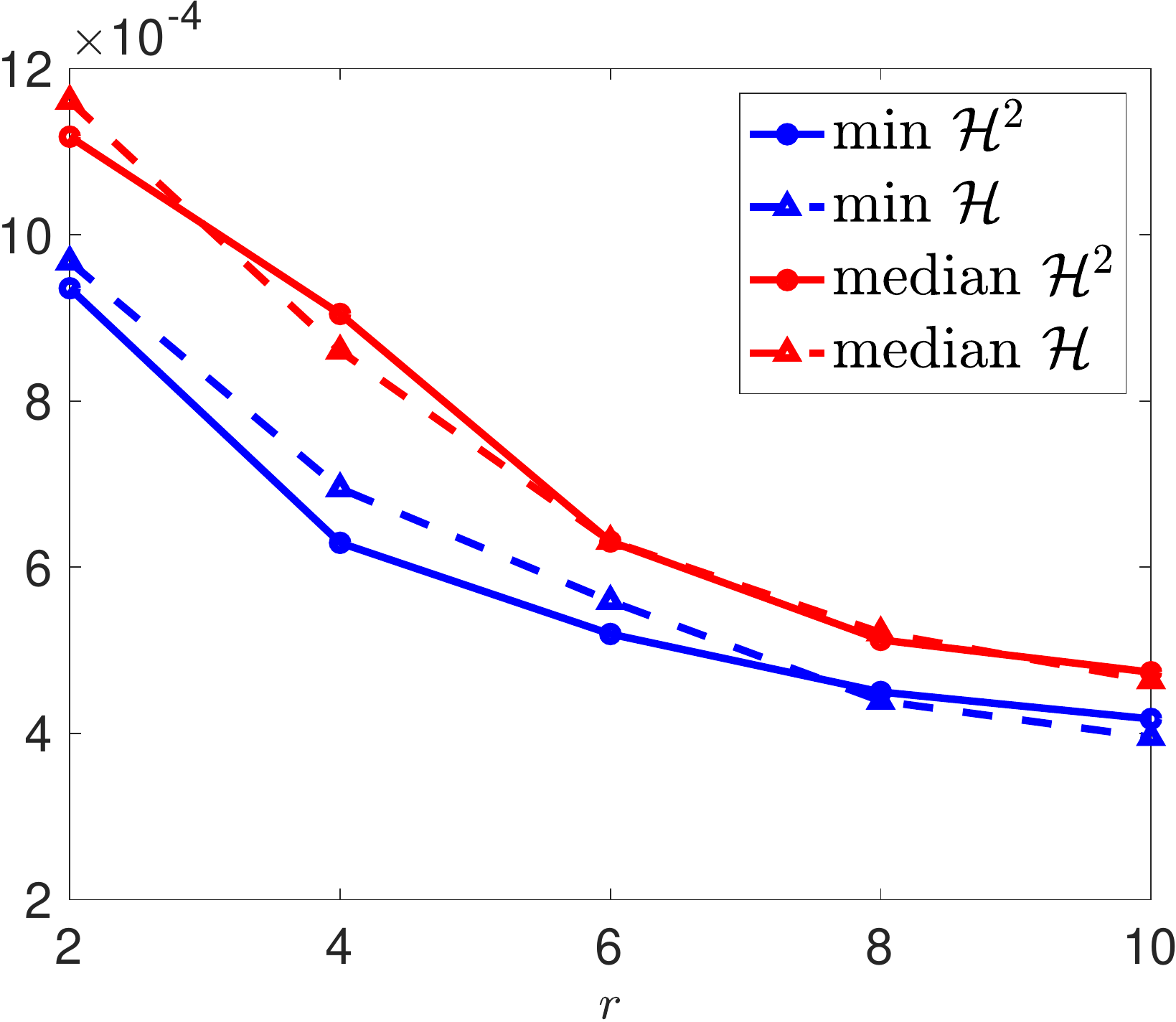}
    }\qquad
    \subfloat[\Nparams]{
    \includegraphics[width=0.3\textwidth,clip]{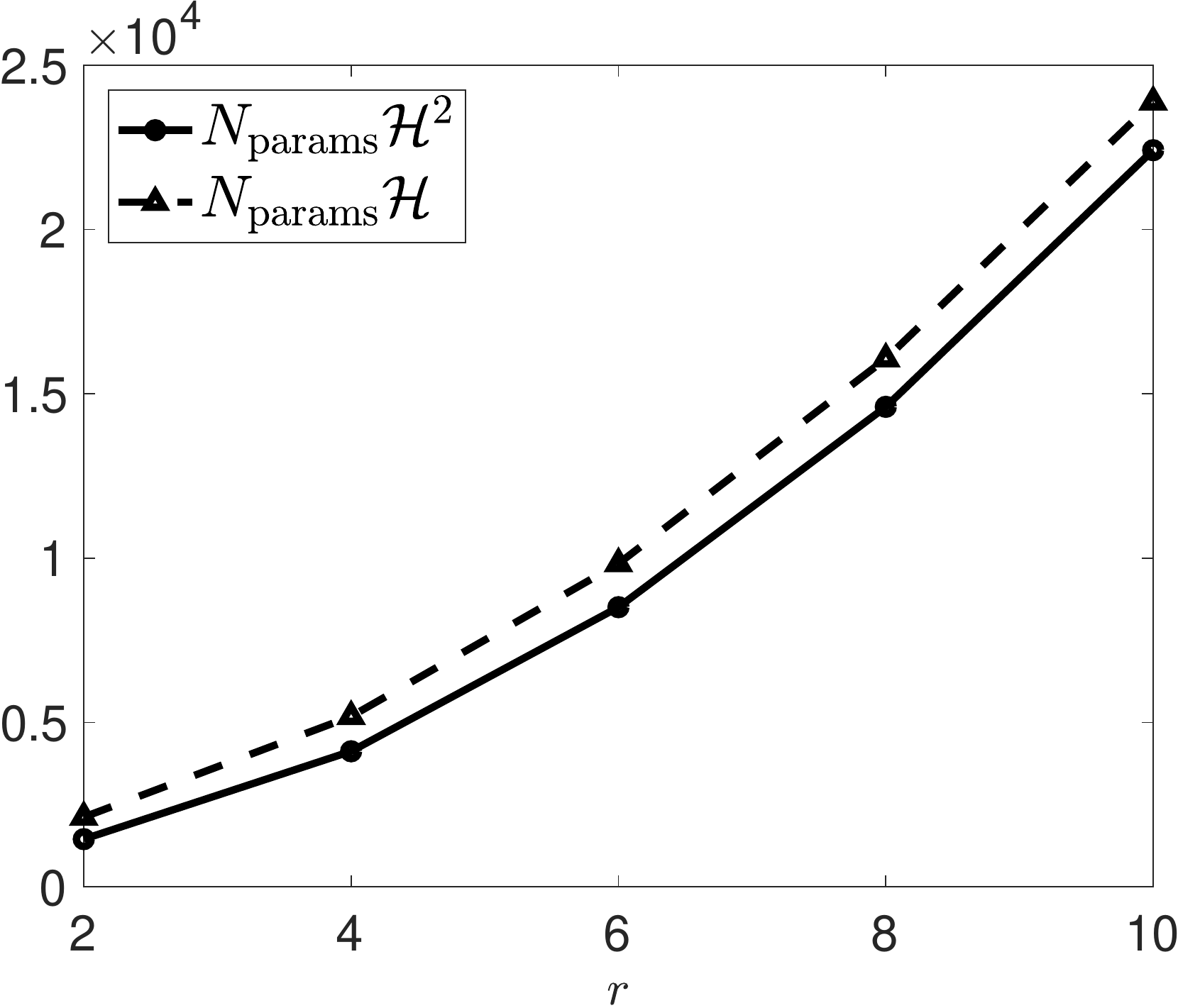}
    }
    \caption{\label{tab:ksmap1d_com}
    Numerical results of MNN-$\cH$ / MNN-$\cH^2$ for the minimum and median $\epsilon_{\mathrm{train}}$ for 1D
    Kohn-Sham map with random initial seed. The ``min'' and ``median'' stand for the test error
    corresponding to the minimum and median training data cases, respectively, and $H$ and $H^2$
    stand for MNN-$\cH$ and MNN-$\cH^2$, respectively. The setup of MNN-$\cH^2$ is $K=5$, $n_g=2$
    and {\Ntrainsample} = \Ntestsample $=5000$.}
\end{figure}

From \cref{tab:ksmap1d_com} we observe that both architectures, MNN-$\mathcal{H}^2$ and
MNN-$\mathcal{H}$, provide comparable results even as the MNN-$\mathcal{H}^2$ has fewer parameters
to fit. Both architectures show the same trends. As the number of channels, $r$, increases the error
decreases sharply, and then stagnates rapidly as shown in \cref{tab:ks_1d_channel}. On the other
hand, as the number of layers, $K$, increases the error decreases sharply, and then stagnates as $K$
becomes large as shown in \cref{tab:ks_1d_K}.
Finally, \cref{tab:ks_1d_ng} shows that the accuracy of MNN-$\mathcal{H}^2$ is relatively
insensitive to the number of wells. In addition, as shown before, we do not observe overfitting for
this example.

\subsubsection{Two-dimensional case}

\begin{figure}
    \centering
    \includegraphics[width=0.3\textwidth]{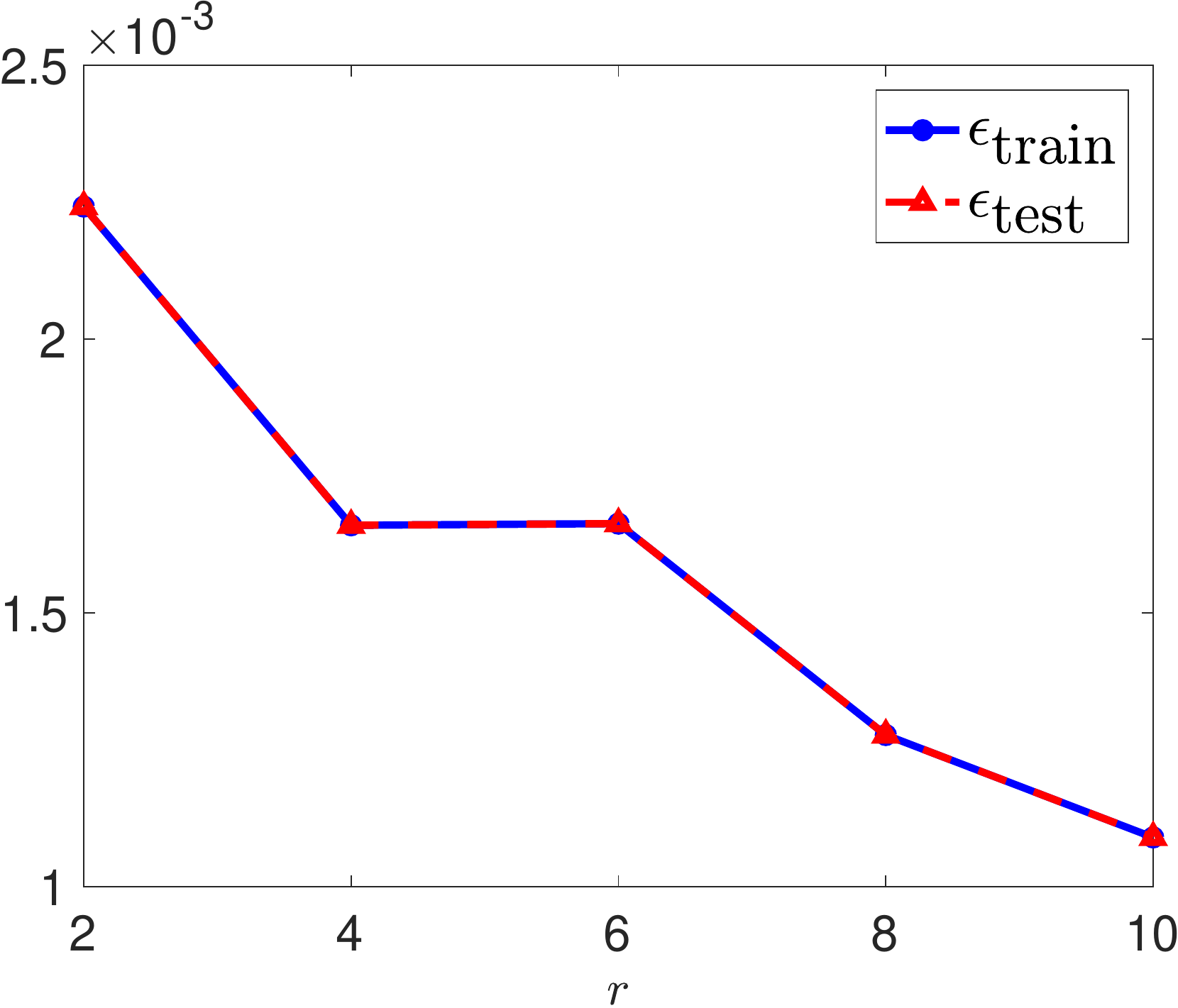}
    \caption{  Relative test error on the approximation of the 2D Kohn-Sham map for different 
    number of channels $r$, and {\Ntrainsample} = $16000$.}
    \label{tab:ks_2d_channel}
\end{figure}

The discretization is the standard extension to 2D using tensor products, using a $64\times 64$ grid.
We consider $n_e=2$ and follow the same number of training and test
samples as that in the $1$D case. We fixed $K = 6$, $L=4$ and $m=4$, and we trained both networks
for different number of channels, $r$. The results are displayed in \cref{tab:ks_2d_channel}, which
shows the same behavior as for the 1D case, comparable errors for both architectures with the error
decreasing as $r$ increases, with virtually no overfitting.

\section{Conclusion}\label{sec:conclusion}
In this paper, motivated by the fast multipole method (FMM) and $\cH^2$-matrices, we developed a
multiscale neural network architecture (MNN-$\cH^2$) to approximate nonlinear maps arising from
integral equations and partial differential equations. Using the framework of neural networks,
MNN-$\cH^2$ naturally generalizes $\cH^2$-matrices to the nonlinear setting.  Compared to the
multiscale neural network based on hierarchical matrices (MNN-$\cH$), the distinguishing feature of
MNN-$\cH^2$ is that the interpolation and restriction layers are represented using a set of nested
layers, which reduces the computational and storage cost for large systems. Numerical results
indicate that MNN-$\cH^2$ can effectively approximate complex nonlinear maps 
arising from the nonlinear Schr\"odinger equation, the steady-state radiative
transfer equation, and the Kohn-Sham density functional theory.
The MNN-$\cH^2$ architecture can be naturally extended. For instance, the $\LCR$ and $\LCI$ networks
can involve nonlinear activation functions and can be extended to networks with more than one layer.
The $\LCK$ network can also be altered to other network structures, such as the sum of two parallel
subnetworks or the ResNet architecture \cite{he2016deep}.


\section*{Acknowledgements}
The work of Y.F. and L.Y. is partially supported by the U.S. Department of Energy, Office of
Science, Office of Advanced Scientific Computing Research, Scientific Discovery through Advanced
Computing (SciDAC) program, the National Science Foundation under award DMS-1818449, 
and the GCP Research Credits Program from Google.
The work of J.F. is partially supported by ``la Caixa'' Fellowship, sponsored by the ``la Caixa''
Banking Foundation of Spain.
The work of L.L and L.Z. is partially supported by the Department of Energy under Grant No.
DE-SC0017867 and the CAMERA project.

\appendix

\section{Comparing MNN-$\cH^2$ with CNN}



In this appendix, by comparing MNN-$\cH^2$ with the classical convolutional neural networks (CNN),
we show that multiscale neural networks not only reduce the number of parameters, but also improve
the accuracy. Since the RTE example is not translation invariant, we perform the comparison using
NLSE and Kohn-Sham map.

\paragraph{NLSE with inhomogeneous background potential}
\begin{figure}[htb]
  \centering \includegraphics[width=0.3\textwidth,clip]{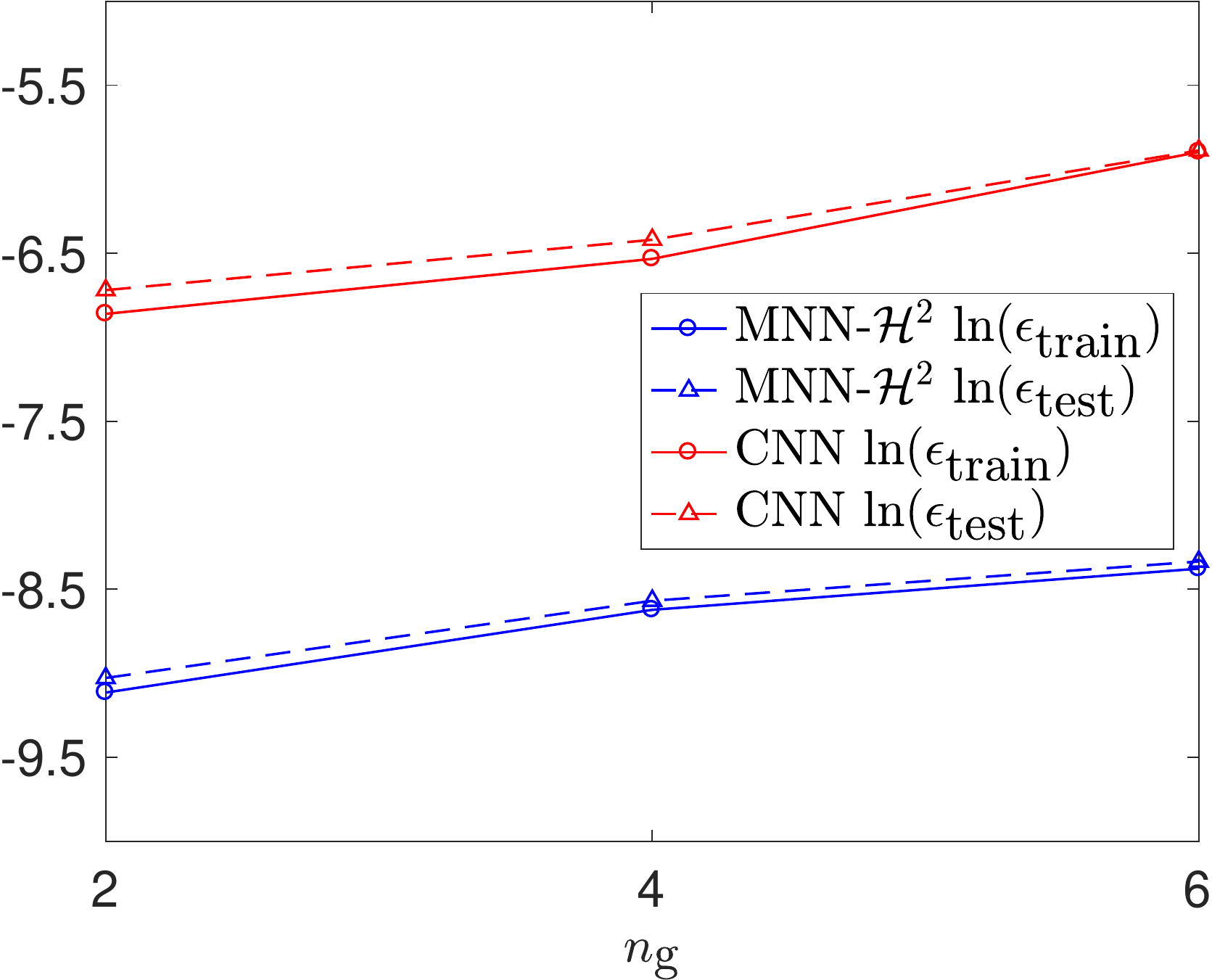} \caption{\label{tab:nlse1d_CNN}
  The training and test errors of MNN-$\cH^2$ with 7209 parameters ($r=6$ and $K=5$) and CNN with
  38161 parameters (15 layers, 10 channels and window size to be $25$) for the one dimensional
    NLSE.}
\end{figure}
Here we study the one-dimensional NLSE using the setup from \cref{sec:nlse1d} for different number
of Gaussians in the potential $V$ \eqref{eq:nlse}. The training and test errors for MNN-$\cH^2$ and
CNN are presented in \cref{tab:nlse1d_CNN}. The channel number, layer number, and window size of CNN
are optimally tuned based on the training error. The figure demonstrates that MNN-$\cH^2$ has fewer
parameters and gives a better approximation to the NLSE.

\begin{figure}[htb]
  \centering
  \includegraphics[width=0.3\textwidth,clip]{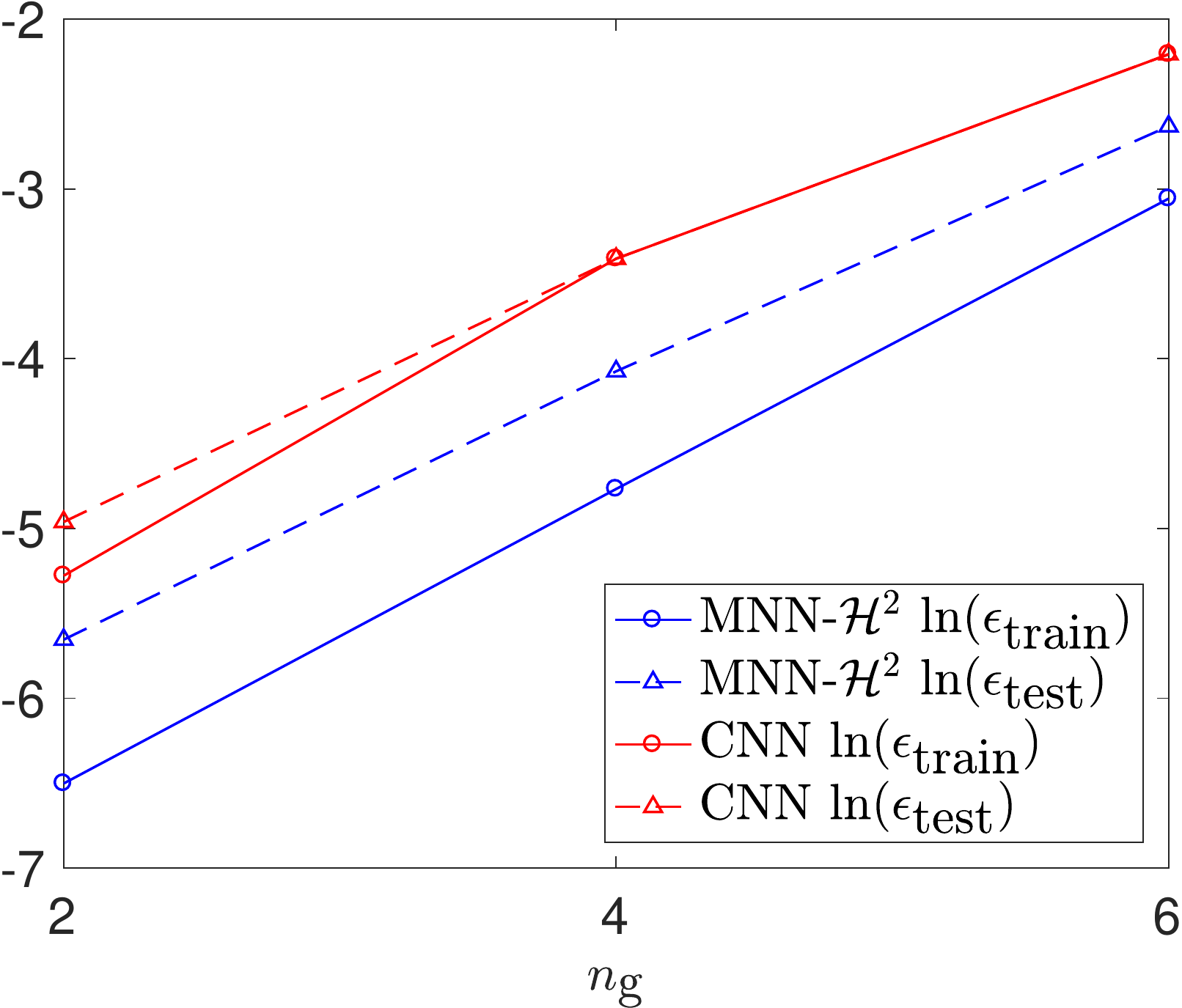}
  \caption{\label{tab:ksmap1d_CNN}
  The training and test errors of MNN-$\cH^2$ with 18985 parameters ($r=10$ and $K=5$) and CNN with 
  25999 parameters (10 layers, 10 channels and window size to be $13$) for the one dimensional
  Kohn-Sham map.}
\end{figure}

\paragraph{Kohn-Sham map}
For the Kohn-Sham map, we consider the one-dimensional setting in \eqref{eqn:gaussian_wells} with
varying number of Gaussian wells. The width of the Gaussian well is set to be $6$. In this case, the
average size of the band gap is $0.01$, and the electron density at point $x$ can depend
sensitively on the value of the potential at a point $y$ that is far away. 
\cref{tab:ksmap1d_CNN} presents the training and test errors of MNN-$\cH^2$ and CNN, where 
MNN-$\cH^2$ outperforms a regular CNN with a comparable number of parameters.

%
%

\bibliographystyle{abbrv}
\bibliography{nn}
\end{document}